\documentclass[reqno,11pt]{amsart}
\usepackage{graphicx}
\usepackage{amssymb}
\usepackage{amsmath}
\usepackage{amsfonts}
\usepackage{amsthm}
\usepackage{color}
\usepackage{hyperref}
\usepackage[mathscr]{eucal}
\usepackage{enumerate}
\usepackage{times}

\numberwithin{equation}{section}
\allowdisplaybreaks[1]

\newtheorem{Def}{Definition}[section]
\newtheorem{Thm}[Def]{Theorem}
\newtheorem{Prop}[Def]{Proposition}
\newtheorem{Lemma}[Def]{Lemma}

\newcommand{\beq}{\begin{equation}}
\newcommand{\eeq}{\end{equation}}
\newcommand{\Proof}{\begin{proof}}
\newcommand{\QED}{\end{proof} \noindent}

\newcommand{\mm}{\hspace{-.08cm}\cdot \hspace{-.08cm}}

\newcommand{\M}{\mathcal{M}}

\newcommand{\R}{\mathbb{R}}

\newcommand{\VM}{\mathcal{V}\hspace{-.05cm}\mathcal{M}}

\newcommand{\A}{\mathcal{A}}
\newcommand{\Ati}{\tilde{\mathcal{A}}}
\newcommand{\Aop}{\mathcal{A}_\textbf{b}}
\newcommand{\X}{X}

\newcommand{\cutoff}{\chi_{_\delta}}

\title[Optimal Regularity for Connections on Vector Bundles]{On the Optimal Regularity Implied by the Assumptions of Geometry II: Connections on Vector Bundles}

\author[M.\ Reintjes]{Moritz Reintjes}
\address{Department of Mathematics\\ City University of Hong Kong \\ Kowloon \\ Hong Kong}
\email{moritzreintjes@gmail.com}

\author[B.\ Temple]{Blake Temple \\ \\ January 28, 2024}
\address{Department of Mathematics\\ University of California\\ Davis, CA 95616\\ USA}
\email{temple@math.ucdavis.edu}

\begin{document} 

\begin{abstract}
We extend authors' prior results on optimal regularity and Uhlenbeck compactness for affine connections to general connections on vector bundles.  This is accomplished by deriving a vector bundle version of the RT-equations, and establishing a new existence theory for these equations.  These new RT-equations,  non-invariant elliptic equations, provide the gauge transformations which transform the fibre component of a non-optimal connection to {\it optimal regularity}, i.e., the connection is one derivative more regular than its curvature in $L^p$. The existence theory handles curvature regularity all the way down to, but not including, $L^1$. Taken together with the affine case, our results extend optimal regularity of Kazden-DeTurck and the compactness theorem of Uhlenbeck, applicable to Riemannian geometry and compact gauge groups, to general connections on vector bundles over {\it non-Riemannian} manifolds, allowing for both compact and non-compact gauge groups.   In particular, this extends optimal regularity and Uhlenbeck compactness to Yang-Mills connections on vector bundles over Lorentzian manifolds as base space, the setting of General Relativity.
\end{abstract}


\maketitle 

\setcounter{tocdepth}{1}
\small
\tableofcontents
\normalsize

\section{Introduction}

In this paper we extend our optimal regularity and Uhlenbeck compactness results for affine connections,  (i.e., connections on the tangent bundle of a manifold $\mathcal{M}$ \cite{ReintjesTemple_ell1, ReintjesTemple_ell2, ReintjesTemple_ell3, ReintjesTemple_Uhl1, ReintjesTemple_ell6}), to general connections on vector bundles.  This is accomplished by deriving a vector bundle version of the {\it RT-equations}\footnote{The \emph{Regularity Transformation} or \emph{Reintjes-Temple} equations \cite{ReintjesTemple_ell1,ReintjesTemple_ell3}.} introduced in \cite{ReintjesTemple_ell1}, and by developing a new existence theory for these equations.  The RT-equations are non-invariant elliptic equations whose solutions provide the gauge and coordinate transformations sufficient to transform a non-optimal connection to optimal regularity. Connections define parallel translation of vectors on the base manifold and in the fibre, and are the starting point of geometry.  A connection is of {\it optimal regularity} in a given coordinate system if its components are one derivative more regular than the components of its curvature in $L^p$.

Our theorems apply to connections on an $n$-dimensional manifold $\mathcal{M}$, whose components are assumed given in a single, global or local, coordinate chart defined on an open set $\Omega\subset\mathcal{M}$.  In this case the connection on the fibre essentially decouples from the connection on the base manifold $\mathcal{M}$, making optimal regularity on the fibre component of a connection independent of the affine component of the connection acting on the base manifold.  In \cite{ReintjesTemple_Uhl1} we established optimal regularity and Uhlenbeck compactness\footnote{By Uhlenbeck compactness we mean weak $W^{1,p}$ compactness derived from uniform $L^p$ bounds on the curvature and the connection alone, without the need to bound all connection derivatives.}  for the affine part of a connection, so to extend these results to the fibres, it suffices to assume, and there is no loss of generality in assuming, the base manifold is Euclidean space, and its affine part is trivial, c.f. \cite{Uhlenbeck}.   With this, we prove  {\it global} optimal regularity and Uhlenbeck compactness for a connection on the fibre whose components and curvature components are bounded in a given coordinate neighborhood $\Omega$, assuming a {\it smallness} condition only on the connection; and we obtain the same result {\it locally} without the smallness assumption.

Specifically, we prove every non-optimal connection $\A$ can be gauge transformed to optimal regularity $\Aop \in W^{1,\hat{p}}$ {\it globally}, i.e., in every compactly contained subset of $\Omega$, under the assumptions that the exterior derivative $d\A$ lies in $L^{\bar{p}}$, $\bar{p} \in (0,\infty)$, and that the connection components of $\A$ are {\it sufficiently small} in $L^{2p}(\Omega)$, $p\in (n/2,\infty]$. This places the curvature 
\beq \label{curvature}
\mathcal{F} \equiv d\A + \A \wedge \A
\eeq 
in $L^{\hat{p}}(\Omega)$ for $\hat{p} \equiv \min\{p,\bar{p}\}$, so $\Aop$ is one derivative more regular than the curvature. We obtain the same result {\it locally}, i.e., in a neighborhood of every point $q\in\Omega$, without the smallness assumption. Note, we assume connections in $L^{2p}$ so the wedge product in \eqref{curvature} is in $L^p$ by H\"older's inequality.   Our results also apply when $\bar{p}=\infty$, $p <\infty$; and when $\bar{p}=p=\infty$ this yields connection regularity arbitrarily close to optimal, since $L^\infty$ is contained in $L^{\bar{p}}$ on bounded sets for any $\bar{p}<\infty$. Our local regularity theorem does not directly imply global regularity because local maps to optimal regularity will not in general preserve the regularity of the transition maps defined on the overlap of local neighborhoods, and hence will not in general preserve the regularity of the atlas of maps which define the manifold.

Optimal regularity supplies an extra derivative bound for the connection, and this implies compactness for a sequence of connections, (Uhlenbeck compactness), so long as each element of the sequence meets the same bound on a uniform neighborhood.  Our optimal regularity theorem establishes this for the global problem, and for the local problem when $p=\infty$.  That is, we deduce Uhlenbeck compactness from optimal regularity on the fibre, globally, (under the smallness assumption on the connection), when $p\in (n/2,\infty]$, $\bar{p}\in(1,\infty)$, and locally (no smallness assumption) when  $p=\infty$, $\bar{p}\in(1,\infty)$. This extends \cite{ReintjesTemple_Uhl1} to curvatures in $L^{\hat{p}}$, $\hat{p} > 1$, and incorporates both compact and non-compact gauge groups. For comparison, Uhlenbeck \cite{Uhlenbeck} assumed compact gauge groups, and addressed connections in $W^{1,p}$, (i.e., already of optimal regularity),  with invariant curvature bounds in $L^p$, $p\in [n/2,\infty]$, where the curvature bounds are expressed in terms of an Riemannian metric of the base manifold.                 Since our results here are restricted to the fibres of a vector bundle, they apply also when the affine part of the connection is non-trivial. Thus putting the results here together with our results for affine connections in \cite{ReintjesTemple_Uhl1}, extends of Uhlenbeck compactness and Kazdan-DeTurck's optimal regularity,\footnote{\label{footnote_Uhl} Uhlenbeck's compactness theorem \cite{Uhlenbeck} is based on invariant curvature bounds inherent to positive definite (Riemannian) metrics. In this case Kazden-DeTurck would provide optimal regularity and compactness for the affine part of the connection when it is Riemannian.  It appears to the authors that a non-invariant version of \cite{Uhlenbeck} would apply here to the fibre component of a general connection by taking the base manifold to be Euclidean space, but optimal regularity and Uhlenbeck compactness for the affine part requires the RT-equations when the base manifold is not Riemannian.  The RT-equations provide new results, and unify the treatment of optimal regularity and Uhlenbeck compactness in the fibre and affine components of a connection.}                     from connections on vector bundles over Riemannian base manifolds, to connections on vector bundles over {\it non-Riemannian} base manifolds. This provides a unified framework incorporating Yang-Mills connections, Lorentzian manifolds and General Relativity.\footnote{The Yang-Mills equations in Lorentzian geometry have been extensively studied, (c.f. \cite{ChoquetChristodoulou, EardleyMoncrief, Finster1, Finster2, Friedrich, KlainermanTataru, KriegerTataru, SmollerWassermannYau, Stuart} and references therein), but to our knowledge none of these results, including the local and global existence results in \cite{ChoquetChristodoulou, EardleyMoncrief, Friedrich, KlainermanTataru, KriegerTataru}, address nor establish optimal regularity for Yang-Mills connections. Note that the notion of ``optimal regularity'' in \cite{KlainermanTataru} is different from ours, referring to well-posedness of the Cauchy problem close to critical Sobolev exponents.}         Our methods, essentially different from \cite{Uhlenbeck}, handle compact and non-compact gauge groups, and allow for curvature regularities as low as $L^{\hat{p}}$, $\hat{p}>1,$ on the fibre. To develop the ideas in a concrete and relevant setting, we address here the group $SO(r,s)$ acting on $N$-dimensional fibres, ($r,s\geq0$, $r+s=N$). 

Our prior results in \cite{ReintjesTemple_Uhl1} are based on the discovery of the RT-equations for affine connections, a Poisson-type nonlinear elliptic system of partial differential equations (PDE's) whose solutions determine coordinate transformations which lift a non-optimal affine connection $\Gamma$ to optimal regularity, one derivative more regular than its curvature ${\rm Riem}(\Gamma)$ in $L^p$. Specifically, in  \cite{ReintjesTemple_Uhl1}, authors prove existence for the RT-equations at low regularity $\Gamma \in L^{2p}$, $d\Gamma\in L^p$, $p>\max\{n/2,2\}$, implying the existence of coordinate transformations which (locally) lift the regularity of the connection on an $n$-dimensional manifold to optimal regularity, $\Gamma\in W^{1,p}$.    This result is general enough to resolve the authors'  original problem of establishing the optimal regularity of GR shock waves, \cite{GroahTemple, Israel, Reintjes, ReintjesTemple_first, ReintjesTemple_geo, SmollerTemple}, and it extends theorems of Kazdan-DeTurck \cite{DeTurckKazdan} and Uhlenbeck \cite{Uhlenbeck} for connections over {\it positive} definite Riemannian base manifolds, to the case of affine connections on arbitrary manifolds---but our prior work on affine connections is not general enough to handle the case of connections $(\Gamma,\A)$ on {\it vector bundles}, the setting of Uhlenbeck's work in \cite{Uhlenbeck}. The numerous applications of Uhlenbeck compactness in the positive definite case \cite{Donaldson,Taubes,Wehrheim} naturally begged the question as to whether our theory of the RT-equations, together with the optimal regularity and Uhlenbeck compactness they imply, extend from the case of affine connections, to connections $(\Gamma,\A)$ on vector bundles over arbitrary base manifolds.  It is this extension which is accomplished in the present paper.

Putting our results for the affine part  $\Gamma$ of a connection together with the results established here for the fibre part  $\A$ of the connection,  we obtain optimal regularity for general connections $(\Gamma,\A)$ on vector bundles under the assumptions $(\Gamma,\A) \in L^{2p}$ with $d\Gamma\in L^{p}$, $d\A\in L^{\bar{p}}$ for $p\in (n/2,\infty]$, $\bar{p}\in(1,\infty)$, $n\geq 2$. This places the curvature in $L^{\hat{p}}$ for $\hat{p}=\min\{p,\bar{p}\}$. Again $L^p$-norms apply to components in a single coordinate chart defined on $\Omega\subset\mathcal{M}$.\footnote{We expect our results extend to complex vector bundles with Lie groups $U(r,s)$ and $SU(r,s)$, the Lie groups important in Physics, and to more general matrix Lie groups, by adapting our analysis to the modification of the RT-equations entailed by different choices of gauge group.}  Our results for the vector bundle version of the RT-equations together with the affine results in \cite{ReintjesTemple_Uhl1},  then imply coordinate and gauge transformations which lift the regularity of non-optimal connections up to optimal regularity $W^{1,\hat{p}}$, one derivative above the curvature in $L^{\hat{p}}$.  Again, by the extra derivative, this extends Uhlenbeck compactness to the general setting of connections $(\Gamma,\A)$ on vector bundles over arbitrary base manifolds with affine connection $\Gamma$. 
 
The RT-equations only involve $(\Gamma,\A)$ and the uncoupled leading order part $(d\A,d\Gamma)$ of the curvature, and decouple in $\Gamma$ and $\A$, even though the curvature does couple $\Gamma$ and $\A$ through the lower order commutator term. The vector bundle version of the RT-equations derived here is given by
\begin{eqnarray}
\Delta \Ati &=& \delta d \A -\delta\big( dU^{-1}\wedge dU \big),   \label{RT_1_intro} \\
\Delta U &=&   U \delta\A   -  (U^T \eta )^{-1}  \langle dU^T ; \eta dU\rangle , \label{RT_2_intro}
\end{eqnarray}
where $U$ is the sought after regularizing gauge transformation, and $\Ati$ is an auxiliary matrix valued $1$-form from which the connection of optimal regularity can be derived, c.f. \eqref{RT_1} - \eqref{RT_2} below.  We prove an existence theorem based on elliptic regularity in $L^p$-spaces for \eqref{RT_1_intro} - \eqref{RT_2_intro}, analogous to our theory in \cite{ReintjesTemple_ell3,ReintjesTemple_Uhl1}. This establishes the existence of global and local gauge transformations which lift connections to optimal regularity in transformed components $\Aop \in W^{1,\hat{p}}$, as stated above. Unlike the RT-equations in the affine case,  equation \eqref{RT_2_intro} for $U$ decouples from $d\A$. This allows us to extend the optimal regularity result to curvatures in $L^{\hat{p}}$ for $\hat{p} \in (1,\infty)$ so long as we require that $\A$ stays in $L^{2p}$ for $p\in (n/2,\infty]$, c.f. Theorem \ref{Thm_opt} below.  The RT-equations in the vector bundle version \eqref{RT_1_intro} - \eqref{RT_2_intro} are actually simpler than the affine version in the sense that they do not couple to a first order Cauchy-Riemann type equation, required in \cite{ReintjesTemple_Uhl1} to impose integrability of Jacobians to coordinates.\footnote{The tangent bundle is a vector bundle with gauge group $GL(n)$, but with the additional (more complicated) constraint that the gauge transformations (Jacobians) are integrable to coordinates.}    So the reader trying to understand the theory of the RT-equations in detail might do well to consider the vector bundle case here first.  

Uhlenbeck's pioneering theorem in \cite{Uhlenbeck} addresses the compactness of sequences of connections $(\Gamma,\A_i)$ on vector bundles with compact gauge groups, where $\Gamma$ is the affine connection of a fixed Riemannian metric on the base manifold, (c.f. footnote \ref{footnote_Uhl}).   Uhlenbeck's theorem assumes optimal regularity for  $(\Gamma,\A_i)$ at the start, and then derives compactness of $\A_i$ from a uniform bound on the curvature, formulated in terms of invariant norms associated with the Riemannian structure of the base manifold. Our original purpose here was to extend the local version of Uhlenbeck compactness to sequences of  {\it non-optimal} connections $(\Gamma_i,\A_i)$ on vector bundles where $\Gamma_i$ are general affine connections on the base manifold, including {\it non-Riemannian} metric connections. Invariant norms analogous to the invariant norms of Riemannian geometry are problematic for Lorentzian metrics or general affine connections. So alternatively, our theory of the RT-equations employs coordinate dependent $L^p$ norms on $(\Gamma,\A)$ and $(d\Gamma,d\A)$ instead, and establishes that uniform bounds in these norms suffice to bypass the problem of formulating invariant norms in non-Riemannian geometries. This also allows us to incorporate non-optimal connections at the start, as well as non-compact gauge groups, c.f. \cite{ReintjesTemple_Uhl1}. Keep in mind that the gauge transformations which regularize $\A_i$ are independent of the affine part $\Gamma_i$ even in the non-Riemannian case, a consequence of the fact that $\Gamma$ uncouples from $\A$, and $\A$ is coupled to $\Gamma$  in the curvature only through the lower order commutator terms. So it appears that Uhlenbeck's argument in \cite{Uhlenbeck}, based on Coulomb gauge, could in principal be applied to the coordinate Euclidean metric to obtain a compactness result analogous to Theorem \ref{Thm_compactness}, under Uhlenbeck's stated assumptions; namely, for {\it compact} gauge groups, and under the assumption that the connection has optimal regularity at the start with a uniform bound only on the curvature. Our proof based on the RT-equations relaxes these assumptions, is tailored to coordinate based norms, unifies the treatment of optimal regularity and Uhlenbeck compactness in the fibre and affine components of a connection, and provides a new route to both.  

Our proof in the case of non-compact Lie groups requires an interesting new twist to the analysis of the affine case in \cite{ReintjesTemple_ell2} and \cite{ReintjesTemple_Uhl1}. To prove that our iteration scheme underlying the existence theory only produces solutions lying within the gauge group requires introducing an auxiliary elliptic equation in the auxiliary variable $w=U^T\eta U- \eta$. We then prove by a spectral argument that solutions $U$ generated by our iteration scheme must satisfy $w=0$. By this we impose the condition that solutions $U$ of the RT-equations lie within the Lie group $SO(r,s)$, in both the compact ($s=0$) and non-compact case ($s>0$), (c.f. (\ref{definew}) -  (\ref{defineeta}) below). The interesting point here is that, even though we are unable to prove $w\to 0$ by estimates based on the iteration scheme alone, we are able to deduce $w=0$ by an auxiliary elliptic equation which $w$ only satisfies in the limit. In the case when the Lie group is compact, we show the auxiliary equation for $w$ is strongly elliptic, and from this it is straightforward to prove that $w=0$ is the only solution, implying that solutions of the RT-equations always lie within the Lie group $SO(N)$. But in the case of {\it non-compact} Lie groups $SO(r,s)$, the auxiliary equation for $w$ may not be strongly elliptic, and could have an associated non-trivial spectrum with non-zero solutions (eigenfunctions) by the Fredholm alternative. We circumvent this problem by introducing an additional spectral parameter $\lambda,$  and by proving that our iteration scheme produces non-trivial solutions at most on a set of measure zero in $\lambda$.  We then prove that the uniform convergence of our iteration scheme by which we generate solutions, implies the continuity of solutions $w$ of the auxiliary equation with respect to $\lambda$. Thus, by continuity of $w$ with respect to $\lambda$, we can extend $w=0$ from almost every $\lambda$, to $w=0$ for all $\lambda$. By this we establish that solutions of the RT-equations generated by our iteration scheme always lie within the Lie group $SO(r,s)$, even when the Lie group is non-compact and the auxiliary equation in $w$ admits non-trivial solutions.  

Most interesting to us is that the RT-equations reproduce a fundamental cancellation first observed for affine connections in \cite{ReintjesTemple_ell1}. Namely, a cancellation of the ``bad terms'' $d\delta \A$ which lie at a regularity too low for the Laplacian in \eqref{RT_1_intro} to lift $\Ati$ to optimal regularity. In the setting of vector bundles, this cancellation is due to an interesting interplay between the Lie group and Lie algebra of the fibres as expressed through the RT-equations. Thus the assumed bound on $d\A$ alone is sufficient to dominate the uncontrolled derivatives $\delta\A$, and this appears to be a principle built into the RT-equations through the geometry of connections.   Also, to get to the low level of $L^p$ regularity in the case of affine connections \cite{ReintjesTemple_Uhl1}, it was necessary to use invariance transformations of the RT-equations to decouple the equation for the Jacobian from the equation for the regularized connection. But here, to our surprise, no such invariance transformations of the RT-equation is required, and the equation for the regularizing gauge transformation ``automatically'' de-couples from the equation for the regularized connection.  Moreover, in the vector bundle case here we do not require an additional condition to impose the requirement that solutions of the RT-equations lie in the gauge group,\footnote{In \cite{ReintjesTemple_Uhl1} we need an auxiliary equation to impose integrability of the Jacobians to coordinates.} but rather we build this into the RT-equations themselves by solving explicitly an equation for the free parameter $\alpha$ associated with the invariance transformations of the equations, (c.f. Section \ref{Sec_RT-eqn_derivation}). Finally note that our derivation of the vector bundle version of the RT-equations from the connection transformation law involves educated guesses based on analogy with the RT-equations for affine connections. In the affine case the fortuitous choices of unknown variables and differential operators were motivated by the Riemann-flat condition in \cite{ReintjesTemple_geo}.

\section{Statement of Results}  \label{Sec_Results}

The problem of optimal regularity and Uhlenbeck compactness for a connection $\A$ on the fibres of a vector bundle, \emph{uncouples} from the problem of optimal regularity for the affine connection $\Gamma$ on the tangent bundle of the base manifold, because the connection transformation laws for $\Gamma$ and $\A$ are independent of each other. This carries through to the RT-equations, because the RT-equations employ only the leading order terms of the curvature, $d\A$ and $d\Gamma$ respectively, but not the full curvature (which couples $\A$ and $\Gamma$ through the lower order commutator term \cite{Husemoller}). Moreover, both the affine and the vector bundle versions of the RT-equations are based on the Euclidean metric in an arbitrary but fixed coordinate system, and this auxiliary Riemannian structure acts independently of the affine connection on the base manifold.  Taken on whole, to extend our results to vector bundles it suffices to assume, without loss of generality, that the base manifold is Euclidean.  The extension to non-trivial non-optimal affine connections on the base manifold is then a straightforward application of our results in \cite{ReintjesTemple_Uhl1}. In Section \ref{Sec_results_1} we state our results for vector bundles in the case of Euclidean base manifolds, and in Section \ref{Sec_results_2} we explain how to incorporate our earlier results for affine connections in \cite{ReintjesTemple_Uhl1} into the statements of the theorems.

\subsection{Connections on the fibre of a vector bundle}   \label{Sec_results_1}

Let $\A_{\VM}$ denote a given connection on a vector bundle $\VM$ of a $n$-dimensional differentiable base manifold $\M$ with $N$-dimensional real valued fibres on which the Lie group $SO(r,s)$ acts as the gauge group, $r+s=N$. As above, we assume without loss of generality that the base manifold is endowed with the Euclidean metric, and thus that the affine connection components vanish.  Since our problem is local and can be considered in any fixed coordinate system on $\M$, assume without loss of generality that the vector bundle is trivial, 
$$
\VM \simeq \R^N \times \Omega,
$$ 
for some $\Omega \subset  \R^n$ open and bounded with smooth boundary ($C^{1,1}$ suffices). One can view $\Omega$ as the image of a coordinate patch $(x,\mathcal{U})$ of $\M$, i.e., $\Omega= x(\mathcal{U})$ with resulting Cartesian coordinates $x$ on $\Omega$. The gauge group $SO(r,s)$ is defined by the condition that an $N\times N$ matrix $U$ is an element of $SO(r,s)$ if and only if
\beq\label{definew}
U^T\eta U=\eta  
\hspace{1cm} \text{and} \hspace{1cm}   
{\rm det}(U)=1,
\eeq
where $\eta$ is the diagonal matrix with $r$ entries $1$ and $s$ entries $-1$, 
\beq \label{defineeta}
\eta={\rm diag}(1,...,1,-1,...,-1),
\eeq 
for integers $r\geq 0$ and $s\geq 0$ with $r+s=N$. We let $\A_\textbf{a}$ denote the connection components of $\A_{\VM}$ on $\VM$ with respect to a choice of frame $\textbf{a}$ on $\VM$, a so-called \emph{gauge}, i.e., $\textbf{a}$ assigns a basis of $\R^N$ at each point of $\Omega$. The connection components $\A_\textbf{a}$ are matrix valued $1$-forms, which lie in the Lie Algebra of $SO(r,s)$, consisting of $N\times N$ matrices $A$ satisfying the condition
\beq \label{LieAlgebra}
A^T \cdot \eta + \eta \cdot A = 0 \hspace{1cm} \text{and} \hspace{1cm}  {\rm tr}(A) =0,
\eeq  
where ${\rm tr}(\cdot)$ denotes the matrix trace and ``$\cdot$'' matrix multiplication (which we may omit or include). Assume now $\textbf{b}$ is another frame of $\VM$ resulting from $\textbf{a}$ by a gauge transformation 
\beq
U : \Omega \rightarrow SO(r,s),
\eeq
that is, $\textbf{b} = U\cdot \textbf{a}$ in the sense of pointwise matrix multiplication with each basis vector in $\textbf{a}$. The components of $\A_{\VM}$ in the new frame $\textbf{b}$ are related to $\A_\textbf{a}$ by the connection transformation law   
\beq \label{connection_transfo_VB}
\A_\textbf{a} = U^{-1} dU + U^{-1} \Aop U,
\eeq
where $\Aop$ denotes the connection in the gauge $\textbf{b} = U\cdot \textbf{a}$, $d$ denotes the exterior derivative on the matrix valued differential $0$-form $U$, $dU = \frac{\partial U}{\partial x^i} dx^i$, and we use Einstein's convention of summing over repeated upper and lower indices.\footnote{An abstract connection $\A_{\VM},$ by {\it definition}, assigns connection components $\A_\textbf{a}$ to every frame $\textbf{a}$ such that \eqref{connection_transfo_VB} holds, (c.f. \cite{Husemoller,MaratheMartucci} for an introduction to vector bundles and \cite{Taubes,Wehrheim} for more background on Uhlenbeck compactness).}

We now fix an arbitrary gauge $\textbf{a}$ and assume the connection components $\A_\textbf{a}$ of $\A_{\VM}$ are non-optimal, i.e.,  $\A_\textbf{a} \in L^{2p}(\Omega)$ and $d\A_\textbf{a} \in L^{\bar{p}}(\Omega)$, in the sense that the components of $\A_\textbf{a}$ have this regularity with respect to fixed Cartesian coordinates on $\Omega$, which we view as the coordinates resulting from a coordinate patch $(x,\mathcal{U})$ on the base manifold $\M$ such that $\Omega = x(\mathcal{U})$. The goal is then to find a gauge transformation 
$$
U : \Omega \rightarrow SO(r,s)
$$ 
such that the connection components $\Aop$ in the resulting gauge $\textbf{b} = U\cdot \textbf{a}$ exhibit optimal regularity, $\A_\textbf{b} \in W^{1,\hat{p}}(\Omega)$ for $\hat{p}=\min\{p,\bar{p}\}$, one derivative of regularity above the curvature in $L^{\hat{p}}$. Here we establish optimal regularity for $p \in (n/2,\infty]$, $\bar{p} \in (1,\infty)$, $n\geq 2$.  We subsequently always denote the connection components $\A_\textbf{a}$ in the fixed incoming gauge $\textbf{a}$ simply by $\A$, and refer to the collection of connection components $\A$ simply as a connection.

In this paper, we begin by deriving the following system of elliptic Poisson type PDE's for the gauge transformations $U$ in $SO(r,s)$ which transform a non-optimal connection $\A \equiv \A_\textbf{a}$ to a connection $\Aop$ of optimal regularity:
\begin{eqnarray}
\Delta \Ati &=& \delta d \A -\delta\big( dU^{-1}\wedge dU \big)  , \label{RT_1} \\
\Delta U &=&   U \delta\A   -  (U^T \eta )^{-1}  \langle dU^T ; \eta dU\rangle . \label{RT_2}
\end{eqnarray}
Equations \eqref{RT_1} - \eqref{RT_2} are the vector bundle version of the RT-equations, (already stated in the Introduction as equations \eqref{RT_1_intro} - \eqref{RT_2_intro}), and for simplicity we refer to them here as the {\it RT-equations}. The unknowns in equations \eqref{RT_1} - \eqref{RT_2} are $(\Ati,U)$, where $\Ati$ is a matrix valued $1$-form associated to $\Aop$, and $U$ is the sought after gauge transformation interpreted as a matrix valued $0$-form. Equations \eqref{RT_1} - \eqref{RT_2} are elliptic, with $\Delta$ being the standard Laplacian  in $\R^n$  acting component-wise, $\Delta \equiv d\delta + \delta d = \partial^2_{x^1} +... + \partial^2_{x^n}$, and $d$ is the exterior derivative and $\delta$ its co-derivative based on the Euclidean metric in Cartesian $x$-coordinates. The matrix valued ``inner product'' $\langle \cdot \: ; \cdot \rangle$ is introduced in \eqref{def_inner-product} in Section \ref{Sec_Prelim} below.  Equation \eqref{RT_2} is what we interpret as the vector bundle version of the {\it reduced} RT-equations, an equation decoupled from \eqref{RT_1}. Solutions $U$ of \eqref{RT_2} are the $SO(r,s)$ gauge transformations which map the connection $\A$ to optimal regularity. 

The logic of our proof is as follows. Starting with a non-optimal connection $\A$, we prove below that solutions of the reduced RT-equations  \eqref{RT_2} with Dirichlet boundary data $U = U_0$ on $\partial\Omega$ for some $U_0$ in $SO(r,s)$, yield a gauge transformation $U$ in $SO(r,s)$, (by which we mean that $U(q) \in SO(r,s)$ pointwise for every $q \in \Omega$). A solution $U$ in $SO(r,s)$ of the reduced RT-equations \eqref{RT_2} gives rise to a solution of  \eqref{RT_1}, defined by 
\beq \label{def_Ati'}
\Ati' = \A - U^{-1} dU,
\eeq 
that is, we prove below that $\Ati'$ satisfies   
\begin{eqnarray}\label{RT_Gammati'}
\Delta \Ati' &=&  \delta d \A - \delta (dU^{-1} \wedge dU) .
\end{eqnarray}
The connection of optimal regularity in the gauge $\textbf{b} = U\cdot \textbf{a}$ is then given by
\beq \label{OptimalGamma}
\Aop = U \Ati' U^{-1}.
\eeq
Since the right hand side of (\ref{RT_Gammati'}) has regularity just one derivative {\it below} $\A$, elliptic regularity will imply that $\Ati',$ and then by (\ref{OptimalGamma}) also $\Aop,$ are one derivative of regularity {\it above} $\A$,  thus establishing the optimal regularity of the connection $\Aop$ in the resulting gauge $\textbf{b} = U\cdot \textbf{a}$.  Thus the role of \eqref{RT_1} is to raise the regularity of $\Ati'$, and hence $\Aop$, by one derivative, to optimal regularity. This result is stated precisely in the following theorem. The idea of proof together with the derivation of the RT-equations \eqref{RT_1} - \eqref{RT_2} is the subject of Section \ref{Sec_RT-eqn}.  We then prove in Section \ref{Sec_weak} that a weak formulation of the RT-equations exists which preserves the logic of the theory for $L^p$ connections and distributional derivatives.

\begin{Thm} \label{Thm_equiv}
Let $\A \equiv \A_\textbf{a}$ be the connection components in a gauge $\textbf{a}$ of a connection $\A_{\VM}$ on an $SO(r,s)$ vector bundle $\VM$ with base manifold $\M \equiv \Omega \subset \R^n$ open and bounded.    Assume $\A \in L^{2p}(\Omega)$ with $d\A \in L^{\bar{p}}(\Omega)$, for $p\in (n/2,\infty]$ and $\bar{p} \in (1,\infty)$, $n\geq 2$.  Then the following equivalence holds:  \vspace{.1cm}

\noindent{\bf (i)} If there exists a solution  $U \in W^{1,2p}(\Omega)$ pointwise in $SO(r,s)$ of the reduced RT-equations \eqref{RT_2}, then the gauge transformed connection $\Aop$ in \eqref{OptimalGamma} has optimal regularity $\Aop \in W^{1,\hat{p}}(\Omega)$ for $\hat{p}=\min\{p,\bar{p}\}$. \vspace{.1cm}

\noindent{\bf (ii)} Conversely, if there exists a gauge transformation $U\in W^{1,2p}(\Omega)$ pointwise in  $SO(r,s)$, such that the gauge transformed connection $\Aop$ in \eqref{connection_transfo_VB} has optimal regularity $\Aop \in W^{1,\hat{p}}(\Omega)$, then $\Ati \equiv U^{-1} \Aop U \in W^{1,\hat{p}}(\Omega)$ and $U$ solve the RT-equations \eqref{RT_1} and \eqref{RT_2}, respectively. 
\end{Thm}

Theorem \ref{Thm_equiv} reduces the proof of our main results on optimal regularity and Uhlenbeck compactness, (Theorems \ref{Thm_opt} and \ref{Thm_compactness} below), to proving the following existence theorem, which is the main technical effort of this paper.  We distinguish between the global and local regularity problem in the neighborhood $\Omega$. That is, the global problem applies to a connection with sufficiently small $L^p$-norm throughout $\Omega$ and yields optimal connection regularity on any compact subset of the entire set $\Omega$, while the local case applies to general connections without a smallness assumption on their $L^p$ norms, and establishes optimal regularity in a subneighborhood of any point in $\Omega$. The local result does not appear to imply a global result, because after a connection is lifted locally to optimal regularity in different coordinate charts, the transition maps between them can loose any higher regularity assumed at the start, c.f. \cite{Donaldson}.\footnote{In \cite{Uhlenbeck} connection are assumed in $W^{1,p}$, which yields global results from local ones.} We treat the cases $p \in (n/2,\infty)$ and $p=\infty$ separately, since the case $p=\infty$ is a singular case of elliptic regularity theory, c.f. \cite{Evans}.   Note that the curvature \eqref{curvature} lies in $L^p$ when $d\A$ lies in $L^p$ and $\A$ in $L^{2p}$, (by applying the H\"older inequality to $\A \wedge \A$), and so to keep the notation consistent we assume that $\A$ lies $L^{2p}$. We now state our existence theorem for the reduced RT-equations.

\begin{Thm}  \label{Thm_existence}
{\bf Global (i):} Let $\A \in L^{2p}(\Omega)$, for $p\in (n/2,\infty)$. Then there exists an $\epsilon>0$ depending only on $\Omega $ and $p$, such that if
\beq \label{initial_bound_small}
 \|\A\|_{L^{2p}(\Omega)}  \leq \; \epsilon ,
\eeq    
then there exists a solution $U \in W^{1,2p}(\Omega)$ of the reduced RT-equations \eqref{RT_2} which is pointwise in $SO(r,s)$, and satisfies
\beq \label{uniform_bound_U_global} 
\|U - I\|_{W^{1,2p}(\Omega)}  + \|U^{-1}-I\|_{W^{1,2p}(\Omega)} \; \leq\; C\; \|\A\|_{L^{2p}(\Omega)} 
\eeq
for some constant $C > 0$ depending only on $\Omega$ and $p$.

\vspace{.1cm} \noindent {\bf Global (ii):}  If \eqref{initial_bound_small} holds for $p=\infty$, then, for any $\tilde{p} <\infty$, there exists a solution $U \in W^{1,2\tilde{p}}(\Omega)$ of \eqref{RT_2} in $SO(r,s)$ which satisfies \eqref{uniform_bound_U_global} for $\tilde{p}$ in place of $p$.\footnote{We do not establish $U\in W^{1,\infty}$, because $p=\infty$ (like $p=1$) is a singular case of elliptic regularity theory.}

\vspace{.2cm} \noindent {\bf Local (i):} Assume $\A \in L^{2p}(\Omega)$, for $p\in (n/2,\infty)$, and let $M>0$ be a constant such that
\beq \label{initial_bound_exist}  
 \|\A\|_{L^{2p}(\Omega)}  \leq \; M.
\eeq    
Then for any point $q \in \Omega$ there exists a neighborhood $\Omega' \subset \Omega$ of $q$, (which in general depends on $\A$), and a solution $U \in W^{1,2p}(\Omega')$ of the reduced RT-equations \eqref{RT_2}, which is pointwise in $SO(r,s)$ and which satisfies
\beq \label{uniform_bound_U_local}
\|U - I\|_{W^{1,2p}(\Omega')}  + \|U^{-1}-I\|_{W^{1,2p}(\Omega')} \; \leq\; C(M)\; \|\A\|_{L^{2p}(\Omega')} 
\eeq
for some constant $C(M) > 0$ depending only on $\Omega', p$ and $M$. 

\vspace{.1cm} \noindent {\bf Local (ii):}  If \eqref{initial_bound_exist} holds for $p=\infty$, then for any $q \in \Omega$ there exists a neighborhood $\Omega' \subset \Omega$, and for any $\tilde{p} <\infty$ there exists a solution $U \in W^{1,2\tilde{p}}(\Omega')$ of \eqref{RT_2} in $SO(r,s)$, which satisfies \eqref{uniform_bound_U_local} for $\tilde{p}$ in place of $p$. In this case, $\Omega'$ depends only on $\Omega, \tilde{p}$ and $M$, but is independent of $\A$. 
\end{Thm}

For concreteness, in the local part of Theorem \ref{Thm_existence}, the neighborhood $\Omega' \subset \R^n$ can be taken to be $B_r(q)$, the Euclidean ball of radius $r$ with center $q$, where in Local (i) the $r>0$ depends on $\Omega, p$ and  $M,$ but also varies with the connection $\A$, depending on the $p$-norm of $\A$ near $q$; while in Local (ii), $r$ depends only on $\Omega, p$ and $M$, but not on $\A$. 

By the existence result in Theorem \ref{Thm_existence}, we obtain our optimal regularity result for connections on $SO(r,s)$ vector bundles by transformation with the gauge transformations $U$ given in Theorem \ref{Thm_existence}. For example, assuming as in \cite{ReintjesTemple_Uhl1} $\A \in L^{2p}$ and $d\A \in L^p$, it follows that the curvature lies in $L^p$, and our optimal regularity result below establishes $\Aop \in W^{1,p}$, $p\in (n/2,\infty)$, exactly one derivative above the curvature in $L^p$. Note that if $\A$ is more regular than $L^{2p}$, it will not improve the optimal regularity, but it will improve the regularity of the gauge transformation $U$ by Theorem \ref{Thm_existence}. More importantly, because the RT-equation \eqref{RT_2} for $U$ decouples from $d\A$, it turns out that elliptic regularity applied to the first RT-equation \eqref{RT_1} allows us to lower the regularity of $d\A$ from $L^p$ to $L^{\bar{p}}$ for $\bar{p} \in (1,\infty)$, while keeping $\A \in L^{2p}$ for  $p\in (n/2,\infty)$ to guarantee existence of $U$. This places the curvature \eqref{curvature} in $L^{\hat{p}}$ for $\hat{p} = \min\{p,\bar{p}\}$, and Theorem \ref{Thm_opt} places the optimal connection $\Aop$ in $W^{1,\hat{p}}$, exactly one derivative above the curvature.  The following theorem is formulated to incorporate all of these cases.

\begin{Thm}  \label{Thm_opt}  
The solutions $U$ established in Theorem \ref{Thm_existence} transform a connection $\A \in L^{2p}(\Omega)$ to a connection $\Aop$ of optimal regularity by the connection transformation law \eqref{connection_transfo_VB} as follows:

\vspace{.1cm} \noindent {\bf Global (i):}
Let $p \in (n/2,\infty)$, $\bar{p} \in (1,\infty)$ and $n\geq 2$. Assume
\beq \label{initial_bound_opt} 
\|(\A,d\A)\|_{L^{2p,\bar{p}}(\Omega)}  \equiv \|\A \|_{L^{2p}(\Omega)} + \|d\A \|_{L^{\bar{p}}(\Omega)} \; \leq \; M,
\eeq  
for some constant $M>0$, (so the curvature \eqref{curvature} lies in $L^{\hat{p}}(\Omega)$ for $\hat{p} = \min\{p,\bar{p}\}$). Assume further $\A$ satisfies the smallness condition \eqref{initial_bound_small}. Then the connection of optimal regularity $\Aop$ lies in $W^{1,\hat{p}}(\Omega'')$ on every neighborhood $\Omega''$ compactly contained in $\Omega$, and satisfies
\beq \label{uniform_bound_A_global}  
\|\Aop \|_{W^{1,\hat{p}}(\Omega'')}  \; \leq\; C(M)\; \|(\A,d\A)\|_{L^{2p,\bar{p}}(\Omega)} 
\eeq
for some constant $C(M) > 0$ depending only on $\Omega'', \Omega, p, \bar{p}$ and $M$.\footnote{Note that we require an $L^{2p}$ bound on $\A$ for our existence theory to handle the product term $U\delta\A$ in \eqref{RT_2} by Morrey's inequality, which works when $p>n/2$, but might fail when $p=n/2$. In contrast, no products on $d\A$ appear in \eqref{RT_1} - \eqref{RT_2}, and for this reason we only require an $L^{\bar{p}}$ bound on $d\A$ for $\bar{p}>1$. In \cite{Uhlenbeck}, a connection in $W^{1,p}$ with a uniform bound on the curvature in $L^p$ for $p\geq n/2$ is assumed.}         

\vspace{.1cm} \noindent {\bf Global (ii):} 
Let $p =\infty$ and $\bar{p} \in (1,\infty)$. Assume $\A$ satisfies \eqref{initial_bound_opt} and the smallness condition \eqref{initial_bound_small}. Then $\Aop$ satisfies \eqref{uniform_bound_A_global} with $\hat{p}=\bar{p}$ and $p =\infty$.

\vspace{.1cm} \noindent {\bf Local (i):}
Let $p \in (n/2,\infty)$, $\bar{p} \in (1,\infty)$ and $n\geq 2$. Assume $\A$ satisfies \eqref{initial_bound_opt}, but not the smallness condition \eqref{initial_bound_small}.   Then for any point $q\in \Omega$, there exists an $\A$-dependent neighborhood $\Omega' \subset \Omega$, such that $\Aop$ lies in $W^{1,\hat{p}}(\Omega'')$ on every neighborhood $\Omega''$ compactly contained in $\Omega'$, such that   
\beq \label{uniform_bound_A_local} 
\|\Aop \|_{W^{1,\hat{p}}(\Omega'')}  \; \leq\; C(M)\; \|(\A,d\A)\|_{L^{2p,\bar{p}}(\Omega')} 
\eeq
for some constant $C(M) > 0$ depending only on $\Omega'', \Omega', p, \bar{p}$ and $M$.  

\vspace{.1cm} \noindent {\bf Local (ii):} 
Let $p =\infty$ and $\bar{p} \in (1,\infty)$.  Assume $\A$ satisfies \eqref{initial_bound_opt}, but not \eqref{initial_bound_small}. Then for any point $q\in \Omega$, there exists a neighborhood $\Omega' \subset \Omega$, independent of $\A$, such that $\Aop$ lies in $W^{1,\hat{p}}(\Omega'')$ and satisfies \eqref{uniform_bound_A_local} with $\hat{p}=\bar{p}$ and $p =\infty$.
\end{Thm}

Theorem \ref{Thm_opt} establishes the existence of local gauge transformations which transform non-optimal connections to optimal regularity, one derivative more regular than the curvature. We cannot use Local (i) of Theorem \ref{Thm_opt} to establish Uhlenbeck compactness, because the domain of the regularized connection depends on the connection, so we cannot maintain a uniform domain for a sequence of connections that meet estimate \eqref{uniform_bound_A_local}. In contrast, the domains in Local (ii) and Global (i) and (ii) are independent of the connection and hence amenable for use in Uhlenbeck compactness. In all of these three cases the extra derivative of regularity for the connection provided by the $W^{1,\hat{p}}$ bounds on $\Aop$, together with the bounds on $U$ in Theorem \ref{Thm_existence}, provide the uniform bounds from which we infer our versions of Uhlenbeck compactness.  We thus obtain a global version  of Uhlenbeck compactness from Theorem \ref{Thm_opt} - Global (i) and (ii), and we  obtain a local version from Theorem \ref{Thm_opt} - Local (ii):

\begin{Thm}   \label{Thm_compactness}
{\bf Global (i):}  
Let $p \in (n/2,\infty)$ and $\bar{p} \in (1,\infty)$. Let $(\A_i)_{i\in \mathbb{N}}$ be a sequence of connections on $\VM$ in a fixed gauge $\textbf{a}$ on a fixed domain $\Omega$. Assume the $\A_i$ satisfy the uniform bound 
\beq \label{initial_bound_2} 
\|(\A_i,d\A_i)\|_{L^{2p,\bar{p}}(\Omega)} \equiv \|\A_i \|_{L^{2p}(\Omega)} + \|d\A_i \|_{L^{\bar{p}}(\Omega)} \; \leq \; M
\eeq   
for a constant $M>0$, and assume each $\A_i$ satisfies the smallness condition \eqref{initial_bound_small}.  Then the following holds on every neighborhood $\Omega''$ compactly contained in $\Omega$:  

\begin{enumerate}
\item[{\bf (a)}] For each $\A_i$ there exists an $SO(r,s)$ gauge transformation $U_i \in W^{1,2p}(\Omega'')$ such that the components $\A_{\textbf{b}_i}$ of $\A_i$ in gauge $\textbf{b}_i = U_i \cdot \textbf{a}$ have optimal regularity $\A_{\textbf{b}_i}\in W^{1,\hat{p}}(\Omega'')$ for $\hat{p}=\min\{p,\bar{p}\}$, with uniform bound 
\beq \label{uniform_bound_A_2}
\|\A_{\textbf{b}_i} \|_{W^{1,\hat{p}}(\Omega'')}   \; \leq \;  C(M) M ,
\eeq
for a constant $C(M) > 0$ depending only on $\Omega'', \Omega, p, \bar{p}$ and $M$, independent of $\A_i$. 

\item[{\bf (b)}] The sequence $(U_i)_{i\in \mathbb{N}}$ is uniformly bounded in $W^{1,2p}(\Omega'')$ by \eqref{uniform_bound_U_global}, and a subsequence converges weakly in $W^{1,2p}$ to some $U \in W^{1,2p}(\Omega'')$ in $SO(r,s)$. 

\item[{\bf (c)}] {\rm Main Conclusion:}   There exists a subsequence of $\A_i$, (denoted again by $\A_i$), such that the components of $\A_{\textbf{b}_i}$ converge to some $\A_{\textbf{b}}$ weakly in $W^{1,\hat{p}}(\Omega'')$, strongly in $L^{\hat{p}}(\Omega'')$. The limit $\A_{\textbf{b}}$ is the connection $\A$ in gauge $\textbf{b} = U\cdot \textbf{a}$, where $\A$ is the weak limit of $\A_i$ in $L^{2p}(\Omega'')$ in gauge $\textbf{a}$.
\end{enumerate}

\vspace{.1cm} \noindent {\bf Global (ii):} Assume a sequence of $L^\infty$ connections $(\A_i)_{i\in \mathbb{N}}$ satisfies \eqref{initial_bound_2} for $p=\infty$, $\bar{p}  \in (1,\infty)$, and the smallness condition \eqref{initial_bound_small} for $p=\infty$. Then, on every compactly contained neighborhood $\Omega''$ in $\Omega$, the conclusions (a) - (c) of case Global (i) hold for $p=\bar{p}$ and $\hat{p}=\bar{p}$.\footnote{Theorem \ref{Thm_existence} implies for the connection regularity $\A_i \in L^\infty$ assumed in case (ii) existence of gauge transformations $U_i\in W^{1,2p}$ for any $p<\infty$. However, because Uhlenbeck compactness concerns only the convergence of $\A_{\textbf{b}_i}$, it suffices in the cases (ii) of Theorem \ref{Thm_compactness} to restrict to $U_i\in W^{1,2\bar{p}}$ without any relevant loss.} 

\vspace{.1cm} \noindent {\bf Local (ii):} Assume a sequence of $L^\infty$ connections $(\A_i)_{i\in \mathbb{N}}$ satisfies \eqref{initial_bound_2} for $p=\infty$ and $\bar{p} \in (1,\infty)$, but not the smallness condition \eqref{initial_bound_small}. Then for any point $q \in \Omega$ there exists a neighborhood $\Omega'' \subset \Omega$ of $q$, independent of $\A_i$, such that (a) - (c) of case Global (i) hold for $p=\bar{p}$ and $\hat{p}=\bar{p}$.  
\end{Thm}

The weak convergence of connections in $W^{1,\hat{p}}(\Omega'')$ asserted by Theorem \ref{Thm_compactness} implies strong convergence in $L^{\hat{p}}(\Omega'')$, and this is a convergence regular enough to pass limits through non-linear products, a property inherently useful for non-linear analysis.   The assumption in Theorem \ref{Thm_compactness} of either an $L^\infty$ bound or a sufficiently small $L^{2p}$ bound on the connection is required to ensure a uniform domain on which optimal regularity is established. We do not view this as a major concession, because our main concern is to address $L^{\hat{p}}$ curvatures ranging all the way down to $\hat{p}>1$. 

The proofs of Theorems \ref{Thm_equiv} -  \ref{Thm_compactness} are given in Sections \ref{Sec_proof_existence} - \ref{Sec_proofs_rest}.  These theorems also apply at higher levels of non-optimal regularity, for example when $\A, d\A \in W^{m,p}(\Omega)$, $m\geq 1$, $n<p<\infty$, as in \cite{ReintjesTemple_ell2}. For this more regular case Theorem \ref{Thm_equiv} gives the equivalence between the existence of $SO(r,s)$ gauge transformations $U \in W^{m+1,p}(\Omega)$ which smooth a connection to optimal regularity $\Aop \in W^{m+1,p}(\Omega)$, and the existence of solutions $(U,\Ati) \in W^{m+1,p}(\Omega)$ of the RT-equations. The existence result corresponding to Theorem \ref{Thm_existence} yields a solution $U \in W^{m+1,p}(\Omega)$ in $SO(r,s)$ of the RT-equation satisfying the estimate 
\beq 
\|U\|_{W^{m+1,p}(\Omega')}  \; \leq\; C(M)\; \|(\A,d\A)\|_{W^{m,p}(\Omega')} ,
\eeq
under the assumption that
\beq 
 \|(\A,d\A)\|_{W^{m,p}(\Omega)}  \equiv \|\A \|_{W^{m,p}(\Omega)} + \|d\A \|_{W^{m,p}(\Omega)} \; \leq \; M.
\eeq   
Based on this existence result, the higher regularity version of Theorem \ref{Thm_opt} establishes optimal regularity in $W^{m+1,p}$; and the higher regularity version of Theorem \ref{Thm_compactness} of Uhlenbeck compactness establishes convergence of sequences of connections $\A_i$ weakly in $W^{m+1,p}$ and hence strongly in $W^{m,p}$.   The proofs for the more regular case $\A, d\A \in W^{m,p}$, a setting which is simpler essentially because $W^{1,p}$ is closed under multiplication by Morrey's inequality, and which we addressed in \cite{ReintjesTemple_ell2} for affine connections, are omitted in this paper.  

Note finally that Theorem \ref{Thm_opt} in combination with \cite[Thm 2.1]{Uhlenbeck}, yields the existence of a transformation to Coulomb gauge in the compact group $SO(n)$, which lowers the $W^{1,p}$ connection regularity assumed in \cite[Thm 2.1]{Uhlenbeck} to non-optimal connections in $L^{2p}$, $p>n/2$.\footnote{ Assuming only a non-optimal connection in $L^{2p}$ with curvature in $L^p$, ($p>2n$), Theorem \ref{Thm_opt} establishes the optimal $W^{1,p}$ regularity assumed in \cite[Thm 2.1]{Uhlenbeck}, and Theorem 2.1 in \cite{Uhlenbeck} then asserts the existence of a transformation in $SO(n)$ to Coulomb gauge. To our knowledge, this is the first existence result for the Coulomb gauge not requiring $W^{1,p}$ connection regularity.}

\subsection{Incorporating non-trivial affine connections into the base manifold}    \label{Sec_results_2}

Theorems \ref{Thm_opt} and \ref{Thm_compactness} on optimal regularity and Uhlenbeck compactness for vector bundles $\VM$ extend directly to the case of non-trivial affine connections $\Gamma$ defined on the base manifold $\M$, by using our results for affine connections in \cite{ReintjesTemple_Uhl1}.  In \cite{ReintjesTemple_Uhl1} we established local optimal regularity for connections $\Gamma \in L^{2p}$ with $d\Gamma \in L^p$, for $p\in (n/2,\infty)$. The global problem is not addressed in \cite{ReintjesTemple_Uhl1}, even though one could prove it, because the authors were concerned with regularity as a purely local issue. It is straightforward to combine the results in \cite{ReintjesTemple_Uhl1} with those in this paper with the assumptions on $\A$ and $\Gamma$ unchanged, with the same conclusion about existence, optimal regularity and Uhlenbeck compactness. For concreteness we state here the case Local (i) for optimal regularity and the case Local (ii) for Uhlenbeck compactness, with the exact same regularities as addressed in \cite{ReintjesTemple_Uhl1}, as special cases. 

For case Local (i) of optimal regularity, assume $\Gamma$ is a non-optimal affine connection on $\M$ of an $SO(r,s)$ vector bundle $\VM$ with connection $\A$ on the fibres, and assume (\ref{initial_bound_opt}) holds for $\hat{p}=\bar{p}=p \in (n/2,\infty)$ together with
\begin{eqnarray}\label{GenStartGamma}
\|(\Gamma,d\Gamma)\|_{L^{2p,p}}&=&\|\Gamma\|_{L^{2p}}+\|d\Gamma\|_{L^p} \leq M.
\end{eqnarray}
The results in \cite{ReintjesTemple_Uhl1} for affine connections and Theorem \ref{Thm_opt} for connections on the fibre act independently. Thus combining Theorem \ref{Thm_opt} - Local (i) in this paper with Theorem 3.1 in \cite{ReintjesTemple_Uhl1}, we obtain the following theorem on optimal regularity.

\begin{Thm} \label{Thm_opt_general}  
Let ${\mathcal V}{\mathcal M}$ be a vector bundle over an $n$-dimensional manifold ${\mathcal M}$ with $N$-dimensional fibre acted on by $SO(r,s)\subset {\mathbb R}^{N\times N}$ as the gauge group, and let $(x,\textbf{a}( x))$  be a given coordinate system $x$ on ${\mathcal M}$  paired with a given gauge $\textbf{a}$. Assume a connection $(\Gamma,\A)$ on ${\VM}$ which satisfies
\beq \label{bound_1}
\|(\Gamma_x,\A_{(x,\textbf{a})})\|_{p,\Omega} \equiv
\| (\Gamma_x,\A_{(x,\textbf{a})}) \|_{L^{2p}(\Omega)} + \|(d\Gamma_x,d\A_{(x,\textbf{a})}) \|_{L^p(\Omega)} \; \leq \; M,
\eeq
where norms are taken components-wise in $x$-coordinates, $p>\max\{n/2,2\}$, $n \geq 2$, $p<\infty$. Then, for any point in $\Omega$ there exists a neighborhood $\Omega' \subset \Omega$ and a coordinate transformation $x \to y$ with Jacobian $J=\frac{\partial y}{\partial x} \in W^{1,2p}(\Omega')$ and a gauge transformation $U \in W^{1,2p}(\Omega')$ in $SO(r,s)$, such that the connection components $(\Gamma_y,\A_{(y,\textbf{b})})$   in $y$-coordinates and gauge $\textbf{b} = U \cdot \textbf{a}$ have optimal regularity, 
\beq
(\Gamma_y,\A_{(y,\textbf{b})}) \in W^{1,p}(\Omega''),
\eeq 
on every open set $\Omega''$ compactly contained in $\Omega'$. The connection of optimal regularity satisfies
\beq \label{bound_2}
\|(\Gamma_y,\A_{(y,\textbf{b})}) \|_{W^{1,p}(\Omega'')} 
\; \leq\; C(M) \|(\Gamma_x,\A_{(x,\textbf{a})})\|_{p,\Omega'} ,
\eeq
and the regularizing transformations $J$ and $U$ satisfy
\begin{align} \label{bound_3}
&\|J-I\|_{W^{1,2p}(\Omega')}  + \|J^{-1}-I\|_{W^{1,2p}(\Omega')} \; \leq\; C(M) \|(\Gamma_x,\A_{(x,\textbf{a})})\|_{p,\Omega'} \\
&\|U-I\|_{W^{1,2p}(\Omega')}  + \|U^{-1}-I\|_{W^{1,2p}(\Omega')}   \; \leq\; C(M) \|(\Gamma_x,\A_{(x,\textbf{a})})\|_{p,\Omega'} , \label{bound_3b}
\end{align}
where $C(M) > 0$ is a constant depending only on $\Omega'', \Omega', p$ and $M$, and norms on $J$, $U$ and $(\Gamma_y,\A_{(y,\textbf{b})})$ can be taken in $x$- or $y$-coordinates.         
\end{Thm}

\Proof
By the results in \cite{ReintjesTemple_Uhl1}, there exist coordinate transformations $x\to y$ which locally lift the regularity of the components of $\Gamma$ to  optimal regularity, $\Gamma\in W^{1,p}$, (one derivative above the Riemann curvature).   Now, the Jacobians which accomplish this are regular enough so that in $y$-coordinates, the estimate (\ref{initial_bound_opt}) on the components of $\A$ continues to hold in $y$-coordinates with a (uniformly) modified upper bound $M$.   Thus, since the arguments establishing Theorems \ref{Thm_opt} in the prior subsection are based on the (auxiliary) Euclidean coordinate metric, we can apply the same arguments in $y$-coordinates to conclude the existence of a gauge transformations $U : \textbf{a} \to \textbf{b} = U \cdot \textbf{a}$, such that, in the transformed gauge $\textbf{b},$ the connection $\A$ has optimal regularity, $\A\in W^{1,p}$.  Since the gauge transformation $U$ does not affect the connection $\Gamma$, Theorem \ref{Thm_opt} extend in a straightforward way to arbitrary non-optimal affine connections $\Gamma$ on the base manifold $\M$ which satisfy the bound (\ref{bound_1}). Taken on whole, this establishes Theorem \ref{Thm_opt_general}.
\QED

For Uhlenbeck compactness we consider the case Local (ii) with $p=\infty$, the case addressed in \cite{ReintjesTemple_Uhl1, ReintjesTemple_ell6}. We state here a concise version of Uhlenbeck compactness, c.f. Theorem \ref{Thm_compactness} here and Theorem 3.2 in \cite{ReintjesTemple_Uhl1}.

\begin{Thm}\label{Thm_compactness_general}
Let $(\Gamma_i,\A_i) \in L^\infty$ be a sequence of connections on ${\mathcal V}\mathcal{M}$ uniformly bounded by
\beq \label{bound_4}
\| (\Gamma_i,\A_i) \|_{L^{\infty}} + \|(d\Gamma_i,d\A_i) \|_{L^p} \; \leq \; M
\eeq
in a given coordinate system and gauge $(x,\textbf{a})$,  $p>n$, (c.f. \eqref{bound_1}). Then, locally, under coordinate and $SO(r,s)$-gauge transformation to optimal regularity, a subsequence of the transformed connections $(\Gamma_i,\A_i)$ converges  weakly in $W^{1,p}$, strongly in $L^p$, to a limit connection $(\Gamma,\A)$. 
\end{Thm}

\Proof
Given a sequence of connections $(\Gamma_i,\A_i)$ satisfying the uniform bound \eqref{bound_4}, Theorem \ref{Thm_compactness_general} follows by combining \cite[Theorem 2.2]{ReintjesTemple_Uhl1} for affine connections to obtain compactness of $\Gamma_i$, followed by the compactness asserted by Theorem \ref{Thm_compactness} of this paper, to obtain compactness of the fibre part $\A_i$.  Putting these results together implies compactness of the sequence of coupled pairs $(\Gamma_i,\A_i)$, i.e. weak $W^{1,p}$ convergence of a subsequence of $(\Gamma_i,\A_i)$.   More precisely, it follows that after extraction of a subsequence, the coordinate transformations $y_i(x)$ to optimal regularity converge weakly in $W^{2,2p}$ to some coordinate transformation $y(x)$, (as a result of \eqref{bound_3} - \eqref{bound_3b}), and the connection components $(\Gamma_{y_i},\A_{(y_i,\textbf{b})})$ in $y_i$-coordinates converge weakly in $W^{1,p}$ to some connection $(\Gamma_y,\A_{(y,\textbf{b})})$, as $i\to \infty$, (as a result of (\ref{bound_2})), when expressed as functions of the original $x$-coordinates. 
\QED

\section{Preliminaries} \label{Sec_Prelim}

We now introduce the Cartan Calculus for matrix valued differential forms required to formulate the RT-equations on vector bundles, analogous to what we introduced in \cite{ReintjesTemple_ell1} for tangent bundles. We again consider the trivialization of a vector bundle $\VM \simeq \R^N \times \Omega$ with $N$-dimensional fibres and base space $\Omega \subset \R^n$ open and bounded.  We continue to assume fixed Cartesian coordinates $x$ on $\Omega$.  By a matrix valued differential $k$-form $\omega$ on $\VM$ we mean a $k$-form over $\Omega$ with $(N\times N)$-matrix components, 
\beq \label{def_matrixvalued_diff-form}
\omega = \omega_{[i_1...i_k]} dx^{i_1} \wedge ... \wedge dx^{i_k} \equiv \sum_{i_1< ... < i_k} \omega_{i_1...i_k} dx^{i_1} \wedge ... \wedge dx^{i_k},
\eeq 
for $(N\times N)$-matrices  $\omega_{i_1...i_k}$  such that total anti-symmetry holds in the indices $i_1,...,i_k \in \{1,...,n\}$, and we follow Einstein's convention of summing over repeated upper and lower indices from $1$ to $n$. (Note, we never ``raise'' or ``lower'' indices.) We define the wedge product of a matrix valued $k$-form $\omega$ with a matrix valued $l$-form $u = u_{j_1...j_l} dx^{j_1} \wedge ... \wedge dx^{j_l}$ as  
\begin{eqnarray} \label{def_wedge}
\omega \wedge u  
&\equiv & \frac{1}{l!k!} \omega_{i_1...i_k} \cdot u_{j_1...j_l} \; dx^{i_1} \wedge ... \wedge dx^{i_k} \wedge dx^{j_1} \wedge ... \wedge dx^{j_l}, 
\end{eqnarray}
where ``$\cdot$'' denotes matrix multiplication, (so $\omega\wedge \omega \neq 0$ is possible). 

The exterior derivative $d$ is defined component-wise on matrix components,
\beq \label{exterior_derivative}
d \omega \equiv \partial_l \omega_{[i_1...i_k]} dx^l \wedge dx^{i_1} \wedge ... \wedge dx^{i_k} ,
\eeq
and we define the co-derivative $\delta$ on a matrix valued $k$-form $\omega$ as 
$$
\delta \omega  \equiv (-1)^{(k+1)(n-k)} *d* \omega,
$$ 
where $*$ is the Hodge star introduced in terms of the Euclidean metric in $x$-coordinates on $\Omega \subset \R^n$. Both $d$ and $\delta$ act component-wise on matrix components, and all properties of $d$ and $\delta$ for scalar valued differential forms carry over to matrix valued forms. Note, $d$ requires no metric, but $\delta$ requires a metric for its definition and we choose the Euclidean metric in $x$-coordinates. As a result, the Laplacian $\Delta \equiv d \delta + \delta d$ is given by the standard Euclidean Laplacian
$$
\Delta = \partial^2_{x^1} + ... + \partial^2_{x^n},
$$ 
is hence elliptic, and acts component-wise on matrix components and differential form components.

We now introduce some formulas which we require in the subsequent analysis, (the proofs in the case of tangent bundles in \cite[Sec. 3]{ReintjesTemple_ell1} carry over directly by replacing $n$ by $N$). The exterior derivative satisfies the product rule  
\beq \label{Leibnitz-rule-d}
d(\omega\wedge u) = d\omega \wedge u + (-1)^k \omega \wedge du,
\eeq
where $\omega$ is a matrix valued $k$-form and $u$ is a matrix valued $j$-form. Assuming $\omega, u\in W^{1,p}(\Omega)$, $p>n$, the right hand side of \eqref{Leibnitz-rule-d} lies again in $W^{1,p}(\Omega)$ by Morrey's inequality. By $\omega \in W^{1,p}(\Omega)$ we mean its components are in $W^{1,p}(\Omega)$. For invertible matrix valued $0$-forms $U \in W^{1,p}(\Omega)$, using $dU^{-1} = - U^{-1}\cdot dU\cdot U^{-1}$, \eqref{Leibnitz-rule-d} implies    
\beq  \label{Leibniz_rule_J-application} 
d\big( U^{-1} \cdot dU \big) =  d(U^{-1}) \wedge  dU  = -  U^{-1} d U \wedge U^{-1} dU\in W^{1,p}(\Omega).
\eeq
We define the matrix valued ``inner product'' $\langle \omega\; ; u \rangle $ on the space of matrix valued $k$-forms by
\beq \label{def_inner-product}
\langle \omega\; ; u \rangle^\mu_\nu \equiv \sum_{\sigma=1}^n \sum_{i_1<...<i_k} \omega^\mu_{\sigma\: i_1...i_k} u^\sigma_{\nu\: i_1...i_k}.
\eeq
Note, $\langle \omega\; ; u \rangle$ converts two matrix valued $k$-forms into a matrix valued $0$-form, and turns into the standard inner product on scalar forms induced by the Euclidean metric. Under multiplication by matrix valued $0$-forms $U$, $\langle \omega\; ; u \rangle$ satisfies
\beq \label{inner-product_muliplications}
 U \cdot \langle \omega\; ; u \rangle =  \langle  U \cdot \omega\; ; u \rangle,  \hspace{.3cm} 
\langle \omega\cdot U \; ; u \rangle =  \langle \omega\; ;  U \cdot u \rangle, \hspace{.3cm}
\langle \omega \; ; u \cdot U \rangle =  \langle \omega\; ; u \rangle \cdot U  .
\eeq
The co-derivative $\delta$ satisfies the product rule
\beq \label{Leibnitz-rule-delta}
\delta (U \mm w ) = U \mm \delta w  + \langle d U ;  w \rangle, 
\eeq
where $U\in W^{1,p}(\Omega)$ is a matrix valued $0$-form and $w \in W^{1,p}(\Omega)$ a matrix valued $1$-form, and if $p>n$, then the right hand side of \eqref{Leibnitz-rule-delta} lies in $W^{1,p}(\Omega)$.

We finally define the $L^2$-inner product on matrix valued differential forms by                   
\beq \label{inner-product_L2}
\langle \omega, u\rangle_{L^2} 
\equiv  \int_\Omega  {\rm tr} \langle \omega^T ; u \rangle dx  \   
=  \ \sum_{\nu, \sigma=1}^N  \sum_{i_1<...<i_k}  \int_\Omega \omega^\nu_{\sigma\: i_1...i_k} u^\nu_{\sigma\: i_1...i_k} dx, 
\eeq
where $dx$ is the Lebesgue measure on $\R^n$ and ${\rm tr}$ denotes the matrix trace. By the method of proof of Lemma 8.1 in \cite{ReintjesTemple_Uhl1}, one can easily show the following integration by parts formula,
\beq \label{partial_integration_matrix}
\langle d u, \omega \rangle_{L^2} + \langle  u, \delta \omega \rangle_{L^2} =0,
\eeq
which applies to matrix valued $k$-forms $u \in W^{1,p}(\Omega)$ and matrix valued $(k+1)$-forms $\omega \in W^{1,p^*}(\Omega)$, $k\geq 0$ and $\frac{1}{p} + \frac{1}{p^*} =1$, so long as either $u$ or $\omega$ vanishes on the boundary $\partial\Omega$. Here, as in \cite{ReintjesTemple_Uhl1}, all Sobolev norms are taken component-wise on matrix valued differential forms, and we use the standard notation $\|\cdot \|_{m,p} \equiv \|\cdot\|_{W^{m,p}(\Omega)}$ for $m\geq -1$.

\section{The RT-equations for vector bundles}  \label{Sec_RT-eqn}

The vector bundle version of the RT-equations \eqref{RT_1} - \eqref{RT_2} is the basis of our proofs of Theorems \ref{Thm_equiv} - \ref{Thm_compactness}. Recall that we fix an arbitrary gauge $\textbf{a}$, and assume the connection components $\A \equiv \A_\textbf{a}$ of $\A_{\VM}$ in this gauge are non-optimal. In this section we derive the RT-equations \eqref{RT_1} - \eqref{RT_2} from the connection transformation law, starting with the assumption that a gauge transformation $U$ exists which maps the connection $\A$ to optimal regularity. Conversely, we then explain how gauge transformations $U$ which solve the reduced RT-equation \eqref{RT_2} do transform connections to optimal regularity, and we explain why such solutions $U$ always lie in $SO(r,s)$. The purpose of this section is to present the logic of our argument in its essence, assuming the case of smooth connections; for example, take $\A \in W^{1,2p}(\Omega)$ and $d\A\in W^{1,p}(\Omega)$, $U \in W^{2,2p}(\Omega)$, and optimal regularity $\Aop \in W^{2,p}(\Omega)$. Keep in mind that the goal of this paper is to treat the low regularity case $\A  \in L^{2p}(\Omega)$ with $d\A  \in L^{p}(\Omega)$, $U \in W^{1,2p}(\Omega)$, and $\Aop \in W^{1,p}(\Omega)$.  The derivation of the RT-equations accomplished here gives the correct final form \eqref{RT_1} - \eqref{RT_2}, but the RT-equations \eqref{RT_1} - \eqref{RT_2} must be interpreted in the weak sense to handle the case of low regularities in Theorems \ref{Thm_equiv} - \ref{Thm_compactness}. The argument in Sections \ref{Sec_RT-eqn_optimal} and \ref{Sec_SO(r,s)}, establishing that solutions of the RT-equations give optimal regularity for smooth connections, must be modified to incorporate the weak formalism of the equations. This is accomplished in Sections \ref{Sec_proof_existence} - \ref{Sec_proofs_rest}.

\subsection{Derivation of the RT-equations} \label{Sec_RT-eqn_derivation}

Our strategy for deriving the RT-equations associated to connections on vector bundles parallels the one for affine connections in \cite[Ch. 3 - 5]{ReintjesTemple_ell1}, (see also \cite[Ch. 5]{ReintjesTemple_ell3}).  To begin, assume there exists a gauge transformation $U$ in $SO(r,s)$ which maps $\A \equiv \A_\textbf{a}$ to a connection $\Aop$ which has optimal regularity, $\textbf{b} = U \cdot \textbf{a}$. Define the matrix valued $1$-form $\Ati$ by
\beq \label{Gammati_VB}
\Ati \equiv U^{-1} \Aop U,
\eeq
so the  connection transformation law \eqref{connection_transfo_VB} reads
\beq \label{connection_transfo_VB2}
\A = U^{-1} dU + \Ati.
\eeq
Starting from the transformation law \eqref{connection_transfo_VB2}, we now derive solvable elliptic equations in $U$ and $\Ati$. Taking the exterior derivative $d$ of \eqref{connection_transfo_VB2} yields
\beq  \label{der_eqn1}
d\Ati = d\A - dU^{-1} \wedge dU,
\eeq  
where the last term follows by the Leibnitz rule \eqref{Leibniz_rule_J-application}. On the other hand, multiplying \eqref{connection_transfo_VB2} by $U$ and taking the co-derivative yields
\beq \label{der_eqn2}
\Delta U = U \cdot \big( \delta \A - \delta \Ati \big) + \langle dU ; \A - \Ati \rangle,
\eeq
where we used that $\Delta U = \delta d U$ for the matrix valued $0$-form $U$, (recall $\delta U=0$ for all $0$-forms), and we used the Leibnitz-rule \eqref{Leibnitz-rule-delta} for $\delta$ to derive the right hand side in terms of the matrix valued inner product $\langle \cdot\: ; \cdot \rangle$ defined in \eqref{def_inner-product}. Observe now that the connection transformation law \eqref{connection_transfo_VB2} leaves the co-derivative $\delta \Ati$ undetermined, so we are free to choose a matrix valued $0$-form $\alpha$ (at least as regular as $\A$, $d\A$) and set
\beq \label{der_eqn3}
\delta \Ati = U^{-1} \alpha.
\eeq
System \eqref{der_eqn1} and \eqref{der_eqn3} is of Cauchy-Riemann type and would in principal determine $\Ati$, but in analogy to \cite{ReintjesTemple_ell1}, we prefer to write this system as the Poisson type equation
\beq \label{RT_VB_1}
\Delta \Ati = \delta d \A - \delta \big( dU^{-1} \wedge dU \big) + d(U^{-1} \alpha ), 
\eeq
which results from taking $d$ of \eqref{der_eqn3} and $\delta$ of \eqref{der_eqn1}, adding the resulting equations, and using $\Delta = d \delta + \delta d$. Moreover, substituting \eqref{der_eqn3} for $\delta \Ati$ in \eqref{der_eqn2}, we obtain
\beq  \label{RT_VB_2} 
\Delta U =   U \delta\A   + \langle dU ; \A - \Ati \rangle - \alpha .
\eeq
The coupled system formed by \eqref{RT_VB_1} and \eqref{RT_VB_2} is what we take to be the preliminary version of the RT-equations.\footnote{Note, there appears to be a loss of information in going from the first order system \eqref{der_eqn1}, \eqref{der_eqn3} to the Poisson type equation \eqref{RT_VB_1}, but interestingly, the final form of the RT-equations suffices to establish optimal regularity, c.f. Lemma \ref{Lemma_Gammati'} below.} At this stage $\alpha$ is some matrix valued $0$-form which is free to be chosen.  From \eqref{RT_VB_1} and \eqref{RT_VB_2} we next derive the final version \eqref{RT_1} - \eqref{RT_2} of the RT-equations, which, surprisingly, decouple. 

To continue the derivation, note that solutions of \eqref{RT_VB_1} - \eqref{RT_VB_2} could in general allow for arbitrary matrix valued solutions. The critical step now is to obtain solutions $(U,\Ati)$ of  \eqref{RT_VB_1} - \eqref{RT_VB_2} such that the gauge transformation $U(x)$ lie in $SO(r,s)$ for every $x\in \Omega$, i.e.,  ${\rm det}(U)=1$ and $U^T \cdot \eta \cdot U = \eta$. The condition ${\rm det}(U)=1$ is met by solving the RT-equation for $U$ close to the identity $I$, and can hence be assumed for the rest of this section. To arrange for $U^T \eta U = \eta$, we now derive an equation for $\alpha$. To begin, define
\beq \label{w_first}
w \equiv  U^T\eta  U - \eta ,
\eeq
so $w=0$ is equivalent to $U \in SO(r,s)$. We now assume $w=0$, so $\Delta w=0$.  Then substituting \eqref{w_first} into $\Delta w=0$ and applying the Leibnitz rule gives
\begin{eqnarray} \label{der_eqn4}
0 = \Delta w = \Delta \big( U^T \eta U \big) =  (\Delta U)^T \eta U + 2 \langle dU^T;\eta dU \rangle + U^T \eta\: \Delta U,
\end{eqnarray}
and substituting \eqref{RT_VB_2} for $\Delta U$, we obtain
\begin{align}
\alpha^T \cdot \eta U + U^T\eta \cdot \alpha \ =\ \  &   
(\delta \A )^T \cdot U^T \eta U + U^T \eta U \cdot  \delta \A      + 2 \langle dU^T ;\eta dU\rangle   \cr  
& +\langle dU ; \A - \Ati \rangle^T \cdot \eta U + U^T \eta \cdot \langle dU ; \A - \Ati \rangle  .
\end{align}
(Note that $(\Delta U)^T = \Delta U^T$ and $(\delta \A)^T = \delta\A^T$.) Using now that $\delta \A$ lies pointwise in the Lie algebra of $SO(r,s)$, (the linear space of anti-symmetric matrices $X$ with respect to $\eta$, i.e. $X^T \eta + \eta X=0$), together with $U^T \eta U=\eta$, we obtain the cancellation
\beq \label{commutator}
\delta \A^T \cdot U^T \eta U + U^T \eta U \cdot  \delta \A  =0.
\eeq 
This cancellation is crucial because $\delta\A$ is at a regularity too low for the RT-equations to imply optimal regularity. From \eqref{commutator} we conclude
\beq \label{alpha-eqn}
\alpha^T \cdot \eta U + U^T\eta \cdot \alpha =   
2 \langle dU^T ;\eta  dU\rangle 
+ \big(U^T\eta \langle dU ; \A - \Ati \rangle \big)^T 
+  U^T \eta \langle  dU ; \A - \Ati \rangle , 
\eeq
which is our sought after equation for $\alpha$.  Since $\langle dU^T ; \eta dU\rangle^T = \langle dU^T ; \eta dU\rangle$ and $\eta^T~=~\eta$, the solutions of \eqref{alpha-eqn} are given by    
\beq \label{alpha_soln}
\alpha \equiv  (U^T \eta)^{-1} \langle dU^T ; \eta dU\rangle + \langle dU ; \A - \Ati \rangle  +  U \X  ,
\eeq
where $\X$ can be any element in the Lie algebra with the same regularity as $\A$. 

Substituting the expression \eqref{alpha_soln} for $\alpha$ back into \eqref{RT_VB_2}, we obtain
\beq \label{RT_der_1}
\Delta U =   U \delta\A   -  (U^T \eta)^{-1} \langle dU^T ; \eta dU\rangle -  U \X   ,
\eeq
the sought after RT-equation \eqref{RT_2}.  The Lie algebra element $\X$ is free to be chosen and for our purposes here one can set $\X=0$; in particular, the choice of $\alpha$ is not unique. Interestingly, for $U \in SO(r,s)$ we have $(U^T \eta)^{-1} = U \eta$, but to deduce that $U \in SO(r,s)$ we found that taking $(U^T \eta)^{-1}$ instead of $U \eta$ in \eqref{RT_der_1} makes the proof more feasible.   In Section \ref{Sec_SO(r,s)} below we show that solutions $U$ of the RT-equations \eqref{RT_der_1} always lie in $SO(r,s)$, so long that they lie in $SO(r,s)$ on the boundary of $\Omega$.

To complete the derivation of the $\Ati$-equation \eqref{RT_1}, anticipating that solutions $U$ of \eqref{RT_der_1} will lie in $SO(r,s)$, we assume $U^T\eta U = \eta$ in the form $U^{-1} = \eta U^T \eta$, and we compute for the first term in \eqref{alpha_soln} that
\begin{eqnarray}   \nonumber
(U^T \eta)^{-1} \langle dU^T ; \eta dU\rangle 
 &=&  U \eta \langle d(U^T \eta) ;  dU\rangle  
=   \langle U d(\eta U^T \eta) ;  dU\rangle    \cr
&=&  \langle U d(U^{-1}) ;  dU\rangle  
= \langle  dU ; U^{-1} dU\rangle,
\end{eqnarray} 
where we used that $dU^{-1} = U^{-1} dU U^{-1}$ and the multiplication property \eqref{inner-product_muliplications} of the matrix valued inner product. Thus, $\alpha$ in \eqref{alpha_soln} simplifies to
\beq 
\alpha = -  \langle dU ; U^{-1} dU -(\A - \Ati) \rangle  +  U \X  ,
\eeq
and using that $dU-U\cdot( \A - \Ati)=0$ by the connection transformation law \eqref{connection_transfo_VB2}, we get
\beq \label{alpha-zero}
\alpha =  U \X  .
\eeq
Substituting now \eqref{alpha-zero} into \eqref{RT_VB_1} leads to the sought after RT-equation
\beq \label{RT_der_2}
\Delta \Ati = \delta d \A - \delta \big(dU^{-1} \wedge dU \big) + d\X  .
\eeq
Again we set $\X =0$, which is sufficient for our theory.

Note finally that the expressions for $\alpha$ in \eqref{alpha_soln} and \eqref{alpha-zero} are equivalent as long that $U \in SO(r,s)$. The choice of the more complicated expression \eqref{alpha-zero} for $\alpha$ in \eqref{RT_der_1}, turned out to be beneficial for our proof that solutions $U$ of \eqref{RT_der_1} lie indeed in $SO(r,s)$. Moreover, \eqref{alpha-zero} is required for the decoupling of \eqref{RT_VB_1} and \eqref{RT_VB_2}, which would not occur for $\alpha=0$. Since \eqref{RT_der_1} and \eqref{RT_der_2} decouple, given a solution $U$ in $SO(r,s)$ of \eqref{RT_der_1}, one may equivalently take \eqref{alpha_soln} or \eqref{alpha-zero} in the derivation of the first RT-equation \eqref{RT_der_2}, essentially amounting to different choices of $X$.

\subsection{How the RT-equations yield optimal regularity}  \label{Sec_RT-eqn_optimal}

The logic of how to obtain optimal regularity from the RT-equations is encapsulated in \eqref{def_Ati'} - \eqref{OptimalGamma}. Given a solution $U \in SO(r,s)$ to the reduced RT-equation \eqref{RT_2}, the connection transformation law \eqref{connection_transfo_VB} determines the connection $\A_\textbf{b}$ in the gauge $\textbf{b} = U \cdot \textbf{a}$. The problem then is to show that $\A_\textbf{b}$ has optimal regularity. For this, we define 
\beq \label{Ati_Sec4}
\Ati' \equiv \A - U^{-1} dU,
\eeq
then the connection $\A_\textbf{b}$ is given by
\beq
\Aop = U \Ati' U^{-1}.
\eeq 
Thus optimal regularity for $\A_\textbf{b}$ follows from optimal regularity for $\Ati'$, which we obtain by showing that $\Ati'$ is an exact solution of the first RT-equation \eqref{RT_1}, as stated in the following lemma.

\begin{Lemma} \label{Lemma_Gammati'}
Assume $U$ solves the reduced RT-equation \eqref{RT_2} such that $U(x) \in SO(r,s)$ for every $x\in \Omega$, then $\Ati'$ defined in \eqref{Ati_Sec4} solves the first RT-equation \eqref{RT_1} with $\Ati$ replaced by $\Ati'$.
\end{Lemma}

\Proof
Applying $\Delta= d\delta + \delta d$ to \eqref{Ati_Sec4} and using the product rule \eqref{Leibniz_rule_J-application} gives
\begin{eqnarray}  \label{Gammati-eqn_eqn1}
\Delta \Ati' &=& \Delta \A - \Delta(U^{-1} dU)   \cr
&=& \delta d \A - \delta (dU^{-1} \wedge dU) 
+ d\Big( \delta \A - \delta (U^{-1}dU) \Big).
\end{eqnarray}
Assuming $U$ satisfies $U^T \eta U=\eta$ and using $\eta^{-1}=\eta$, we have $U^{-1} = \eta U^T \eta$. Employing this identity, the product rule \eqref{Leibnitz-rule-delta} gives
\beq
\delta (U^{-1}dU) =  \eta^{-1} \langle dU^T ; \eta dU \rangle + \eta^{-1} U^T \eta \cdot  \Delta U,
\eeq 
where we used that $\Delta U \equiv d\delta U$ because $\delta U=0$ for $0$-forms. Then, substituting the RT-equation \eqref{RT_2} for $\Delta U$, i.e. $\Delta U =   U \delta\A   -  (U^T \eta )^{-1}  \langle dU^T ; \eta dU\rangle$, and using $U^T\eta U =\eta$, we find that the last term in \eqref{Gammati-eqn_eqn1} vanishes,\footnote{This resembles the Coulomb gauge, but it is different, since $\Ati'= U^{-1} \A_\textbf{b} U$ is not a connection.}
\begin{eqnarray}  \label{Gammati-eqn_eqn2}
\delta \Ati' = \delta \A - \delta (U^{-1}dU) 
=  0.
\end{eqnarray}
Substitution of \eqref{Gammati-eqn_eqn2} into \eqref{Gammati-eqn_eqn1} then gives
\begin{eqnarray}
\Delta \Ati' &=&  \delta d \A - \delta (dU^{-1} \wedge dU)  ,
\end{eqnarray}
which is the sought after RT-equation \eqref{RT_1}. This completes the proof.
\QED

We prove in Section \ref{Sec_proof_equiv} below, that elliptic regularity theory implies that $\Ati'$ is two derivatives more regular than the right hand side of \eqref{RT_1}, which lies in $W^{-1,p}(\Omega)$ for non-optimal connections when $\A \in L^{2p}(\Omega)$ with $d\A \in L^p(\Omega)$. This gives, after establishing the above lemma at the weak level to account for the low regularity in Section \ref{Sec_weak}, the sought after optimal regularity $W^{1,p}$ for the connection $\Aop = U \Ati' U^{-1}$ in the new gauge $\textbf{b} = U \cdot \textbf{a}$. This is the basis for the proofs of Theorems \ref{Thm_equiv} and \ref{Thm_opt}, once existence of solutions to the RT-equations is proven.

\subsection{How the RT-equations yield solutions in the gauge group}   \label{Sec_SO(r,s)}

We now explain how to prove that solutions $U$ to the reduced RT-equation \eqref{RT_2} always lie in $SO(r,s)$, so long as they lie in $SO(r,s)$ on the boundary of $\Omega$. The argument for compact groups $SO(N)\equiv SO(N,0)$ is significantly simpler than the one for non-compact groups $SO(r,s)$, $r >0$, $s > 0$.

To start recall that $U \in SO(r,s)$ is equivalent to $w=0$, for $w \equiv U^T \eta U-\eta$. We now derive an equation which $w$ satisfies whenever $U$ is a solution of the reduced RT-equation \eqref{RT_2}. As in \eqref{der_eqn4}, the Leibniz rule gives 
\beq
\Delta w =  (\Delta U)^T \eta U + 2 \langle dU^T;\eta dU \rangle + U^T \eta\: \Delta U,
\eeq
and substituting for $\Delta U$ and $(\Delta U)^T$ the right hand side of \eqref{RT_2} yields
\beq \label{SO(n)_eqn1}
\Delta w = \delta \A^T \cdot U^T \eta U + U^T \eta U \cdot  \delta \A ,
\eeq
where we have used the symmetry $\langle dU^T ; \eta dU\rangle^T = \langle dU^T ; \eta dU\rangle$.  Subtracting the Lie algebra condition $\delta \A^T \cdot \eta + \eta \cdot \delta \A =0$ from the right hand side of \eqref{SO(n)_eqn1} gives,
\beq \label{SO(n)_eqn2}
\Delta w = \delta \A^T \cdot w + w \cdot  \delta \A  ,
\eeq
which is our sought after equation for $w$. Equation \eqref{SO(n)_eqn2} is a linear system of elliptic PDE's. Assuming $w$ vanishes on the boundary $\partial \Omega$, we conclude that $w=0$ is a solution of \eqref{SO(n)_eqn2} in $\Omega$. But to prove that $U \in SO(N)$ or $U \in SO(r,s)$, we have to show that $w=0$ is the only solution to \eqref{SO(n)_eqn2} with zero boundary data, a non-trivial problem in view of the Fredholm alternative.  We establish the desired uniqueness in the case of the compact group $SO(N)$ in the following lemma, which is based on a cancellation in the bilinear form associated with \eqref{SO(n)_eqn2} when $\eta=I$. A more subtle argument, which works for both compact and non-compact groups $SO(r,s)$, ($r,s\geq 0$), and on which our proofs of Theorems \ref{Thm_existence} - \ref{Thm_compactness} are based, is developed below the Lemma to establish $w=0$. We complete this argument in Section \ref{Sec_proof_exist_SO(r,s)} at the low regularity of $L^p$ connections. 
The following lemma, which works only in the compact case, shows that the compact case is simpler than the non-compact case.\footnote{We find this difference interesting because the method in \cite{Uhlenbeck} appears to require compactness of the gauge group.}

\begin{Lemma}  \label{Lemma_SO(n)}
Assume $U$ is a solution of the reduced RT-equation \eqref{RT_2} satisfying the condition that $U$ lies in $SO(N)$ on the boundary $\partial\Omega$. Then $U(x) \in SO(N)$ for any $x\in \Omega$.              
\end{Lemma}

\Proof     
Let $\eta=I$, so $SO(r,s)=SO(N)$. To prove that a solution $U$ of \eqref{RT_1} lies in $SO(N)$ whenever it lies in $SO(N)$ on $\partial\Omega$, we now show that $w \equiv U^T  U-I=0$ in $\Omega$. Since $w$ satisfies \eqref{SO(n)_eqn2}, it suffices to prove that any symmetric solution $w \in W^{1,2p}(\Omega)$ of \eqref{SO(n)_eqn2} which vanishes on $\partial\Omega$, also vanishes inside $\Omega$. For this, we define the bilinear functional associated to \eqref{SO(n)_eqn2}
\beq \label{SO(n)_eqn3}
B(w,v) \equiv \int_\Omega {\rm tr} \langle dw^T ; dv \rangle dx  + \int_\Omega {\rm tr}\big( w^T \delta\A \: v + (w \delta \A)^T v\big) dx ,
\eeq
on matrix valued $0$-forms $w \in W^{1,2p}(\Omega)$ and $v \in W^{1,(2p)^*}(\Omega)$. Integration by parts shows that $B$ vanishes on solutions $w$ of \eqref{SO(n)_eqn2} for any test form $v \in W^{1,(2p)^*}(\Omega)$. Thus, by cyclic commutativity of the trace, any symmetric solution $w \in W^{1,2p}(\Omega)$ of \eqref{SO(n)_eqn2} satisfies
\begin{eqnarray} \label{SO(n)_eqn3b}
0 = B(w,w) 
= \|dw\|_{L^2}^2 + \int_\Omega {\rm tr}\big(  (\delta\A + \delta\A^T)w^2\big) dx 
= \|dw\|_{L^2}^2  ,
\end{eqnarray}
where we used that 
\beq \label{SO(n)_eqn3c}
\delta\A+\delta\A^T=0
\eeq 
for the last inequality, and note that $W^{1,2p}(\Omega) \subset H^1(\Omega) \subset W^{1,(2p)^*}(\Omega)$. Applying next the Poincar\'e inequality, using that $w=0$ on $\partial\Omega$, gives 
\beq 
0= B(w,w) \geq C \| w \|_{H^1}^2,
\eeq
which implies that any symmetric $w$ solving \eqref{SO(n)_eqn2} is the trivial solution $w=0$. This completes the proof.
\QED

\noindent Note that for the non-compact case, the cancellation in \eqref{SO(n)_eqn3c} does not take place, since $\delta\A+\delta\A^T= \delta\A - \eta \cdot \delta \A \cdot \eta$ is nonzero in general when $\eta \neq I$. 

Thus we now develop a more subtle argument, different from the one in Lemma \ref{Lemma_SO(n)}, which is sufficient to handle the case of both compact and non-compact groups $SO(r,s)$, ($r,s\geq 0$). The complication in the non-compact case over the compact one is that, in view of the Fredholm alternative, the spectrum associated with equation \eqref{SO(n)_eqn2} could be non-trivial and could allow for non-zero solutions. To prove that solutions $U$ of \eqref{RT_2} lie pointwise in $SO(r,s)$ everywhere in $\Omega$, assuming $U\in SO(r,s)$ on $\partial\Omega$, for a fixed $\A$, we replace $\A$ in \eqref{RT_2} and \eqref{SO(n)_eqn2} by $\lambda \A$ for $\lambda$ in a neighborhood of $1$, say, $\lambda \in (0,1]$; we then show below that for almost every $\lambda \in (0,1],$ the solution $U^\lambda$ of \eqref{RT_2} lies in $SO(r,s)$ because equation \eqref{SO(n)_eqn2} for $w^\lambda$ has no non-trivial eigenvectors.  We extend this to every $\lambda$ by proving in Section \ref{Sec_proof_exist_SO(r,s)} that our iteration scheme converges {\it uniformly} in $\lambda$. This implies continuity of  $U^\lambda$ with respect to $\lambda$, which in turn implies the continuity of $w^\lambda$ with respect to $\lambda$ as well. Since $w^{\lambda}=0$ for almost every $\lambda$, continuity implies $w\equiv0$, and this completes the proof in Section \ref{Sec_proof_exist_SO(r,s)} that $U^\lambda$ lies in $SO(r,s)$ for all $\lambda \in (0,1]$.  

More precisely, the key observation in our argument is that \eqref{SO(n)_eqn2} can be written as an eigenvalue problem when replacing $\A$ by $\lambda\A$ as follows: To begin, write the right hand side of \eqref{SO(n)_eqn2} as 
\beq
M(w) \equiv \delta \A^T \cdot w + w \cdot  \delta \A ,  
\eeq
which is a linear mapping in $w$ pointwise at each $x\in \Omega$. We then write \eqref{SO(n)_eqn2} as
\beq \label{SO(n)_eqn4}
 w = \Delta^{-1} M(w).
\eeq
We prove in Section \ref{Sec_weak} that for $\A \in L^{2p}(\Omega)$, $M: W^{1,2p}\to W^{-1,2p}$ is a bounded linear operator, from which we then conclude that 
\beq
K \equiv \Delta^{-1} M: \  W^{1,2p} \to W^{1,2p}
\eeq 
is a bounded and compact operator, and therefore has (at most) a countable spectrum $\Sigma$, see \cite[Thm 4.25]{Rudin}. 

Now for each fixed $\A,$ consider the family of solutions $U_\lambda$ of the RT-equation \eqref{RT_2} with $\A$ replaced by $\lambda \A$ for $\lambda \in (0,1]$ such that $U_\lambda$ lies in $SO(r,s)$ on the boundary of $\Omega$ for each $\lambda \in (0,1]$. For such a $U_\lambda$, the eigenvalue problem \eqref{SO(n)_eqn4} turns into 
\beq \label{SO(n)_eqn5}
\frac{1}{\lambda}  w_\lambda = K(w_\lambda),
\eeq
where
\beq \label{SO(n)_eqnNew} 
w_\lambda \equiv U_\lambda^T \eta U_\lambda - \eta.
\eeq 
But since $K\equiv \Delta^{-1} M$ is a compact operator, it has a countable spectrum, so the solution $w^{\lambda}$ in (\ref{SO(n)_eqnNew}) must be zero for almost every $\lambda$.   To guarantee the solution $U$ satisfies $w=0$ at $\lambda=1$ when $\lambda=1$ is in the spectrum of $K\equiv \Delta^{-1} M,$ we prove in Section \ref{Sec_proof_exist_SO(r,s)} that the uniform convergence of the iteration scheme implies that $U^\lambda$ is a continuous function of $\lambda$ in $W^{1,2p}$. Thus, since $w^\lambda =0$ for almost every $\lambda$, it follows by continuity that $w^\lambda =0$ for all $\lambda$ including $\lambda=1$. The argument here will be incorporated into our existence theory in Section \ref{Sec_proof_existence} in order to prove rigorously that $U$ lies in $SO(r,s)$ at the low regularities of Theorems \ref{Thm_existence} - \ref{Thm_compactness}. 

Let us finally remark that the above derivation of equation \eqref{SO(n)_eqn2} for $w$ as well as our proof of Lemma \ref{Ati_Sec4} are both based on a pointwise substitution of the reduced RT-equation \eqref{RT_2}. This is not possible at the low regularity of $L^p$ connections addressed in Theorems \ref{Thm_equiv} - \ref{Thm_compactness}, which requires a weak formulation of equation \eqref{RT_2} and a corresponding substitution of functionals in a suitable sense, both established in Section \ref{Sec_weak}. In Section \ref{Sec_proof_existence} we prove the required existence of weak solutions of the reduced RT-equations \eqref{RT_2}, which is based on a weak formulation and an iteration scheme which reduces the existence problem to that of the classical weak Poisson equations on each component.

\section{Existence theory for the RT-equations  - Proof of Theorem \ref{Thm_existence}}  \label{Sec_proof_existence}

We now prove our existence result, Theorem \ref{Thm_existence}, for the reduced RT-equations \eqref{RT_2}, a non-linear elliptic system. For this, as in Theorem \ref{Thm_existence}, we assume $\A \in L^{2p}(\Omega)$, $p \in (n/2,\infty)$, and let $M>0$ be a constant such that
\beq \label{initial_bound_recall}
 \|\A\|_{L^{2p}(\Omega)}  \; \leq \; M.
\eeq  
Theorem \ref{Thm_existence} then establishes, (in the local case), that for any point $q\in \Omega$ there exists a neighborhood $\Omega' \subset \Omega$, and there exists a solution $U \in W^{1,2p}(\Omega')$ of the reduced RT-equations \eqref{RT_2}, such that $U(x) \in SO(r,s)$ for every $x\in\Omega'$, and such that 
$U$ satisfies  
\beq \label{uniform_bound_U_recall}
\|U\|_{W^{1,2p}(\Omega')}  \; \leq\; C(M)\; \|\A\|_{L^{2p}(\Omega')} ,
\eeq
for some constant $C(M) > 0$, depending only on $\Omega, p, n$ and $M$. The neighborhood $\Omega' \subset \R^n$ depends on $\Omega, p, n, M$ as well as the $L^p$-norm of $\A$ near $q$. If \eqref{initial_bound_recall} holds for $p=\infty$, then \eqref{uniform_bound_U_recall} holds for any $p<\infty$, (but not for $p=\infty$, which is a singular case of elliptic regularity theory), and $\Omega'$ depends only on $\Omega, p, n$ and $M$, but not on $\A$. Moreover, if the right hand side of \eqref{initial_bound_recall} is $\epsilon M$ for some $\epsilon >0$ sufficiently small, then the above results holds globally in the sense that $\Omega'=\Omega$.

The proof of Theorem \ref{Thm_existence} is based on the following iteration scheme. The initial iterate is $U_1 = I$, so $U_1 \in SO(r,s)$. Assuming then that $U_k \in W^{1,2p}(\Omega)$ is given, define the subsequent iterate $U_{k+1} \in W^{1,2p}(\Omega)$ as the solution of 
\beq
\begin{cases}
\Delta U_{k+1} = U_k \delta \A - (U^T_k \eta)^{-1} \langle dU^T_k ; \eta dU_k \rangle  \label{iteration_RT1} \cr
U_{k+1} = I \hspace{1cm} \text{on} \ \ \ \partial\Omega.
\end{cases} 
\eeq
The iterates $U_k$ will in general not lie in the group $SO(r,s)$, but once we show that their limit $U$ is a solution of the reduced RT-equation \eqref{RT_2}, one can prove that $U$ lies in $SO(r,s)$ by the spectral argument in Sections \ref{Sec_SO(r,s)} and \ref{Sec_proof_exist_SO(r,s)}. 

Clearly the low regularities addressed here require interpreting the equations in a weak sense. 
The RT-equations admit the following weak formulation, based on the adjoints associated with the exterior derivative $d$ and its co-derivative $\delta$,   
\begin{align}
\Delta U[\phi] =  - \langle U\A ,d\phi \rangle_{L^2}     - \big\langle \langle dU;\A\rangle + (U^T\eta)^{-1} \langle dU^T ;\eta dU \rangle , \phi \big\rangle_{L^2},  \label{weak_RT2_preview}  
\end{align}
in terms of the weak Laplacian $\Delta U[\phi] \equiv - \langle dU,d\phi \rangle_{L^2}$,  matrix valued $0$-forms $\phi \in W^{1,(2p)^*}_0(\Omega)$ for testing, and the $L^2$ inner product $\langle\cdot , \cdot \rangle_{L^2}$ defined in \eqref{inner-product_L2}, and we wrote the first term in \eqref{RT_2} as $U\delta \A =  \delta (U\A) - \langle dU ; \A\rangle$ to arrange for total derivatives. Note that the choice of a weak formulation must preserve the overall logic of the RT-equations, and this is accomplished in Section \ref{Sec_weak} for the weak formulation \eqref{weak_RT2_preview}. In particular, we show that the constructions in Sections \ref{Sec_RT-eqn_optimal} and \ref{Sec_SO(r,s)}, which are based on point-wise substitution of the reduced RT-equation \eqref{RT_2}, carry over to the weak formulation \eqref{weak_RT2_preview} for which the substitution involves linear functionals.\footnote{Keep in mind that not all weak formulations are equivalent, and other formulations might not preserve this logic. This is well-known in the theory of shock waves, e.g. different weak formulations of Burger's equation lead to different weak solutions \cite{Smoller}.}  Note finally that the above iteration scheme \eqref{iteration_RT1} interpreted via the weak RT-equation \eqref{weak_RT2_preview} reduces the existence problem to the one of the classical weak Poisson equations on each component, by restriction to single matrix component of the test functions $\phi$, which gives an orthogonal decomposition of the space of test functions with respect to the Hilbert-Schmidt inner product (employed in our weak formulation).

\subsection{The $\epsilon$-rescaled equations} \label{Sec_proof_exist_scaling}

To handle the non-linearity of the reduced RT-equation \eqref{RT_2} and prove convergence of our iteration scheme, we now introduce a small parameter $\epsilon >0$ into the reduced RT-equations \eqref{RT_2}. In a first step we derive these rescaled RT-equations from a smallness assumption on the connection, and then explain how to arrange for this smallness assumption by localizing. To start, assume the $L^p$ bound \eqref{initial_bound_exist} on the connection is small in the sense that  
\beq \label{initial_bound_smallness}
\|\A\|_{L^{2p}(\Omega)} \leq \epsilon M
\eeq 
for some $\epsilon >0$ sufficiently small. We now define $\A^*$ by
\beq \label{A_star_smallness}
\A \equiv \epsilon \A^*  ,
\eeq
which satisfies the original $L^p$ bound \eqref{initial_bound_exist}, 
\beq \label{initial_bound_A-star}
\|\A^*\|_{L^{2p}(\Omega)} \leq M .
\eeq
We next define a matrix valued function $v$ by
\beq \label{rescaling_U}
U = I + \epsilon \: v,
\eeq
where $I$ denotes the identity on $\R^N$. Then, substituting \eqref{A_star_smallness} and \eqref{rescaling_U} into the reduced RT-equations \eqref{RT_2} yields
\beq \label{rescaling_RT_2}
\Delta v =  U \delta\A^*   -  \epsilon \: (U^T \eta )^{-1}  \langle d v^T ; \eta d v\rangle ,
\eeq
which is the sought-after rescaled RT-equation on the new unknown $v$. From \eqref{initial_bound_A-star}, we prove in Section \ref{Sec_proof_exist_iteration} the following existence result for \eqref{rescaling_RT_2}.

\begin{Prop}  \label{Prop_existence_RTscaled}
Let $n/2 < p < \infty$. Assume $\A^*$ satisfies the bound   
\beq \label{scaled_ext_Prop_eqn1}
\|\A^*\|_{L^{2p}(\Omega)} \leq M.
\eeq 
Then there exist a constant $\bar{\epsilon} >0$, depending only on $p,\Omega$ and $M$, such that for any $0 < \epsilon < \bar{\epsilon}$ there exists a solution $v \in W^{1,2p}(\Omega)$ of the rescaled RT-equations \eqref{rescaling_RT_2} with boundary data $v=0$ on $\partial \Omega$, and such that $U \equiv I + \epsilon v$ is pointwise an invertible matrix.   
\end{Prop}

Proposition \ref{Prop_existence_RTscaled} is an abstract existence theorem which applies to \eqref{rescaling_RT_2} whenever the scaled connection $\A^*$ in \eqref{rescaling_RT_2} satisfies \eqref{scaled_ext_Prop_eqn1}. In case of small connect in the sense of \eqref{initial_bound_smallness}, Proposition yields existence of solutions to the reduced RT-equations \eqref{RT_2} on the full domain $\Omega$.   We now explain how to derive the rescaled RT-equations \eqref{rescaling_RT_2} together with the bound \eqref{scaled_ext_Prop_eqn1} for the case of general ``large'' connections by localization. This requires two different cases, depending on whether our initial $L^p$ bound on $\A$ \eqref{initial_bound_recall} holds for $p=\infty$ or $p< \infty$. For this, assume without loss of generality that $\Omega = B_1(0) \equiv \{x\in \R^n\: | \: x_1^2+...+ x_n^2 < 1 \}$, the ball of radius $1$ (with respect to the Euclidean coordinate norm in $\R^n$) centered at $x(q)=0$.  \\

We begin with the case $p=\infty$. To start, restrict $\A$ to the ball of radius $\epsilon$, $B_{\epsilon}(0)$ for some $0<\epsilon <1$. We now introduce a change to coordinates 
\beq \label{rescaling}
x \to \tilde{x}\equiv x/\epsilon.
\eeq 
Note that coordinate transformations in the base space do not affect $U$ and $\eta$ on the fibres, only their coordinate derivatives, while connections transform as $1$-forms with matrix valued coefficients (acting on the fibre) otherwise unchanged, i.e., the coefficients of  $\A^{\tilde{x}} \equiv \A_\nu^{\tilde{x}} d\tilde{x}^\nu$ and $\A \equiv \A_i^x dx^i$ transform as $\A_\nu^{\tilde{x}} = \frac{\partial x^i}{\partial \tilde{x}^\nu} \A_i^x$. This gives
\beq \label{rescaling_A-star}
\A_i^{\tilde{x}} (\tilde{x}) =  \epsilon \A_i^x(\tilde{x}) \equiv \epsilon \A^*_i(\tilde{x}), 
\eeq
which is what we take for the definition of $\A^*$ in place of \eqref{A_star_smallness}. We now introduce the RT-equations in $\tilde{x}$-coordinates instead of $x$-coordinates and substitute \eqref{rescaling_A-star} and \eqref{rescaling_U}, from which the rescaled RT-equations \eqref{rescaling_RT_2} follow. The bound \eqref{Prop_existence_RTscaled} follows directly from the incoming bound \eqref{initial_bound_recall}, since the $L^\infty$-norm is preserved under coordinate transformations, which implies
\beq \label{scaling_L-infty}
\|\A^*(\tilde{x})\|_{L^{\infty}(\Omega)} = \|\A^x(x)\|_{L^{\infty}(B_\epsilon(0))}  \leq  \|\A^x(x)\|_{L^{\infty}(\Omega)} \leq M.
\eeq
Proposition \ref{Prop_existence_RTscaled} then yields existence of a solution $U \in W^{1,p}(\Omega)$ to the RT-equations in $\tilde{x}$-coordinates on $\Omega = B_1(0)$ for any $p<\infty$, and, by scaling invariance of the RT-equations \eqref{RT_2} under \eqref{rescaling}, a solution to the RT-equations in $x$-coordinates,  (c.f. Section \ref{Sec_Proof_Thm-ext}). From this we obtain (in Section \ref{Sec_proof_opt}) optimal connection regularity in $x$-coordinates on $\Omega' = B_\epsilon(0)$.   \\

When $p<\infty$, we use a cut off to derive the rescaled equations \eqref{rescaling_RT_2} in a different way and to arrange for the bound \eqref{scaled_ext_Prop_eqn1}. That is, we replace $\A$ in the RT-eqautions in the original $x$-coordinates by its cut off $\cutoff \A$, where $\cutoff$ is the characteristic function which vanishes outside $B_\delta(0)$ and is equal to $1$ inside $B_\delta(0)$, and we choose $\delta>0$ to be the largest $\delta>0$ such that $B_\delta(0) \subset \Omega$ and 
\beq \label{scaled_ext_Prop_eqn2}
\| \cutoff \A\|_{L^{2p}(\Omega)} \leq  \frac{\bar\epsilon}{2} \; M .
\eeq 
The resulting function $\delta=\delta(\bar\epsilon)$ exists and is strictly positive for $\bar\epsilon>0$, because
\beq \label{scaled_ext_Prop_eqn3}
\|\cutoff \A\|_{L^{2p}(\Omega)} = \|\A\|_{L^{2p}(B_{\delta}(0))}  \longrightarrow 0, \hspace{1cm} \text{as} \ \ \ \ \delta\to 0,
\eeq
by definition of the $L^p$-norm based on the Lebesgue integral. Moreover, $\delta(\bar\epsilon)$ depends only on the $L^{2p}$ norm of $\A$ near the coordinate center $x(q)=0$ as well as $M$, when taking the $\bar\epsilon$ from Proposition \ref{Prop_existence_RTscaled}, (where $\bar\epsilon$ is determined by $M$ alone).  Now, replacing  $\A$ by its cut off $\cutoff \A$ in the reduced RT-equations \eqref{RT_2}, in the original $x$-coordinates on $\Omega$, and defining $\A^*$ by
\beq \label{A-star_Lp-case}
\cutoff \A \equiv \epsilon \A^* ,
\eeq
in place of \eqref{A_star_smallness}, we obtain again the rescaled RT-equations \eqref{rescaling_RT_2} addressed by Proposition \ref{Prop_existence_RTscaled}. Finally, introducing the lower bound $\frac12 \bar{\epsilon} <\epsilon < \bar{\epsilon}$, it follows from \eqref{A-star_Lp-case} and \eqref{scaled_ext_Prop_eqn3} that 
\beq \label{scaled_ext_Prop_eqn4}
\|\A^*\|_{L^{2p}(\Omega)} = \frac{1}{\epsilon} \|\cutoff \A \|_{L^{2p}(\Omega)}  \leq \frac{\bar\epsilon}{2\epsilon} \; M   \leq M.
\eeq
This is the sought after bound \eqref{scaled_ext_Prop_eqn1}, and Proposition \ref{Prop_existence_RTscaled} yields existence of a solution $U$ of the reduced RT-equation for $\cutoff \A$ in place of $\A$ in the original $x$-coordinates on $\Omega$. The true geometric solution $U$ of the reduced RT-equation \eqref{RT_2} for $\A$ is then obtained by restriction to $B_\delta(0)$, which is a neighborhood only depending on $M$ and $\A$ near the coordinate center. From this we prove optimal connection regularity on $B_\delta(0)$ in Section \ref{Sec_proof_opt}.\footnote{Note that we work with a cut off function $\cutoff$ over the fixed domain $\Omega$, rather than restricting to the ball radius $\delta$ in our analysis, since the constant of Morrey's inequality depends on the domain.} 

Note that, for $p<\infty$, one can also establish the bound \eqref{scaled_ext_Prop_eqn1} on $\A^*$ defined in \eqref{rescaling} with $\cutoff \A$ in place of $\A$, by choosing $\delta$ such that $\|\cutoff \A\|_{L^{2p}(\Omega)} \leq \Big( \frac{\bar\epsilon}{2}\Big)^{\frac{n}{2p}} \; M$ in place of the bound \eqref{scaled_ext_Prop_eqn2} to compensate for the change from $x$ to $\tilde{x}$-coordinates in the $L^p$ norm of the connection $1$-form $\A$, i.e.,
\beq 
\|\A^{\tilde{x}} \|_{L^{2p}(B_1(0))}  = \epsilon^{1-\frac{n}{2p}} \|\A^x \|_{L^{2p}(B_{\epsilon}(0))} .
\eeq

Let us remark, that the above method of choosing $\delta$  to control the incoming bound \eqref{initial_bound_recall} works for our existence theory and the optimal regularity result in Theorem \ref{Thm_opt}. However, for Uhlenbeck compactness this method of choosing $\delta$ no longer works, because the cut off $\delta$ depends on the connection $\A$ near the origin, and thus can in general not be chosen uniformly over a sequence of connections. This is the reason why in Theorem \ref{Thm_compactness} we assume the uniform $L^\infty$ bound $\|\A_i\|_{L^\infty(\Omega)} \leq M$ on the sequence $\A_i$, which allows us to solve the RT-equations for each connection $\A_i$ in the same $\tilde{x}$-coordinates on the fixed domain $\Omega = B_1(0)$. This suffices to establish Uhlenbeck compactness in Section \ref{Sec_proof_Uhl}, since, under transform from $\tilde{x}$- to $x$-coordinates, the uniform $L^p$ bound on $d\A_i$ is preserved for $p \geq n/2$, that is,
\beq
\| d_{\tilde{x}}\A_i^{\tilde{x}} \|_{L^p(B_1(0))}   
= \epsilon^{2-\frac{n}{p}} \| d_{x}\A_i^x \|_{L^p(B_{\epsilon}(0))} 
\leq \bar\epsilon^{\; 2-\frac{n}{p}} M.
\eeq

\subsection{Existence of solutions and proof of Proposition \ref{Prop_existence_RTscaled}} \label{Sec_proof_exist_iteration}

For the proof of Proposition \ref{Prop_existence_RTscaled} it does not matter whether the rescaled connection $\A^*$ is defined in $\tilde{x}$-coordinates by \eqref{rescaling} or in $x$-coordinates by \eqref{A_star_smallness} or \eqref{A-star_Lp-case}, since the coordinate systems is kept fixed (denoted here by $x$) over the fixed domain $\Omega$ (which we may take to be $B_1(0)$). So assume, as in Proposition \ref{Prop_existence_RTscaled}, that $\A^*$ satisfies the bound \eqref{scaled_ext_Prop_eqn1} for some $M>0$ on the fixed neighborhood $\Omega$ in fixed $x$-coordinates. 

We now introduce an iteration scheme to prove the existence of a solution $v(x)$ of the $\epsilon$-rescaled reduced RT-equations \eqref{rescaling_RT_2} on $\Omega$, for $\epsilon>0$ sufficiently small. Start with the iterate $v_1 = 0$. Then $U_1 \equiv I + \epsilon v_1=I$ is the identity, which lies in $SO(r,s)$.  Assume now that $v_k \in W^{1,2p}(\Omega)$ is a given $(N\times N)$-matrix valued function with $v_k =0$ on $\partial \Omega$, $p>n/2$, so $U_k \equiv I + \epsilon v_k$ lies in $SO(r,s)$ on $\partial\Omega$, but not  necessarily in $SO(r,s)$ everywhere in $\Omega$. Define the next iterate $U_{k+1} = I + \epsilon v_{k+1} \in W^{1,2p}(\Omega)$ by solving  
\beq \label{it_rescal_RT_1}
\Delta v_{k+1} =  U_{k} \delta\A^*   -  \epsilon \: (U_k^T \eta )^{-1}  \langle d v_k^T ; \eta d v_k\rangle, 
\eeq
with Dirichlet boundary data $v_{k+1} =0$ on $\partial \Omega$. Throughout this section we interpret \eqref{it_rescal_RT_1} implicitly in the weak sense \eqref{weak_RT2_preview}, to account for $\Delta v_{k+1}$ and $\delta \A^*$ being only defined in a distributional sense. As shown in Lemma \ref{Lemma_ell_estimates} below, this defines a sequence of iterates $(v_k)_{k\in \mathbb{N}}$ in $W^{1,2p}(\Omega)$, lying in $SO(r,s)$ on $\partial\Omega$, but not necessarily in $SO(r,s)$ everywhere in $\Omega$.  

We now derive estimates in terms of our incoming curvature bound \eqref{initial_bound_recall}, for the iterates $v_k$ and differences of iterates $\overline{v_{k+1}} \equiv v_{k+1} - v_k$. We begin by clarifying the existence of the inverse of $U_k = I + \epsilon v_k$ when $\epsilon>0$ sufficiently small.

\begin{Lemma}  \label{Lemma_inv}
Let $p>n/2$ and assume $0< \epsilon < \epsilon_k$ for
\beq \label{def_epsilon_k}
\epsilon_k \equiv \frac{1}{2C_0 \|v_k\|_{1,2p}},
\eeq 
where $C_0>0$ is the constant from Morrey's inequality \eqref{Morrey_textbook}. Then the iterate $U_k \equiv I + \epsilon v_k$ is invertible with inverse $U_k^{-1} = I - \epsilon u_k$, and there exists a constant $C_{-1}>0$ depending only on $p,n,\Omega$ and $C_0$ such that                    
\beq  \label{estimate_inverse}
\| u_k\|_{W^{1,2p}}  \leq  C_{-1} \|v_k\|_{W^{1,2p}};
\eeq
and for $0<\epsilon< \min(\epsilon_k,\epsilon_{k-1})$, the difference $\overline{u_{k}} \equiv u_{k} - u_{k-1}$ satisfies the estimate  
\beq \label{estimate_inv_diff}
\| \overline{u_k}\|_{W^{1,2p}}  \leq  C_{-1} \|\overline{v_k}\|_{W^{1,2p}}.
\eeq
\end{Lemma}

\Proof
This is proven in \cite{ReintjesTemple_Uhl1}, Lemmas 6.1 and 6.3.
\QED

We now prove existence of iterates $v_{k+1}$ in $W^{1,2p}$ for $\epsilon>0$ sufficiently small.

\begin{Lemma}   \label{Lemma_ell_estimates}
Let $v_k \in W^{1,2p}_0(\Omega)$, $p>n/2$, and assume $0< \epsilon < \epsilon_k$ with $\epsilon_k$ given in \eqref{def_epsilon_k}. Then there exists a solution $v_{k+1} \in W^{1,2p}(\Omega)$ of \eqref{it_rescal_RT_1} with boundary data $v_{k+1}=0$ on $\partial\Omega$, and $v_{k+1}$ satisfies the elliptic estimate
\beq \label{estimate_iterates}
\|v_{k+1}\|_{W^{1,2p}} \ \leq \ C_1 \Big(  \| \A^* \|_{L^{2p}} + \epsilon \|\A^*\|_{L^{2p}} \|v_k\|_{1,2p} + \epsilon \big( 1 + \epsilon \|v_k\|_{1,2p} \big) \|v_k\|^2_{1,2p} \Big) ,
\eeq
where $C_1>0$ is some constant depending only on $p,n,\Omega$, $M$, $C_0$ and $C_{-1}$.
\end{Lemma}

\Proof
The existence of a solution $v_{k+1}$ to \eqref{it_rescal_RT_1} together with estimate \eqref{estimate_iterates} for any $k \in \mathbb{N}$ follows from Theorem \ref{Thm_Poisson}, (a standard result from elliptic PDE theory included with references in  the appendix), applied component-wise in the sense of restricting to single matrix component of the test functions $\phi$, (giving an orthogonal decomposition of the space of test functions with respect to the Hilbert-Schmidt inner product). To apply Theorem \ref{Thm_Poisson}, we need to show that the right hand side of \eqref{it_rescal_RT_1} is in $W^{-1,2p}(\Omega)$, the space of functionals over $W^{1,(2p)^*}_0(\Omega)$, (the conjugate exponent $(2p)^*$ satisfying $\frac{1}{2p} + \frac{1}{(2p)^*} =1$).

We begin by writing the first term in \eqref{it_rescal_RT_1}, $U_k \delta \A^*$, as a functional in $W^{-1,2p}(\Omega)$. For this we use Leibniz rule to write 
\beq
U_k \delta \A^* = \delta(U_k \A^*) - \langle dU_k ; \A^* \rangle, 
\eeq
and since $\langle dU_k ; \A^* \rangle \in L^p(\Omega) \subset W^{-1,2p}(\Omega)$ by Sobolev embedding for $p\geq n/2$, (as proven in Lemma 8.6 in \cite{ReintjesTemple_Uhl1} for $p>n/2$ and $p=n/2$), we find that the functional     
\beq \label{estimate_iterate_def0}
U_k \delta \A^*[\phi] = - \big\langle U_k \A^*, d\phi \big\rangle_{L^2} - \big\langle \langle dU_k ; \A^* \rangle, \phi \big\rangle_{L^2}
\eeq
is finite for any matrix valued $0$-form $\phi \in W^{1,(2p)^*}_0(\Omega)$. Integration by parts \eqref{partial_integration_matrix} then shows that \eqref{estimate_iterate_def0} indeed defines $U_k \delta \A^*$ as a functional in $W^{-1,2p}(\Omega)$. Now, to bound the $W^{-1,2p}$-norm of $U_k \delta \A^*$, we use H\"older's inequality twice to estimate
\beq  \label{estimate_iterate_bound1}
\big| \big\langle \langle dU_k ; \A^* \rangle, \phi \big\rangle_{L^2}  \big| 
\leq   \big\| \langle dU_k ; \A^* \rangle \big\|_{L^p} \|\phi\|_{L^{p^*}}
\leq C \|dU_k\|_{L^{2p}} \|\A^*\|_{L^{2p}} \|\phi\|_{W^{1,(2p)^*}},
\eeq
where the last inequality follows by Sobolev embedding, (with a generic constant $C>0$),  as in Lemma 8.6 in \cite{ReintjesTemple_Uhl1}. Using first H\"older's inequality and then Morrey's inequality, we bound the first terms in \eqref{estimate_iterate_def0} as
\beq \label{estimate_iterate_bound2}
\big| \big\langle U_k \A^*, d\phi \big\rangle_{L^2} \big| 
\leq  \|U_k\A^*\|_{L^{2p}}  \|d\phi\|_{L^{(2p)^*}}
\leq  C \|U_k\|_{W^{1,2p}} \|\A^*\|_{L^{2p}} \|\phi\|_{W^{1,(2p)^*}}.
\eeq
Combining \eqref{estimate_iterate_bound1} with \eqref{estimate_iterate_bound2}, substituting $U_k = I + \epsilon v_k$, and using that the operator norm is taken over test functions $\phi$ with $\|\phi\|_{W^{1,(2p)^*}}=1$, we conclude that the first source term in \eqref{it_rescal_RT_1} is bounded by
\beq \label{estimate_iterate_eqn1}
\| U_k \delta \A^*\|_{-1,2p} \leq   C \big(\| \A^* \|_{L^{2p}} + \epsilon\:  \|\A^*\|_{L^{2p}} \|v_k\|_{1,2p}  \big). 
\eeq

To handle the non-linear term in \eqref{it_rescal_RT_1}, which is one derivative more regular than $\delta\A^*$, we use again the Sobolev embedding $L^{p}(\Omega) \subset W^{-1,2p}(\Omega)$, (Lemma 7.6 in \cite{ReintjesTemple_Uhl1}), from which we derive the estimate\footnote{To be more precise, the embedding $L^{p}(\Omega) \subset W^{-1,2p}(\Omega)$ follows from applying Sobolev embedding to show that $\phi \in W^{1,(2p)^*} \subset L^{p^*}$ for $p\geq n/2$. Thus the dual paring of  $(U_k^T\eta)^{-1} \langle dv_k^T ;\eta dv_k \rangle$ with $\phi \in W^{1,(2p)^*}$ is finite, which gives the desired embedding. (Note, Sobolev embedding states here that $\phi \in L^q$ for $\frac1{q} = \frac1{(2p)^*} - \frac1{n}$, c.f. \eqref{Sobolev}. Now the required embedding $L^q \subset L^{p^*}$ holds if and only if $q \geq p^*$, which in turn is equivalent to the condition $p\geq n/2$.)} 
\begin{eqnarray} \label{estimate_iterate_eqn2}
\| (U_k^T\eta)^{-1} \langle dv_k^T ;\eta dv_k \rangle \|_{-1,2p} 
& \leq & C \|(U_k^T\eta)^{-1} \langle dv_k^T ;\eta dv_k \rangle \|_{L^p} \cr
& \leq & C \|(U_k^T\eta)^{-1}\|_{1,2p} \; \| dv_k \|_{L^{2p}}^2
\end{eqnarray}
in terms of a generic constant $C>0$, where the first inequality follows from Sobolev embedding, and the second inequality follows from Morrey's inequality \eqref{Morrey_textbook}, (used to bound the $L^\infty$ norm of $(U_k^T\eta)^{-1}$), in combination with the H\"older inequality, (used to bound the matrix valued inner product on $dv_k$). Using finally the bound \eqref{estimate_inverse} of Lemma \ref{Lemma_inv} on $U_k^{-1} = I - \epsilon u_k$, we obtain
\begin{eqnarray} \label{estimate_iterate_eqn3}
\| (U_k^T\eta)^{-1} \langle dv_k^T ;\eta dv_k \rangle \|_{-1,2p} 
& \leq & C \big( 1 + \epsilon \|v_k\|_{1,2p} \big) \|v_k\|^2_{1,2p} .
\end{eqnarray}
In summary, the bounds \eqref{estimate_iterate_eqn1} and \eqref{estimate_iterate_eqn3} show that the right hand side of \eqref{weak_RT_2} defines a functional over $W^{1,(2p)^*}$.

Theorem \ref{Thm_Poisson} now implies the existence of $v_{k+1} \in W^{1,2p}(\Omega)$ which solves \eqref{it_rescal_RT_1} with Dirichlet data $v_{k+1} =0$ on $\partial\Omega$. The standard elliptic estimate 
\eqref{Poisson-3} of Theorem \ref{Thm_Poisson} implies that the $W^{1,2p}$ norm of $v_{k+1}$ is bounded by the $W^{-1,2p}$ norm of the right hand side of \eqref{it_rescal_RT_1}, which in turn is bounded by \eqref{estimate_iterate_eqn1} and $\epsilon$ times \eqref{estimate_iterate_eqn3}. This proves the sought after estimate \eqref{estimate_iterates} with a generic constant $C_1>0$ depending only on $p,n,\Omega$ and $M$. 
\QED

We now introduce an induction assumption which is maintained by our iteration scheme and which provides uniform bounds. 

\begin{Lemma} \label{Lemma_ind_ass}
Let $k\in \mathbb{N}$ and assume 
\beq  \label{ind_ass}
0< \epsilon < \min\Big\{ \tfrac{1}{4M C_1 C_0} ,\tfrac{1}{4 C_1^2 M}, \tfrac{1}{4C_1 M ( 1  +  C_1  )}  \Big\}.
\eeq 
Then, if  
\beq  \label{ind_ass_uniform_bound_0}
\|v_k\|_{1,2p} \leq 2C_1 \|\A^*\|_{L^{2p}},
\eeq 
the subsequent iterate satisfies 
\beq \label{ind_ass_uniform_bound}
\|v_{k+1}\|_{1,2p} \leq 2 C_1 \|\A^*\|_{L^{2p}}.       
\eeq               
\end{Lemma}

\Proof
The induction assumption \eqref{ind_ass_uniform_bound_0} together with \eqref{ind_ass} imply
$$
0\; < \; \epsilon \; < \;  \frac{1}{2 C_0} \frac{1}{2 C_1 M}  
\; \leq \;   \frac{1}{2 C_0} \frac{1}{2 C_1 \|\A^*\|_{L^{2p}}}   
\; \leq \;  \frac{1}{2C_0} \frac{1}{\|v_k\|_{1,2p}} 
\; = \; \epsilon_k.
$$  
Thus the assumption $0<\epsilon<\epsilon_k$ of Lemma \ref{Lemma_ell_estimates} holds, by which the elliptic estimate \eqref{estimate_iterates} applies, to give
\begin{align}  \nonumber
\|v_{k+1}\|_{1,2p} 
& \leq  C_1 \Big(  \| \A^* \|_{L^{2p}} + \epsilon \|\A^*\|_{L^{2p}} \|v_k\|_{1,2p} + \epsilon \big( 1 + \epsilon \|v_k\|_{1,2p} \big) \|v_k\|^2_{1,2p} \Big)  \cr
& \leq  C_1 \|\A^*\|_{L^{2p}} +  \epsilon 2C_1^2 \| \A^* \|_{L^{2p}}^2 \big( 1  + 2 C_1 + \epsilon 4 C_1^2 \| \A^* \|_{L^{2p}} \big) .
\end{align}
We now use $\epsilon < \frac{1}{4 C_1^2 M}$ and $\epsilon < \frac{1}{4C_1(1+C_1)M}$, by assumption \eqref{ind_ass}, in combination with $\|\A^*\|_{L^{2p}} \leq M$, to bound $\epsilon\: 4 C_1^2 \| \A^* \|_{L^{2p}} \leq 1$ and $\epsilon\: 4C_1( 1  + C_1 ) \|\A^*\|_{L^{2p}} \leq 1$. From this and the above estimate, we obtain that
\begin{align} \nonumber
\|v_{k+1}\|_{1,2p} 
\ \leq\  C_1 \|\A^*\|_{L^{2p}} \Big(  1 +  \epsilon\: 4C_1( 1  + C_1 ) \|\A^*\|_{L^{2p}}  \Big) 
\ \leq\  2C_1 \|\A^*\|_{L^{2p}},
\end{align}
which is the sought-after bound \eqref{ind_ass_uniform_bound} for maintaining the induction assumption.
\QED

The uniform bound \eqref{ind_ass_uniform_bound} provided by Lemma \ref{Lemma_ind_ass} already implies weak $W^{1,2p}$ convergence of the iterates, and hence strong $L^{2p}$ convergence, but this is not sufficient to prove that the limit function solves the rescaled RT-equation \eqref{rescaling_RT_2}, because of the non-linear product of derivatives in \eqref{rescaling_RT_2}. For this, we now establish strong convergences in the $W^{1,2p}$ norm by deriving estimates on differences of iterates which are of order $\epsilon$. A geometric series argument then implies that the sequence of iterates converges strongly in $W^{1,2p}(\Omega)$. 

\begin{Lemma} \label{Lemma_est_differences}
For $0<\epsilon< \min(\epsilon_k,\epsilon_{k-1})$, there exists a constant $C_2>0$ depending only on $p,n,\Omega, C_0, C_{-1}$, such that the difference of iterates $\overline{v_{k+1}} \equiv v_{k+1} - v_k$ satisfies the estimate
\beq \label{estimate_differences}
\| \overline{v_{k+1}} \|_{W^{1,2p}} \ \leq \ \epsilon\: C_2 \Big( \|\A^*\|_{L^{2p}} + C(k) \Big)  \|\overline{v_k}\|_{W^{1,2p}},
\eeq
where
\beq \label{estimate_differences_Ck}
C(k) \equiv \big( 1+\epsilon \|v_{k-1}\|_{1,2p} \big) \big( \|v_k\|_{1,2p} + \|v_{k-1}\|_{1,2p} \big) + \epsilon \|v_k\|_{1,2p}.
\eeq
\end{Lemma}

\Proof
From \eqref{it_rescal_RT_1}, we find that 
\beq
\Delta \overline{v_{k+1}} = \epsilon \: \overline{v_k} \cdot \delta\A^*   -  \epsilon \:\eta^{-1} N_k
\eeq
where
\begin{eqnarray}
N_k &\equiv & (U_k^T )^{-1}  \langle d v_k^T ; \eta d v_k\rangle - (U_{k-1}^T )^{-1}  \langle d v_{k-1}^T ; \eta d v_{k-1}\rangle  \cr
&=& \big(1-\epsilon (u_{k-1})^T \big) \Big( \langle d \overline{v_k}^T ; \eta d v_k\rangle -  \langle d v_{k-1}^T ; \eta d \overline{v_k}\rangle \Big)  - \epsilon  \overline{u_k}^T \langle d v_k^T ; \eta d v_k\rangle , \nonumber
\end{eqnarray}
which is linear in differences of iterates and their inverses. The elliptic estimate \eqref{Poisson-3} gives us
\beq 
\|\overline{v_{k+1}} \|_{W^{1,2p}}     \leq     \epsilon\: C \big( \|\A^*\|_{L^{2p}}\: \|\overline{v_k}\|_{1,2p} +  \|N_k\|_{L^p} \big),
\eeq
where we used that $\|N_k\|_{-1,2p} \leq C \|N_k\|_{L^p}$  for $p> n/2$ by Sobolev embedding \cite[Lemma 8.6]{ReintjesTemple_Uhl1}, and where $C>0$ is some constant depending only on $p,n,\Omega$. Using now the Morrey inequality \eqref{Morrey_textbook} as well as estimates \eqref{estimate_inverse} and \eqref{estimate_inv_diff} on inverses, we bound $\|N_k\|_{L^p}$ as
\beq
\|N_k\|_{L^p} \leq C \Big( \big( 1+\epsilon \|v_{k-1}\|_{1,2p} \big) \big( \|v_k\|_{1,2p} + \|v_{k-1}\|_{1,2p} \big) + \epsilon \|v_k\|_{1,2p} \Big) \|\overline{v_k}\|_{1,2p},
\eeq
where $C>0$ is some constant depending only on $p,n,\Omega$. This establishes \eqref{estimate_differences} and proves the lemma.
\QED

We now prove convergence of the iteration scheme, assuming without loss of generality that $C_0, C_1, C_2 , M \geq 1$.

\begin{Prop} \label{Prop_iteration_conv}
There exist an $\bar{\epsilon} >0$, such that for any $0 < \epsilon < \bar{\epsilon}$ the iterates $v_k$ converge to some $v$ strongly in $W^{1,2p}(\Omega)$. Moreover, $v$ solves the rescaled RT-equations \eqref{rescaling_RT_2} with boundary data $v=0$ on $\partial \Omega$ and $U \equiv I + \epsilon v$ is invertible with $U^{-1} = I - \epsilon u$ for some $u \in W^{1,2p}(\Omega)$, such that
\beq \label{estimate_v_final}
\| v\|_{W^{1,2p}(\Omega)} + \| u\|_{W^{1,2p}(\Omega)} \leq  2 C_1 (1+C_{-1}) \; \|\A^*\|_{L^{2p}(\Omega)} .
\eeq                      
\end{Prop}

\Proof
Choose $\epsilon >0$ small enough to meet the $\epsilon$-bound \eqref{ind_ass}. Since $v_1=0$ meets the induction assumption \eqref{ind_ass_uniform_bound_0}, Lemma \ref{Lemma_ind_ass} implies for all $k \in \mathbb{N}$ that 
\beq \label{convergence_eqn1}
\| v_k\|_{1,2p} \leq 2C_1 \|\A^*\|_{L^{2p}} .
\eeq
Combining \eqref{convergence_eqn1} with the difference estimate \eqref{estimate_differences} of Lemma \ref{Lemma_est_differences} then gives
\begin{align} \label{convergence_eqn2}
\| \overline{v_{k+1}} \|_{W^{1,2p}}   &\leq  \ \epsilon\: C_2 \Big( \|\A^*\|_{L^{2p}} + C(k) \Big)  \|\overline{v_k}\|_{W^{1,2p}} \cr
 &\leq \ \epsilon\: C_2 \tfrac{1 + 2C_1 M + 12 C_1^2 M}{2C_1 M}  \|\A^*\|_{L^{2p}} \ \|\overline{v_k}\|_{W^{1,2p}}, 
\end{align}
since substitution of \eqref{convergence_eqn1} into \eqref{estimate_differences_Ck} gives
\begin{align}
C(k)  & \equiv \big( 1+\epsilon \|v_{k-1}\|_{1,2p} \big) \big( \|v_k\|_{1,2p} + \|v_{k-1}\|_{1,2p} \big) + \epsilon \|v_k\|_{1,2p}  \cr
& \leq \big( 1+\epsilon 2C_1 \|\A^*\|_{L^{2p}} \big) 4 C_1 \|\A^*\|_{L^{2p}} + \epsilon 2C_1 \|\A^*\|_{L^{2p}}   \cr
& \leq  4 C_1 \|\A^*\|_{L^{2p}}  +  \epsilon 8 C_1^2 \|\A^*\|_{L^{2p}}^2   + \epsilon 2C_1 \|\A^*\|_{L^{2p}}   \cr
& \leq (6 + \tfrac{1}{2 C_1^2 M}) C_1 \;\|\A^*\|_{L^{2p}},
\end{align}
where for the last inequality we use \eqref{ind_ass} to bound $\epsilon 8 C_1^2 \|\A^*\|_{L^{2p}}^2 \leq 2 \|\A^*\|_{L^{2p}}$ and $2 \epsilon < \frac{1}{2 C_1^2 M}$.
Restricting $\epsilon$ further to 
\beq \label{final_epsilon_bound}
0<\epsilon < \bar\epsilon \equiv \min\Big\{\tfrac{2C_1}{ C_2 \big( 1 + 2C_1 M + 12 C_1^2 M \big) }, \tfrac{1}{4M C_1 C_0} ,\tfrac{1}{4 C_1^2 M}, \tfrac{1}{4C_1 M ( 1  +  C_1  )}  \Big\}  ,
\eeq
estimate \eqref{convergence_eqn2} implies strong convergence in $W^{1,2p}(\Omega)$ by a standard geometric series argument for proving that the sequence of iterates is Cauchy, c.f. \cite[Prop. 9.11]{ReintjesTemple_Uhl1}. The limit $v$ therefore also solves \eqref{rescaling_RT_2} weakly with boundary data $v=0$ on $\partial\Omega$, because the convergence is strong enough to pass limits through the non-linear term of the right hand side of equations \eqref{it_rescal_RT_1}. By strong $W^{1,2p}$-convergence of $v_k$ to $v$, estimate \eqref{convergence_eqn1} on $v_k$ implies estimate \eqref{estimate_v_final} on $v$. By Lemma \ref{Lemma_inv} and in light of \eqref{final_epsilon_bound}, $U \equiv I + \epsilon v$ is invertible with $U^{-1} = I - \epsilon u$ for some $u \in W^{1,2p}(\Omega)$ with $W^{1,2p}$-norm bounded in terms of $\| v\|_{W^{1,2p}(\Omega)}$ from which we conclude \eqref{estimate_v_final}. This completes the proof of Proposition \ref{Prop_iteration_conv}.
\QED

\noindent Taken together, this completes the proof of Proposition \ref{Prop_existence_RTscaled}.

\subsection{Solutions generated by our iteration scheme always lie in $SO(r,s)$}  \label{Sec_proof_exist_SO(r,s)}  

We now prove that the solution $U= I + \epsilon v$ constructed in Proposition \ref{Prop_existence_RTscaled} through our iteration scheme does indeed lie pointwise in $SO(r,s)$. For this observe first that, for $\epsilon>0$ sufficiently small, $U = I + \epsilon v$ lies close to the identity and hence ${\rm det}(U)$ is close to $1$. Thus, since $U^T \eta U = \eta$ always implies ${\rm det}(U)=\pm 1$, proving $U^T \eta U = \eta$ directly implies $U \in SO(r,s)$. Defining $w \equiv U^T \eta U - \eta$, the sought after condition $U^T \eta U = \eta$ is then equivalent to $w=0$. In Section \ref{Sec_SO(r,s)} we showed that any solution $U$ of the RT-equation \eqref{RT_2} satisfies \eqref{SO(n)_eqn2}, that is,
\beq \label{SO_eqn0}
\Delta w = \delta \A^T \cdot w + w \cdot  \delta \A  ,
\eeq
with Dirichlet boundary data
\beq \label{SO_eqn0_data}
w=0 \ \ \ \ \ \text{on}  \ \ \ \ \partial\Omega,
\eeq 
(see Lemma \ref{Lemma_weak_w} below for a proof at low regularity). Clearly $w=0$ is a solution of \eqref{SO_eqn0}. However, in light of the Fredholm alternative, solutions to \eqref{SO_eqn0} might not be unique, and our iteration scheme above provides no information as to whether $w_k \equiv U_k^T \eta U_k - \eta$ converges to zero. To prove that $w=0$ does indeed hold, we need to incorporate the spectral argument outlined in Section \ref{Sec_SO(r,s)} into the framework of the rescaled RT-equations \eqref{rescaling_RT_2}.

To implement this spectral argument, write the rescaled RT-equation \eqref{rescaling_RT_2} as 
\beq \label{SO_eqn1}
\Delta v = (I + \epsilon v) \delta \A^* - \epsilon H(v),
\eeq
where for this argument
\beq
 H(v) \equiv (U^T \eta )^{-1}  \langle d v^T ; \eta d v\rangle.
\eeq
For $\lambda \in (0,1]$, consider the following modification 
\beq \label{SO_eqn3_RT}
\Delta v^\lambda = (I + \epsilon v^\lambda ) \lambda \delta \A^* - \epsilon H(v^\lambda),
\eeq 
based on replacing $\A^*$ by $\lambda \A^*$. Proposition \ref{Prop_existence_RTscaled} applies for each $\lambda \in (0,1]$ and yields the existence of a solution $U^{\lambda} = I + \epsilon v^{\lambda} \in W^{1,2p}(\Omega)$ of \eqref{SO_eqn3_RT}, because $\lambda \in (0,1]$ does not increase the initial bound \eqref{scaled_ext_Prop_eqn1} on $\A^*$.  Clearly $w^\lambda \equiv (U^\lambda)^T \eta U^\lambda - \eta$ solves \eqref{SO_eqn0} with $\lambda \A$ in place of $\A$, and following the argument between \eqref{SO(n)_eqn4} and \eqref{SO(n)_eqn5} in Section \ref{Sec_SO(r,s)}, we write this equation as 
\beq  \label{SO_eqn4}
\frac{1}{\lambda} w^\lambda = K (w^\lambda)
\eeq
for $K(w) \equiv \Delta^{-1} M(w)$ for $M(w) \equiv  \delta\A^T \cdot w + w \cdot \delta\A$.  It is clear that 
\beq \label{K}
K \equiv \Delta^{-1} M :  W^{1,2p}(\Omega) \longrightarrow W^{1,2p}(\Omega)
\eeq 
is a compact linear operator when $\A$ is smooth, and in Lemma \ref{Lemma_compact_op} below we prove that this also holds true when $\A \in L^{2p}(\Omega)$. The compactness of $K$ implies that the spectrum of $K$ is at most countable, (c.f. Theorem 4.25 in \cite{Rudin}).  We thus conclude that $w^\lambda=0$ is the unique solution of \eqref{SO_eqn4} for {\it almost every} $\lambda \in (0,1]$, i.e., for every $\lambda$ with $1/\lambda$ in the complement of the spectrum of $K$.  

We now prove continuity of the $w^\lambda$ with respect to $\lambda$, for solutions $w^\lambda$ generated by our iteration scheme. Continuity then implies $w^\lambda=0$ for all $\lambda \in (0,1]$, and in particular at $\lambda=1$, and this gives $w=0$. The continuity of $w^\lambda$ with respect to $\lambda$ is a consequence of the following lemma, which asserts that solution $U^\lambda$ generated by our iteration scheme are continuous in $\lambda$ with respect to $W^{1,2p}$, (in both $x$ and $\tilde{x}$-coordinates).

\begin{Lemma} \label{Lemma_est_diff_2}
Let $U^{\lambda} = I + \epsilon v^{\lambda}$ and $U^{\lambda'} = I + \epsilon v^{\lambda'}$ be solutions of the rescaled RT-equation \eqref{SO_eqn3_RT} generated by our iteration scheme. Then there exists $\epsilon_0 >0$ such that for any $0< \epsilon< \epsilon _0$ we have the estimate
\beq \label{estimate_diff_2}
\|U^{\lambda} - U^{\lambda'} \|_{W^{1,2p}} \leq C_3 |\lambda - \lambda'|,
\eeq
for some constant $C_3>0$ depending only on $p,n,\Omega$, $C_0, C_1$ and $C_2$.
\end{Lemma}

\Proof
By Lemma \ref{Lemma_ind_ass} there exists a $\epsilon_0' >0$ such that for all $0<\epsilon < \epsilon_0'$ any solution $v^\lambda$ generated by the iteration scheme with $\lambda \in (0,1]$ satisfies the uniform bound    
\beq \label{SO_eqn5}
\| v^\lambda \|_{1,2p} \leq 2 C_1 M.
\eeq
This provides a uniform bound on $v^\lambda$ for $\lambda \in (0,1]$. 

Next, we establish the Cauchy property \eqref{estimate_diff_2} in the spirit of the proof of Lemma \ref{Lemma_est_differences}. Since $v^\lambda$ and $v^{\lambda'}$ both satisfy the RT-equation \eqref{SO_eqn3_RT}, we have
\beq \label{SO_eqn6}
\Delta (v^{\lambda} - v^{\lambda'})    
=  (\lambda - \lambda')  \delta \A^*  + \epsilon  (\lambda v^\lambda - \lambda' v^{\lambda'}) \delta \A^*  - \epsilon \big( H(v^\lambda)  - H(v^{\lambda'}) \big).
\eeq
Elliptic estimate \eqref{Poisson-3} together with source estimates (as in Section \ref{Sec_proof_exist_iteration}) applied to \eqref{SO_eqn6} then imply
\begin{eqnarray}  \label{SO_eqn7}
\|v^{\lambda} - v^{\lambda'} \|_{1,2p}   
& \leq &    |\lambda - \lambda'| \: \|\delta \A^*\|_{-1,2p}  + \epsilon\:  \| \lambda v^\lambda - \lambda' v^{\lambda'} \|_{1,2p} \: \|\delta \A^*\|_{-1,2p}  \cr
& &  + \: \epsilon\:  \| H(v^\lambda) - H(v^{\lambda'})\|_{-1,2p}.
\end{eqnarray}
As before $\|\delta \A^*\|_{-1,2p}  \leq \|\A^*\|_{L^{2p}} \leq M$ by definition of the operator norm, while   
\beq  \label{SO_eqn7a}
\| \lambda v^\lambda - \lambda' v^{\lambda'} \|_{1,2p}  \leq |\lambda - \lambda'| \: \|v^\lambda\|_{1,2p} + |\lambda'| \: \| v^\lambda - v^{\lambda'} \|_{1,2p} .
\eeq
Moreover, by the corresponding estimate on differences of iterates in the proof of Lemma \ref{Lemma_est_differences}, we can bound the non-linear term by
\beq  \label{SO_eqn7b}
\| H(v^\lambda) - H(v^{\lambda'})\|_{-1,2p} \leq C_2 \Big( \|\A^*\|_{L^{2p}} + C(\lambda,\lambda') \Big)  \| v^\lambda - v^{\lambda'}\|_{1,2p}  ,
\eeq
where
\beq
C(\lambda,\lambda') \equiv \big( 1+\epsilon \|v^{\lambda'}\|_{1,2p} \big) \big( \|v^\lambda\|_{1,2p} + \|v^{\lambda'}\|_{1,2p} \big) + \epsilon \: \|v^\lambda\|_{1,2p}.
\eeq
Using now \eqref{SO_eqn5}, in the form $\| v^\lambda \|_{1,2p} \leq 2 C_1 M$ and $\| v^{\lambda'} \|_{1,2p} \leq 2 C_1 M$, to further bound the right hand sides of \eqref{SO_eqn7a} and \eqref{SO_eqn7b}, we find that \eqref{SO_eqn7} gives the bound
\beq  \label{SO_eqn8}
\|v^{\lambda} - v^{\lambda'} \|_{1,2p}   \ \leq \   
(M + \epsilon \: 2 C_1 M) \: |\lambda - \lambda'|  + \epsilon \: C_2' \| v^\lambda - v^{\lambda'}\|_{1,2p},
\eeq
for some constant $C_2' >0$ depending only on $p,n,\Omega$, $M$, $C_0$, $C_1$ and $C_2$. Choosing now $0< \epsilon < \epsilon_0 \equiv \min(\epsilon_0', \frac{1}{2C_2'})$, subtraction of the last term in \eqref{SO_eqn8} from both sides of the equation yields
\beq  \label{SO_eqn9}
\frac{1}{2} \|v^{\lambda} - v^{\lambda'} \|_{1,2p}   \ \leq \   
(M + \epsilon 2 C_1 M) \: |\lambda - \lambda'|  ,
\eeq
which implies the desired estimate \eqref{estimate_diff_2}.
\QED

The continuity of $U^\lambda$ and $w^\lambda$ with respect to $\lambda$ in the $W^{1,2p}$-norm asserted by Lemma \ref{Lemma_est_diff_2}, implies that $w^\lambda=0$ for all $\lambda$, including the original $\lambda=1$. This proves that $U$ lies in $SO(r,s)$ point-wise everywhere in $\Omega$. This completes the proof of Theorem \ref{Thm_existence}, subject to establishing the compactness of $K$ in \eqref{K} when $\A \in L^{2p}(\Omega)$ in Lemma \ref{Lemma_compact_op} of Section \ref{Sec_weak}. The complete proof of the different cases of Theorem \ref{Thm_existence} is given in Section \ref{Sec_Proof_Thm-ext} below.

\subsection{Proof of Theorem \ref{Thm_existence}} \label{Sec_Proof_Thm-ext}

We now complete the proof of Theorem \ref{Thm_existence}, addressing each case separately: \\ 

\noindent {\bf Global (i):} Let  $p> n/2$ and assume $\A \in L^{2p}(\Omega)$. Theorem \ref{Thm_existence} then asserts the existence of some $\tilde{\epsilon}>0$, depending only on $\Omega$ and $p$, such that, if 
\beq \label{Global_i_proof_1}
\|\A\|_{L^{2p}(\Omega)} \leq \tilde{\epsilon} ,
\eeq
then there exists a solution $U \in W^{1,2p}(\Omega)$ of the reduced RT-equations \eqref{RT_2} which is pointwise in $SO(r,s)$, and satisfies estimate \eqref{uniform_bound_U_global}. 

To begin the proof, note that \eqref{Global_i_proof_1} is equivalent to the smallness assumption \eqref{initial_bound_smallness} for $M=1$. Thus, it suffices to set $\tilde{\epsilon} \equiv \bar\epsilon/2$ in terms of the $\bar\epsilon>0$ established by Proposition \ref{Prop_existence_RTscaled}, which depends only on $\Omega, p$ and $M$. Namely, it now follows that $\A^* = \frac{1}{\epsilon} \A$ meets the bound \eqref{scaled_ext_Prop_eqn1} for $M=1$, for any $\epsilon \in (\bar\epsilon/2,\bar\epsilon)$. Thus, Proposition \ref{Prop_existence_RTscaled} applies and yields the existence of a solution $v \in W^{1,2p}(\Omega)$ to the rescaled RT-equations \eqref{rescaling_RT_2} for all $\epsilon \in (0,\bar\epsilon)$, where $\bar\epsilon$ depends only on $\Omega, p$ and $M$. Defining now $U = I + \epsilon v$ for some $\epsilon \in (0,\bar\epsilon)$, it follows that $U$ solves the reduced RT-equations \eqref{RT_2}, since the right hand side of
\begin{eqnarray}  \nonumber
\Delta U -  U \delta\A   +  (U^T \eta )^{-1}  \langle dU^T ; \eta dU\rangle 
=  \epsilon \Big( \Delta v -  U \delta\A^*   + \epsilon (U^T \eta )^{-1}  \langle dv^T ; \eta dv\rangle \Big) 
\end{eqnarray}  
vanishes by the rescaled RT-equations \eqref{rescaling_RT_2}. Proposition \ref{rescaling_U} yields further that $U$ is invertible with inverse given by $U^{-1} = I - \epsilon u$. Moreover, by the results in Section \ref{Sec_proof_exist_SO(r,s)} it follows that $U = I + \epsilon v$ lies pointwise in $SO(r,s)$, (assuming for the moment the validity of Lemma \ref{Lemma_compact_op}, which is proven in Section \ref{Sec_weak} below). To establish estimate \eqref{uniform_bound_U_global} on $U = I + \epsilon v$ and $U^{-1} = I - \epsilon u$, recall that, by Proposition \ref{Prop_iteration_conv}, $v$ and $u$ satisfy
\beq \label{exist-proof_estimate_v-u}
\| v\|_{W^{1,2p}(\Omega)} + \| u\|_{W^{1,2p}(\Omega)} \leq  2 C_1 (1+C_{-1}) \; \|\A^*\|_{L^{2p}(\Omega)}.
\eeq 
From this, keeping in mind that $\A = \epsilon \A^*$, we obtain directly the bound
\begin{eqnarray} \label{exist-proof_estimate_U}
\|U-I\|_{W^{1,2p}(\Omega)}   + \|U^{-1}-I\|_{W^{1,2p}(\Omega)}   
&=  & \ \epsilon \Big(  \| v\|_{W^{1,2p}(\Omega)} + \| u\|_{W^{1,2p}(\Omega)}  \Big) \cr
&\leq &\  C(M) \; \|\A\|_{L^{2p}(\Omega)},
\end{eqnarray}
where $C(M) \equiv 2 C_1 (1+C_{-1}) >0$. For $M=1$, this is the sought-after estimate \eqref{uniform_bound_U_global} of Theorem \ref{Thm_existence}. This completes the proof of Theorem \ref{Thm_existence} - Global (i). \\

\noindent {\bf Global (ii):} Assume $\A$ meets the smallness condition \eqref{Global_i_proof_1} for $p=\infty$ and $\tilde{\epsilon} \equiv \bar\epsilon/2$. Then Theorem \ref{Thm_existence} asserts that, for any $\bar{p} <\infty$, there exists a solution $U \in W^{1,\bar{p}}(\Omega)$ of \eqref{RT_2} in $SO(r,s)$ which satisfies estimate \eqref{uniform_bound_U_global} for $\bar{p}$ in place of $2p$. 

To prove this statement, note that \eqref{Global_i_proof_1} for $p=\infty$ implies that the $L^p$ norm of $\A$ satisfies \eqref{initial_bound_smallness} for $M=|\Omega|$, (the volume of $\Omega$), since
\beq \label{Global_ii_proof}
\|\A\|_{L^{2p}(\Omega)} \leq \|\A\|_{L^{\infty}(\Omega)} |\Omega| \leq \epsilon |\Omega|. 
\eeq
The proof of Theorem \ref{Thm_existence} - Global (i) now applies for any $\bar{p}\equiv 2p <\infty$, and yields a solution $U \in W^{1,\bar{p}}(\Omega)$ of \eqref{RT_2} in $SO(r,s)$ which satisfies \eqref{exist-proof_estimate_U} for $\bar{p}$ in place of $2p$ over $\Omega$. This completes the proof of Theorem \ref{Thm_existence} - Global (ii).\footnote{The case $\bar{p}=\infty$ is not established, since it is a singular case of elliptic regularity theory. Note that our existence theory in Section \ref{Sec_proof_exist_iteration} might yield different solutions $U$ for different $\bar{p}<\infty$.} \\

\noindent {\bf Local (i):} Let $p\in (n/2,\infty)$ and assume $\A \in L^{2p}(\Omega)$ satisfies
\beq \label{Local_i_proof_1}
\|\A\|_{L^{2p}(\Omega)} \leq M
\eeq
for some constant $M>0$. Then, for any point $q \in \Omega$, Theorem \ref{Thm_existence} asserts the existence of a neighborhood $\Omega' \subset \Omega$ of $q$, (depending on $\A$), and a solution $U \in W^{1,2p}(\Omega')$ of \eqref{RT_2} in $SO(r,s)$ which satisfies estimate \eqref{uniform_bound_U_local} on $\Omega'$, for some constant $C(M) > 0$ depending on $\Omega', p$ and $M$.

To begin the proof, we take for $\Omega'$ the ball $B_\delta(q)$ on which $\cutoff \A = \A$, (assuming for ease $q=0$), and we define $\A^*$ by \eqref{A-star_Lp-case} as $\cutoff \A \equiv \epsilon \A^*$. For $\A^*$ to meet the bound \eqref{scaled_ext_Prop_eqn1} required in Proposition \ref{Prop_existence_RTscaled} for existence, we choose the cut off $\delta>0$ so small such that \eqref{scaled_ext_Prop_eqn2}  holds, that is, small enough such that $\|\cutoff \A\|_{L^{2p}(\Omega)} \leq  \frac{\bar\epsilon}{2} \; M$.  Since $\bar\epsilon$ in Proposition \ref{Prop_existence_RTscaled} depends only on $\Omega, p$ and $M$, it follows that the cut off $\delta>0$ can be chosen in the first step in terms of $\Omega, p$ and $M$ alone. So, for $\frac12 \bar{\epsilon} <\epsilon < \bar{\epsilon}$, \eqref{scaled_ext_Prop_eqn2} implies that $\A^*$ satisfies \eqref{scaled_ext_Prop_eqn1}, by which Proposition \ref{Prop_existence_RTscaled} yields the existence of a solution $v \in W^{1,2p}(\Omega)$ of the rescaled RT-equations. 

Now, the proof of Theorem \ref{Thm_existence} - Global (i) applies, and implies that $U=I + \epsilon v \in  W^{1,2p}(\Omega)$ solves the reduced RT-equations \eqref{RT_2} for $\cutoff \A$ in place of $\A$, is invertible, lies in $SO(r,s)$ and satisfies estimate \eqref{exist-proof_estimate_U} on $\Omega$ with $\cutoff \A$ in place of $\A$. We now obtain the true geometric solution $U$ of the reduced RT-equation \eqref{RT_2} by restricting $U$ to $\Omega' \equiv B_\delta(0)$, (where $\cutoff \A = \A$), and we obtain the sought-after estimate \eqref{initial_bound_exist} from \eqref{exist-proof_estimate_U}, (with $\A$ replaced by $\cutoff \A$), since
\small
\beq \nonumber
\|U-I \|_{W^{1,2p}(\Omega')} + \|U^{-1}-I\|_{W^{1,2p}(\Omega')}  
\leq C(M) \; \|\cutoff\A\|_{L^{2p}(\Omega)} 
 =  C(M) \; \|\A\|_{L^{2p}(\Omega')}.
\eeq
\normalsize
This completes the proof of Theorem \ref{Thm_existence} - Local (i).\\

\noindent {\bf Local (ii):} Assume $\A \in L^{\infty}(\Omega)$ satisfies \eqref{Local_i_proof_1} for $p=\infty$ in $x$-coordinates on $\Omega \equiv B_1(0)$. Theorem \ref{Thm_existence} then asserts that, for any $\bar{p} <\infty$, there exists a solution $U \in W^{1,\bar{p}}(\Omega')$ of \eqref{RT_2} in $SO(r,s)$, which satisfies \eqref{uniform_bound_U_local} for $\bar{p}$ in place of $2p$; and $\Omega'$ depends only on $\Omega, p$ and $M$, but is independent of $\A$. 

To prove this statement, we restrict $\A$ to $B_\epsilon(0) \subset \Omega$ and take $\Omega'\equiv B_\epsilon(0)$ as the neighborhood of $q$, for some $\epsilon \in (\bar\epsilon/2, \bar\epsilon)$, (again assuming $q=0$). Note that $\Omega'$ depends only on $\Omega, p$ and $M$, but is independent of $\A$, since $\epsilon \in (\bar\epsilon/2, \bar\epsilon)$ and $\bar\epsilon$ depends only on $\Omega, p$ and $M$.  We next change to coordinates $x \to \tilde{x}\equiv x/\epsilon$, and define $\A^*$ by \eqref{rescaling_A-star} as the original connection $\A^x$ in $x$-coordinates with components transformed as scalar functions to $\tilde{x}$-coordinates, $\A^*(\tilde{x}) \equiv \A^x(\tilde{x})$. Since the $L^\infty$-norm is preserved under coordinate transformations, as shown in \eqref{scaling_L-infty}, it follows from \eqref{Global_ii_proof} that $\A^*$ meets the bound \eqref{scaled_ext_Prop_eqn1} of Proposition \ref{Prop_existence_RTscaled} in $\tilde{x}$-coordinates, for any $\bar{p}\equiv 2p <\infty$, (with $|\Omega| M$ in place of $M$). Proposition \ref{Prop_existence_RTscaled} now yields the existence of a solution $v \in W^{1,\bar{p}}(\Omega)$ of the rescaled RT-equations in $\tilde{x}$-coordinates on $\Omega=B_1(0)$. Thus, as before, $U = I + \epsilon v$ is invertible, lies in $SO(r,s)$, solves the reduced RT-equations \eqref{RT_2} and satisfies \eqref{uniform_bound_U_local} for $\bar{p}$ in place of $2p$.  Under transformation from $\tilde{x}$ to $x = \epsilon \tilde{x}$, the reduced RT-equations transform by 
\begin{eqnarray} \label{RT_2_scalinginvariance}
&& \Delta_x U  -   U \delta_x\A^x   +  (U^T \eta )^{-1}  \langle d_xU^T ; \eta d_xU\rangle    \cr
&&= \frac{1}{\epsilon^2} \Big( \Delta_{\tilde{x}} U -   U \delta_{\tilde{x}}\A^{\tilde{x}}   +  (U^T \eta )^{-1}  \langle d_{\tilde{x}} U^T ; \eta d_{\tilde{x}} U\rangle  \Big),
\end{eqnarray}  
since $\A$ transforms as a $1$-form, $\A_i^x = \A_\nu^{\tilde{x}} \frac{\partial \tilde{x}^\nu}{\partial x^i} =\frac{1}{\epsilon} \A_i^{\tilde{x}}$. From the scale invariance \eqref{RT_2_scalinginvariance}, we conclude  that $U'(x) \equiv U(x/\epsilon)$ solves the RT-equations in $x$-coordinates on $\Omega' = B_\epsilon(0)$. We denote this solution $U'(x)$ again as $U(x)$. Finally, we find that $U$ satisfies the sought-after estimate \eqref{uniform_bound_U_local} for $\bar{p}=2p$, in $x$-coordinates, because transformation from $\tilde{x}$- to $x$-coordinates modifies the bound by factors of $\epsilon \in (\bar\epsilon/2, \bar\epsilon)$, which, since $\bar\epsilon$ depends only on $\Omega, p$ and $M$, can be absorbed into the constant $C(M)$. This completes the proof of  Theorem \ref{Thm_existence} - Local (ii). \\

Recall, to establish $U$ in $SO(r,s)$ by the results of Section \ref{Sec_proof_exist_SO(r,s)}, we assumed in the above proof that compactness of the operator $K$ defined in \eqref{K} extends to the low regularity $\A \in L^{2p}(\Omega)$. This requires a careful analysis of the weak equations, accomplished in Lemma \ref{Lemma_compact_op} in Section \ref{Sec_weak} below. Assuming here the validity of Lemma \ref{Lemma_compact_op}, the proof of Theorem \ref{Thm_existence}  is complete.  \hfill $\Box$

\section{Weak formalism}   \label{Sec_weak}

We finally address the weak formulation of the full RT-equations \eqref{RT_1} - \eqref{RT_2} required for the low regularity $\A \in L^{2p}(\Omega)$, $d\A \in L^{\bar{p}}(\Omega)$ and $U \in W^{1,2p}(\Omega)$.   The weak form of the reduced RT-equation \eqref{RT_2} was already given in \eqref{weak_RT2_preview}, and our existence theory in Proposition \ref{Prop_existence_RTscaled} establishes existence of weak solutions $U \in W^{1,2p}(\Omega)$ of \eqref{RT_2} when $\A \in L^{2p}(\Omega)$ for $p\in (n/2,\infty)$. It remains only to extend the derivation of equation \eqref{SO(n)_eqn2} for $w = U^T \eta U - \eta$ in Section \ref{Sec_SO(r,s)} and the proof of Lemma \ref{Lemma_Gammati'} to the weak setting, because both involve point-wise substitutions of the RT-equation \eqref{RT_2}, but at the level of distributions such substitutions require justification. In Lemma \ref{Lemma_compact_op}, we use the weak form of \eqref{SO(n)_eqn2} to show that $K$ in \eqref{K} is compact.  

Based on the integration by parts formula \eqref{partial_integration_matrix}, we define the weak Laplacian on $1$-forms as 
\beq
\Delta \Ati [\psi] \equiv - \langle d\Ati , d\psi \rangle_{L^2} - \langle \delta \Ati, \delta \psi \rangle_{L^2},
\eeq
acting on matrix valued $1$-forms $\psi \in W^{1,\hat{p}^*}_0(\Omega)$ for $\hat{p} = \min\{p,\bar{p}\}$, and we define the weak Laplacian on $0$-forms as 
\beq
\Delta U[\phi] \equiv - \langle dU,d\phi \rangle_{L^2},
\eeq
for matrix valued $0$-forms $\phi \in W^{1,(2p)^*}_0(\Omega)$, where $\frac1{\hat{p}} + \frac1{\hat{p}^*} =1$ and $\frac1{2p} + \frac1{(2p)^*} =1$, and where $\langle\cdot , \cdot \rangle_{L^2}$ is the $L^2$ inner product defined in \eqref{inner-product_L2}. (Note that $\delta \phi=0$ for all $0$-forms.) In terms of this, we define the weak form of the RT-equations \eqref{RT_1} - \eqref{RT_2} as
\begin{align}
\Delta \Ati [\psi] &= - \langle d\A , d\psi \rangle_{L^2} + \langle dU^{-1} \wedge dU , d\psi \rangle_{L^2} \label{weak_RT_1} \\
\Delta U[\phi] &=  - \langle U\A ,d\phi \rangle_{L^2}     - \big\langle \langle dU;\A\rangle + (U^T\eta)^{-1} \langle dU^T ;\eta dU \rangle , \phi \big\rangle_{L^2},  \label{weak_RT_2}  
\end{align}
where we wrote the first term in \eqref{RT_2} as $U\delta \A =  \delta (U\A) - \langle dU ; \A\rangle$ to arrange for total derivatives, c.f. \eqref{weak_RT2_preview}. From the H\"older inequality and Sobolev embedding it follows that the right hand side of \eqref{weak_RT_2} defines a functional in $W^{-1,2p}(\Omega)$, (c.f. the proof of Lemma \ref{Lemma_ell_estimates}), and that the right hand side of \eqref{weak_RT_1} defines a functional in $W^{-1,\hat{p}}(\Omega)$, which is consistent with $\Ati \in W^{1,\hat{p}}(\Omega)$ and $U \in W^{1,2p}(\Omega)$.  

In the next lemma we extend the derivation of equation \eqref{SO(n)_eqn2} in Section \ref{Sec_SO(r,s)} to the weak setting here. That is, for solutions $U \in W^{1,2p}(\Omega)$ of the reduced RT-equation \eqref{RT_2}, we show that $w \equiv U^T \eta U - \eta$ solves 
$$
\Delta w = \delta \A^T \cdot w + w \cdot  \delta \A 
$$ 
in the weak sense that, for any matrix valued $0$-form $\phi \in W^{1,(2p)^*}_0(\Omega)$,
\beq \label{weak_w_eqn0}
\Delta w [\phi] \equiv - \langle dw,d\phi \rangle_{L^2} = - \langle \A^T ,d(\phi w) \rangle_{L^2} + \langle \A ,d(w \phi) \rangle_{L^2} .
\eeq

\begin{Lemma} \label{Lemma_weak_w}
Assume $U \in W^{1,2p}(\Omega)$ is a weak solution of \eqref{weak_RT_2}, $p>n/2$. Then $w \equiv U^T \eta U - \eta$ is in $W^{1,2p}(\Omega)$ and satisfies \eqref{weak_w_eqn0} for any matrix valued $0$-form $\phi \in W^{1,(2p)^*}_0(\Omega)$. 
\end{Lemma}

\Proof
The $W^{1,2p}$ regularity of $w \equiv U^T \eta U - \eta$ follows from $U \in W^{1,2p}(\Omega)$, since $W^{1,2p}(\Omega)$ is closed under multiplication. To show that $w$ satisfy \eqref{weak_w_eqn0}, it suffices to consider test functions $\phi \in C_0^\infty(\Omega)$, since $C_0^\infty(\Omega)$ is contained in $W^{1,(2p)^*}(\Omega)$. Substitution of $dw = dU^T \cdot \eta U + U^T \eta \cdot dU$ into the weak Laplacian in \eqref{weak_w_eqn0} gives 
\begin{eqnarray} \label{weak_w_eqn1}
- \Delta w [\phi] \equiv  \langle d w , d\phi \rangle_{L^2} 
= \langle dU^T , d\phi \cdot U^T \eta \rangle_{L^2}  +  \langle dU , \eta U \cdot d\phi \rangle_{L^2},
\end{eqnarray}
where we used cyclic commutativity of the matrix trace in $\langle\cdot , \cdot \rangle_{L^2}$; (recall that $\langle A, B\rangle_{L^2} \equiv  \int_\Omega  {\rm tr} \langle A^T ; B \rangle dx$ on matrix valued $k$-forms $A$ and $B$, c.f. \eqref{inner-product_L2}). Writing $\eta U\cdot  d\phi  = d(\eta U \phi) - d(\eta U) \cdot \phi$, the second term in \eqref{weak_w_eqn1} is
\begin{eqnarray} \label{weak_w_eqn2}
\langle dU , \eta U \cdot d\phi \rangle_{L^2}
&=&  \langle dU , d(\eta U \phi) \rangle_{L^2} - \langle dU , d(\eta U) \cdot \phi \rangle_{L^2}. 
\end{eqnarray}
Since $\phi' \equiv \eta U \phi \in W^{1,(2p)^*}(\Omega)$  is a test function, we can substitute the reduced RT-equation \eqref{weak_RT_2} for the first term in \eqref{weak_w_eqn2}, which gives
\begin{align} \label{weak_w_eqn3}
\langle dU , \eta U \cdot d\phi \rangle_{L^2}
&=   \langle U\A ,d\phi' \rangle_{L^2}    + \big\langle \langle dU;\A\rangle , \phi' \big\rangle_{L^2}  \cr
& + \big\langle  (U^T\eta)^{-1} \langle dU^T ;\eta dU \rangle , \phi' \big\rangle_{L^2} -  \langle dU , d(\eta U) \cdot \phi \rangle_{L^2}.    
\end{align}
By cyclic commutativity of the trace, and since 
\beq \label{weak_techeqn}
\big\langle\langle A ; B \rangle , C \big\rangle_{L^2} = \big\langle B , A^T C \big\rangle_{L^2}
\eeq 
for matrix valued $1$-forms $A, B$ and matrix valued $0$-forms $C$, the terms in the last line of \eqref{weak_w_eqn3} cancel, while the remaining terms simplify to
\begin{eqnarray} \label{weak_w_eqn4}
\langle dU , \eta U \cdot d\phi \rangle_{L^2}
&=&   \langle \A ,d(U^T \eta U \phi) \rangle_{L^2}  .
\end{eqnarray}

To compute the first term on the right hand side of \eqref{weak_w_eqn1}, we use that $ {\rm tr}(A^T) =  {\rm tr}(A)$ for any matrix $A$, so $\langle dU^T , d\phi \cdot U^T \eta \rangle_{L^2}   =  \langle dU , \eta U d\phi^T  \rangle_{L^2}$. Thus, applying \eqref{weak_w_eqn4} with $\phi$ replaced by $\phi^T$, we obtain
\begin{eqnarray} \label{weak_w_eqn5}
\langle dU^T , d\phi \cdot U^T \eta \rangle_{L^2}   =   \langle \A ,d(U^T \eta U \phi^T) \rangle_{L^2}  =    \langle \A^T ,d( \phi U^T \eta U) \rangle_{L^2} ,
\end{eqnarray}
where the last equality follows by symmetry of $U^T\eta U$. Substituting \eqref{weak_w_eqn4} and \eqref{weak_w_eqn5} back into \eqref{weak_w_eqn1}, and writing $U^T \eta U = w +\eta$ gives
\begin{eqnarray} \label{weak_w_eqn6}
- \Delta w [\phi] 
&=&  \langle \A^T ,d(\phi (w +\eta)) \rangle_{L^2} + \langle \A ,d((w +\eta) \phi) \rangle_{L^2} .
\end{eqnarray}
Using finally the Lie Algebra property $\A^T \eta + \eta \A =0$ in the form 
\beq \nonumber
0 = \langle ( \A^T \eta + \eta \A^T ) , d\phi \rangle_{L^2} 
= \langle \A^T , d(\phi \eta) \rangle_{L^2}  + \langle \A , d(\eta \phi) \rangle_{L^2},
\eeq
the terms containing $\eta$ in \eqref{weak_w_eqn6} cancel. By this cancellation, \eqref{weak_w_eqn6} reduces to the sought after equation \eqref{weak_w_eqn0}, which completes the proof.
\QED

We now prove that $K \equiv \Delta^{-1} M$ is a compact operator when $\A \in L^{2p}(\Omega)$, as claimed in \eqref{K}. In this case $M(w) \equiv  \delta\A^T \cdot w + w \cdot \delta\A $ must now be interpreted in the weak sense as the bilinear form on the right hand side of \eqref{weak_w_eqn0},
\beq \label{M}
M(w)[\phi] \equiv - \langle \A^T ,d(\phi w) \rangle_{L^2} + \langle \A ,d(w \phi) \rangle_{L^2}.
\eeq

\begin{Lemma} \label{Lemma_compact_op}
For $\A \in L^{2p}(\Omega)$, the linear operator $K$ is compact as an operator from $W^{1,2p}(\Omega)$ to $W^{1,2p}(\Omega)$.
\end{Lemma}

\Proof
So let $M(w)[\phi] \equiv - \langle \A^T ,d(\phi w) \rangle_{L^2} + \langle \A ,d(w \phi) \rangle_{L^2}$ and let $\A \in L^{2p}(\Omega)$. We first show that $M(w) \in W^{-1,2p}(\Omega)$ for any $w \in W^{1,2p}(\Omega)$, and that $M: W^{1,2p}(\Omega) \to W^{-1,2p}(\Omega)$ is a bounded linear operator. For this use the Leibnitz rule to write
\beq \label{compact_op_eqn1}
|\langle \A ,d(w \phi) \rangle_{L^2}| \leq  |\langle \A ,dw \cdot  \phi \rangle_{L^2}| + |\langle \A ,w\cdot d\phi \rangle_{L^2}|.
\eeq
We then combine the Morrey and H\"older inequalities to estimate
\beq \label{compact_op_eqn2}
|\langle \A ,w d\phi \rangle_{L^2}|  
\leq C_0 \|w\|_{W^{1,2p}} |\langle \A , d\phi \rangle_{L^2}|
\leq C_0 \|w\|_{W^{1,2p}} \|\A\|_{L^{2p}} \|d\phi\|_{L^{(2p)^*}},
\eeq
and we use H\"older's inequality (twice) and Sobolev embedding to estimate
\beq \label{compact_op_eqn3}
 |\langle \A ,dw \cdot  \phi \rangle_{L^2}|  
 \leq  \| \langle A; dw\rangle \|_{L^p} \| \phi \|_{L^{p^*}}
 \leq C_S \| A\|_{L^{2p}} \|dw\|_{L^{2p}}  \| \phi \|_{W^{1,(2p)^*}},
\eeq
(see \cite[Lemma 7.6]{ReintjesTemple_Uhl1} for a proof of the inequality $\| \phi \|_{L^{p^*}} \leq C_S  \| \phi \|_{W^{1,(2p)^*}}$). From \eqref{compact_op_eqn1} - \eqref{compact_op_eqn3} it follows that 
\beq \label{compact_op_eqn4}
|M(w)[\phi]| \leq \max\{C_0,C_S\}   \| A\|_{L^{2p}} \|w\|_{W^{1,2p}} \| \phi \|_{W^{1,(2p)^*}},
\eeq
which implies that $M: W^{1,2p}(\Omega) \to W^{-1,2p}(\Omega)$ is a bounded linear operator. Now, to prove that $K \equiv \Delta^{-1} M: W^{1,2p}(\Omega) \to W^{1,2p}(\Omega)$ is a compact operator, observe that $\Delta^{-1}$ is a compact operator from $W^{-1,2p}(\Omega)$ to $W^{1,2p}(\Omega)$ by elliptic regularity theory, see for instance \eqref{Poisson-4}. Thus, since the composition of a compact operator with a bounded operator is again compact \cite{Rudin}, $K$ is indeed a compact operator. This completes the proof.
\QED

We now extend Lemma \ref{Lemma_Gammati'} to the weak setting, as required for the proofs of Theorems \ref{Thm_equiv} and \ref{Thm_opt} in Section \ref{Sec_proofs_rest} below. We here place $\A$ and $d\A$ in different $L^p$ spaces, a slight generalization convenient for the proof of Theorem \ref{Thm_opt}.

\begin{Lemma} \label{Lemma_weak_Ati}
Let $p\in (n/2,\infty)$ and $\bar{p} \in (1, \infty)$. Assume $\A \in L^{2p}(\Omega)$ and $d\A \in L^{\bar{p}}(\Omega)$, and assume $U \in W^{1,2p}(\Omega)$ in $SO(r,s)$ is a solution of the reduced RT-equations \eqref{weak_RT_2}. Then $\Ati' \equiv \A - U^{-1} dU$ solves \eqref{weak_RT_1} with $\Ati$ replaced by $\Ati'$, and $\Ati' \in W^{1,\hat{p}}(\Omega'')$ for any $\Omega''$ compactly contained in $\Omega$, and $\hat{p} = \min\{\bar{p},p\}$.
\end{Lemma}

\Proof
Without loss of generality, assume $1< \bar{p} \leq p$, so $\hat{p} = \bar{p}$. To begin observe $\Ati' \equiv \A - U^{-1} dU \in L^{2p}(\Omega)$, since $\A \in L^{2p}$ and since $U^{-1} \in W^{1,2p}(\Omega)$ is H\"older continuous by Morrey's inequality. By \eqref{Leibniz_rule_J-application}, we have $d\Ati' = d \A - dU^{-1} \wedge d U$, which implies by the H\"older inequality applied to the last term that $d\Ati' \in L^{\bar{p}}(\Omega)$, keeping in mind that $d\A \in L^{\bar{p}}(\Omega)$ and $\bar{p} \leq p$. From $d\A \in L^{\bar{p}}(\Omega)$, $d\Ati' \in L^{\bar{p}}(\Omega)$ and $dU^{-1} \wedge d U \in L^{\bar{p}}(\Omega)$, we now obtain that
\beq \label{weak_Ati_eqn1}
- \langle d\Ati' , d\psi \rangle_{L^2} =  - \langle d\A , d\psi \rangle_{L^2} + \langle dU^{-1} \wedge dU , d\psi \rangle_{L^2}
\eeq
is well-defined for any matrix valued $1$-form $\psi \in W^{1,\bar{p}^*}_0(\Omega)$. Thus, to show that $\Ati'$ solves the first RT-equation \eqref{weak_RT_1} in a weak sense, it remains to show that
\beq \label{weak_Ati_eqn2}
 \langle \delta \Ati', \delta \psi \rangle_{L^2} =0.
\eeq
In analogy to the proof of Lemma \ref{Lemma_Gammati'}, we accomplish this by showing that
\beq \label{weak_Ati_eqn3}
\delta \Ati' [\phi]  \equiv - \langle \Ati' , d \phi \rangle_{L^2} =0
\eeq
for any matrix valued $0$-form $\phi \in C^\infty_0(\Omega)$, which then yields \eqref{weak_Ati_eqn3} since $C^\infty_0(\Omega)$ is dense in $L^{\bar{p}^*}(\Omega)$. For this, using that $U^T \eta U = \eta$ implies $U^{-1} = \eta U^T \eta$, we compute that
\begin{eqnarray} \label{weak_Ati_eqn4}
- \delta \Ati' [\phi] 
&=& \langle \A , d\phi \rangle_{L^2} - \langle dU , \eta U \eta \cdot d\phi \rangle_{L^2}   \cr
&= & \langle \A , d\phi \rangle_{L^2} + \langle dU , d(\eta U \eta) \phi \rangle_{L^2} 
- \langle dU , d(\eta U \eta  \phi) \rangle_{L^2} .
\end{eqnarray}
We then substitute the weak RT-equation \eqref{weak_RT_2} for the last term in \eqref{weak_Ati_eqn4}, which is possible because $\phi' \equiv \eta U \eta  \phi$ lies in the space of test functions, $\phi' \in W^{1,(2\bar{p})^*}_0(\Omega)$, for $\phi \in C^\infty_0(\Omega)$ and $U \in W^{1,2p}(\Omega) \subset W^{1,(2\bar{p})^*}(\Omega)$, since $2p>n\geq 2$ and $(2\bar{p})^*<2$. 
This substitution then gives
\begin{eqnarray} \label{weak_Ati_eqn5}
- \delta \Ati' [\phi] 
&= & \langle \A , d\phi \rangle_{L^2}  - \langle U\A ,d\phi' \rangle_{L^2} - \big\langle \langle dU;\A\rangle  , \phi' \big\rangle_{L^2}  \cr 
&& + \langle dU , d(\eta U \eta) \phi \rangle_{L^2}   - \big\langle (U^T\eta)^{-1} \langle dU^T ;\eta dU \rangle , \phi' \big\rangle_{L^2} .
\end{eqnarray}
The terms in the first line mutually cancel, as well as those in the second line, which can be shown using $U^T \eta U = \eta$, cyclic commutativity of the trace in \eqref{def_inner-product}, and identity \eqref{weak_techeqn} above. This shows that \eqref{weak_Ati_eqn3} holds and proves that $\Ati'$ solves the sought after RT-equation \eqref{weak_RT_1} in the weak sense.

To prove the regularity gain from $\Ati' \in L^{2p}(\Omega)$ to $\Ati' \in W^{1,\bar{p}}(\Omega'')$, recall the weak form \eqref{weak_RT_1} of the first RT-equation,
\beq \label{weak_RT_1_proof}
\Delta \Ati [\psi] = - \langle d\A , d\psi \rangle_{L^2} + \langle dU^{-1} \wedge dU , d\psi \rangle_{L^2} .
\eeq
Observe that the right hand side in \eqref{weak_RT_1_proof} defines a functional over $W^{1,\bar{p}^*}(\Omega)$. Namely, since $d\A \in L^{\bar{p}}(\Omega)$ and $dU^{-1} \wedge dU \in L^{p}(\Omega) \subset L^{\bar{p}}(\Omega)$, (for $\bar{p} \leq p$), the $L^2$ inner products on the right hand side of \eqref{weak_RT_1} are both finite by H\"older's inequality and hence define functionals over $W^{1,\bar{p}^*}(\Omega)$.  Interior elliptic regularity theory, as stated in Theorem \ref{Thm_Poisson_interior} below, (keeping in mind that $\Ati'  \in L^{2p}(\Omega)$ by construction), then yields $\Ati' \in W^{1,\bar{p}}(\Omega'')$ any $\Omega''$ compactly contained in $\Omega$. This completes the proof.
\QED

\section{Proof of optimal regularity and Uhlenbeck compactness} \label{Sec_proofs_rest}

Theorem \ref{Thm_existence}, our basic existence result, was proven  in Section \ref{Sec_proof_existence}. We now complete the proofs of Theorems \ref{Thm_equiv}, \ref{Thm_opt} and \ref{Thm_compactness}, which are based on Theorem \ref{Thm_existence}, together with Lemmas \ref{Lemma_Gammati'} and \ref{Lemma_weak_Ati} which give the regularity boost for $\Ati'.$ For this, let $\A \equiv \A_\textbf{a}$ be the connection components in a gauge $\textbf{a}$ of a connection $\A_{\VM}$ on an $SO(r,s)$ vector bundle $\VM$ with base manifold $\M \equiv \Omega \subset \R^n$ open and bounded.

\subsection{Proof of Theorem \ref{Thm_equiv}}   \label{Sec_proof_equiv}

Assume throughout that $\A \in L^{2p}(\Omega)$ and $d\A \in L^{\bar{p}}(\Omega)$, for $p \in (n/2,\infty]$, $\bar{p} \in (1,\infty)$. 

Theorem \ref{Thm_equiv}, part (ii), states that if there exists a gauge transformation $U\in W^{1,2p}(\Omega)$ pointwise in  $SO(r,s)$, such that the gauge transformed connection $\Aop$ in \eqref{connection_transfo_VB} has optimal regularity $\Aop \in W^{1,\hat{p}}(\Omega)$ for $\hat{p} = \min\{p,\bar{p}\}$, then $\Ati \equiv U^{-1} \Aop U \in W^{1,\hat{p}}$ together with $U$ solve the RT-equations \eqref{RT_1} and \eqref{RT_2}, respectively. The proof of this statement follows directly from the derivation of the RT-equation in Section \ref{Sec_RT-eqn_derivation}, and it is straightforward to extend this derivation to the weak formalism of Section \ref{Sec_weak} and the regularities addressed here.

Theorem \ref{Thm_equiv}, part (i), states that if there exists a solution $U \in W^{1,2p}(\Omega)$ in $SO(r,s)$ of the reduced RT-equations \eqref{RT_2}, then the gauge transformed connection $\Aop$ in \eqref{connection_transfo_VB} has optimal regularity $\Aop \in W^{1,\hat{p}}(\Omega)$. The proof of this statement follows from Lemma \ref{Lemma_weak_Ati}. That is, Lemma \ref{Lemma_weak_Ati} asserts that $\Ati' \equiv \A - U^{-1} dU$ has regularity $W^{1,\hat{p}}(\Omega)$ as a result of $\Ati'$ being a solution of the first RT-equation \eqref{RT_1}. Thus, since the connection in the gauge $\textbf{b} = U\cdot \textbf{a}$ satisfies
\beq   \label{product_eqn_A}
\Aop = U \Ati' U^{-1},
\eeq
it follows that $\Aop \in W^{1,\hat{p}}(\Omega)$, which is optimal regularity. (That the products in \eqref{product_eqn_A} do indeed stay in $W^{1,\hat{p}}(\Omega)$ is proven in equations \eqref{optimal_reg_eqn_proof} - \eqref{optimal_estimate_proof_techeqn4} below). This completes the proof of Theorem \ref{Thm_equiv}.  \hfill $\Box$

\subsection{Optimal regularity - Proof of Theorem \ref{Thm_opt}}  \label{Sec_proof_opt}

Let $p \in (n/2,\infty)$, $\bar{p} \in (1,\infty)$ and $n\geq 2$. Assume that the connection $\A$ in $x$-coordinates satisfies the curvature type bound
\beq \label{bound_opt_proof} 
\|(\A,d\A)\|_{L^{2p,\bar{p}}(\Omega)}  \equiv \|\A \|_{L^{2p}(\Omega)} + \|d\A \|_{L^{\bar{p}}(\Omega)} \; \leq \; M,
\eeq  
for some constant $M>0$.  
Then for any point $q\in \Omega$, Theorem \ref{Thm_existence} yields the existence of a solution $U\in W^{1,2p}(\Omega')$ of the reduced RT-equations \eqref{RT_2} in $SO(r,s)$, on some neighborhood $\Omega' \subset \Omega$ of $q$, satisfying
\beq \label{uniform_bound_U_proof}
\|U-I\|_{W^{1,2p}(\Omega')} + \|U^{-1}-I\|_{W^{1,2p}(\Omega')}  \; \leq\; C(M)\; \|\A\|_{L^{2p}(\Omega')} ,
\eeq
for some constant $C(M) > 0$ depending only on $\Omega', p, n$ and $M$; in the case of connections which satisfy the smallness condition \eqref{initial_bound_small} we can take $\Omega' = \Omega$. We now prove that the $SO(r,s)$ gauge transformation $U$ establishes optimal regularity $W^{1,\hat{p}}$ for $\hat{p} = \min\{p,\bar{p}\}$. For ease of presentation we assume without loss of generality that $1< \bar{p} \leq p$, by which $\hat{p} = \bar{p}$, and we address the four cases of Theorem \ref{Thm_opt} in the end.

Under the above assumptions, we find that $\Ati' \equiv \A - U^{-1} dU \in L^{2p}$, and Lemma \ref{Lemma_weak_Ati} further yields that $\Ati'$ solves the first RT-equation \eqref{RT_1}. We now show that interior elliptic regularity theory (Theorem \ref{Thm_Poisson_interior}), applied to \eqref{RT_1}, yields on any neighborhood $\Omega''$ compactly contained in $\Omega'$ the estimate 
\beq \label{uniform_bound_A_proof2}
\|\Ati' \|_{W^{1,\bar{p}}(\Omega'')}   \; \leq\; C(M)\; \|(\A,d\A)\|_{L^{2p,\bar{p}}(\Omega')} ,
\eeq
for some constant $C(M) > 0$ depending only on $\Omega'', \Omega', \bar{p}, p, n$ and $M$. Throughout the proof we denote with $C(M)>0$ a universal constant. By Lemma \ref{Lemma_weak_Ati}, $\Ati'$ is in $ W^{1,\bar{p}}(\Omega'')$ and solves the RT-equation \eqref{RT_1}, which we state here as
\beq \label{RT_U-I}
\Delta \Ati' = \delta d \A -\delta\big( d(U^{-1}-I)\wedge d(U-I) \big),
\eeq 
by replacing $U$ and $U^{-1}$ by $U-I$ and $U^{-1}-I$. Applying the interior elliptic estimate  \eqref{Poisson-4} of Theorem \ref{Thm_Poisson_interior} to \eqref{RT_U-I}, we obtain 
\beq \label{uniform_bound_A_proof3}
\|\Ati'\|_{W^{1,\bar{p}}(\Omega'')} \leq C \Big( \|\Ati'\|_{L^{\bar{p}}(\Omega')} +   \big\|\delta d \A -\delta\big( d(U^{-1}-I)\wedge d(U-I) \big) \big\|_{W^{-1,\bar{p}}(\Omega')} \Big).
\eeq
We now bound the $L^{\bar{p}}$ norm of $\Ati' = \A - U^{-1} d(U-I) \in L^{2p}(\Omega)$, using  \eqref{uniform_bound_U_proof} and Morrey's inequality for $2p>n$, by
\begin{eqnarray}\label{uniform_bound_A_proof4}
\|\Ati'\|_{L^{\bar{p}}(\Omega')}  
&\leq & \|\Ati'\|_{L^{2p}(\Omega')}  \cr
&\leq & \|\A\|_{L^{2p}(\Omega')} + \|U^{-1}\|_{W^{1,2p}(\Omega')} \|U-I\|_{W^{1,2p}(\Omega')} \cr
&\leq & C(M) \|\A\|_{L^{2p}(\Omega')},
\end{eqnarray}
for some suitable constant $C(M)>0$ bounding  $\|U^{-1}\|_{W^{1,2p}(\Omega')}$. We next bound the $W^{-1,\bar{p}}$-norm on the right hand side of \eqref{RT_U-I}, using the H\"older inequality \eqref{Hoelder} and estimate \eqref{uniform_bound_U_proof},  by 
\begin{eqnarray}  \label{uniform_bound_A_proof5}
&&\| \delta d \A -\delta\big( d(U^{-1}-I)\wedge d(U-I) \big) \|_{W^{-1,\bar{p}}(\Omega')}   \cr
&&\leq \|d \A\|_{L^{\bar{p}}(\Omega')} +\|U^{-1}-I\|_{W^{1,2p}(\Omega')} \|U-I\|_{W^{1,2p}(\Omega')}  \cr
&&\leq   C(M) \|(\A,d\A)\|_{L^{2p,\bar{p}}(\Omega')},
\end{eqnarray}
keeping in mind that $1< \bar{p} \leq p$ by assumption.  Combining \eqref{uniform_bound_A_proof3}, \eqref{uniform_bound_A_proof4} and \eqref{uniform_bound_A_proof5}, we obtain the sought-after estimate \eqref{uniform_bound_A_proof2}.

Based on estimate \eqref{uniform_bound_A_proof2} for $\Ati'$, we now prove that the connection in gauge $\textbf{b} = U\cdot \textbf{a}$ has optimal regularity
\beq \label{optimal_reg_eqn_proof}
\Aop \equiv U \Ati' U^{-1}  \in W^{1,\hat{p}}(\Omega''),
\eeq 
by establishing the estimate 
\beq \label{optimal_estimate_proof}
\|\Aop \|_{W^{1,\hat{p}}(\Omega'')}   \leq  C(M) \|(\A,d\A)\|_{L^{2p,\bar{p}}(\Omega')},
\eeq
where $\hat{p} = \min\{p,\bar{p}\} = \bar{p}$ by our assumption here.  For this we estimate the product $U \Ati'$, using first the Leibnitz rule and then Morrey's inequality \eqref{Morrey_textbook}, by
\begin{eqnarray} \label{optimal_estimate_proof_techeqn1}
\|U\Ati'\|_{W^{1,\bar{p}}} 
&\leq & \| U\Ati'\|_{L^{\bar{p}}} + \| U \partial\Ati' \|_{L^{\bar{p}}}  + \| \partial U \Ati' \|_{L^{\bar{p}}} \cr
&\leq & C_0  \|U\|_{W^{1,\bar{p}}}  \Big( \|\Ati'\|_{L^{\bar{p}}}  +  \| \partial\Ati'\|_{L^{\bar{p}}}  \Big)  + \| \partial U \Ati' \|_{L^{\bar{p}}},
\end{eqnarray}
where $\partial$ denotes component-wise differentiation. We estimate the last term in \eqref{optimal_estimate_proof_techeqn1} using H\"older's inequality as
\beq  \label{optimal_estimate_proof_techeqn2}
\|\partial U \Ati' \|_{L^{\bar{p}}} \leq \| \partial U\|_{L^{2\bar{p}}} \|\Ati'\|_{L^{2\bar{p}}} 
\leq  \|U\|_{W^{1,2p}} \|\Ati'\|_{L^{2\bar{p}}}   ,
\eeq
keeping in mind that $\bar{p}\leq p$. Combining \eqref{optimal_estimate_proof_techeqn1} and \eqref{optimal_estimate_proof_techeqn2} gives
\begin{eqnarray} \label{optimal_estimate_proof_techeqn3}
\|U\Ati'\|_{W^{1,\bar{p}}(\Omega'')} 
\leq  (C_0 +1) \|U\|_{W^{1,\bar{p}}(\Omega'')} \Big( \|\Ati'\|_{W^{1,\bar{p}}(\Omega'')}  + \| \Ati' \|_{L^{2p}(\Omega'')} \Big),
\end{eqnarray}
and using \eqref{uniform_bound_A_proof4} to bound $\|\Ati'\|_{L^{2p}(\Omega'')} \leq  C(M) \|\A\|_{L^{2p}(\Omega')}$,  together with \eqref{uniform_bound_A_proof2} and \eqref{uniform_bound_U_proof} to bound $\|\Ati'\|_{W^{1,\bar{p}}(\Omega'')}$ and $\|U\|_{W^{1,\bar{p}}(\Omega'')} \leq \|U\|_{W^{1,p}(\Omega'')}$, we obtain
\beq \label{optimal_estimate_proof_techeqn4}
\|U\Ati'\|_{W^{1,\bar{p}}(\Omega'')} 
\leq  C(M)  \|(\A,d\A)\|_{L^{2p,\bar{p}}(\Omega')}.
\eeq
The sought-after estimate \eqref{optimal_estimate_proof} follows by a derivation analogous to that of \eqref{optimal_estimate_proof_techeqn4}. To complete the proof, we now address each case in Theorem \ref{Thm_opt} separately.\\

\noindent {\bf Global (i):} We assume \eqref{bound_opt_proof} for $p \in (n/2,\infty)$ and $\bar{p}\in (1,\infty)$, such that the $L^{2p}$ norm of $\A$ is small in the sense of \eqref{initial_bound_small}. By Theorem \ref{Thm_existence} the regularizing $SO(r,s)$ gauge transformation $U \in W^{1,2p}(\Omega)$ is then defined on all of $\Omega$ and satisfies \eqref{uniform_bound_U_proof} for $\Omega' =\Omega$. By the above proof, we obtain estimate \eqref{optimal_estimate_proof}  for $\Omega'=\Omega$ which directly implies $\Aop\in W^{1,\hat{p}}(\Omega'')$ together with the sought-after estimate \eqref{uniform_bound_A_global}. 

\vspace{.1cm} \noindent {\bf Global (ii):}  We assume \eqref{bound_opt_proof} for $p =\infty$ and $\bar{p}\in (1,\infty)$, and we assume $\A$ satisfies the smallness condition \eqref{initial_bound_small} for $p=\infty$. For any $\tilde{p} <\infty$, Theorem \ref{Thm_existence} yields an $SO(r,s)$ gauge transformation $U \in W^{1,2\tilde{p}}(\Omega)$ which solves the reduced RT-equations \eqref{RT_2} and satisfies \eqref{uniform_bound_U_proof} on $\Omega$. Choosing $\tilde{p} \geq \bar{p}$, the above proof yields $\Ati' \in W^{1,\hat{p}}(\Omega'')$ and \eqref{optimal_estimate_proof} for $p =\infty$, which is the sough-after estimate of Theorem \ref{Thm_opt} - Global (ii). 

\vspace{.1cm} \noindent {\bf Local (i):} We assume \eqref{bound_opt_proof} holds for $p \in (n/2,\infty)$  and $\bar{p} \in (1,\infty)$, but without the smallness condition \eqref{initial_bound_small}. The argument of the case Global (i) now applies on an $\A$-dependent neighborhood $\Omega' \subset \Omega$ of $q$ and yields $\Aop\in W^{1,\hat{p}}(\Omega'')$ and \eqref{optimal_estimate_proof} for any  $\Omega''$ compactly contained in $\Omega'$.

\vspace{.1cm} \noindent {\bf Local (ii):} We assume \eqref{bound_opt_proof} for $p =\infty$ and $\bar{p}\in (1,\infty)$, but without the smallness condition \eqref{initial_bound_small}. The argument of case Global (ii) now applies on a neighborhood $\Omega' \subset \Omega$ of $q$ and yields $\Aop\in W^{1,\hat{p}}(\Omega'')$ and \eqref{uniform_bound_A_global} for $p=\infty$, on any  $\Omega''$ compactly contained in $\Omega'$. The neighborhood $\Omega'$, determined in Theorem \ref{Thm_existence} - Local (ii), is now independent of $\A$. This completes the proof of Theorem \ref{Thm_opt}.  \hfill $\Box$

\subsection{Uhlenbeck compactness - Proof of Theorem \ref{Thm_compactness}}    \label{Sec_proof_Uhl}  

We now prove the three different cases of Theorem \ref{Thm_compactness}.

\vspace{.1cm} \noindent {\bf Global (i):} Let $p \in (n/2,\infty)$ and $\bar{p} \in (1,\infty)$, and assume a sequence of connections $(\A_i)_{i\in \mathbb{N}}$ on $\VM$ in fixed gauge $\textbf{a}$ satisfies the uniform bound 
\beq \label{initial_bound_2_proof}
\|(\A_i,d\A_i)\|_{L^{2p,\bar{p}}(\Omega)} \equiv \|\A_i \|_{L^{2p}(\Omega)} + \|d\A_i \|_{L^{\bar{p}}(\Omega)} \; \leq \; M
\eeq   
for some constant $M>0$. Assume further each $\A_i$ satisfies the smallness condition \eqref{initial_bound_small}.  

\vspace{.1cm} \noindent {\bf (a)} Applying our optimal regularity result in Theorem \ref{Thm_opt} - Global (i) to each connection $\A_i$, it follows that on every neighborhood $\Omega''$ compactly contained in $\Omega$ there exists an $SO(r,s)$ gauge transformation $U_i \in W^{1,2p}(\Omega'')$ such that the components $\A_{\textbf{b}_i}$ of $\A_i$ in gauge $\textbf{b}_i = U_i \cdot \textbf{a}$ have optimal regularity $\A_{\textbf{b}_i}\in W^{1,\hat{p}}(\Omega'')$ for $\hat{p}=\min\{p,\bar{p}\}$, with uniform bound 
\beq \label{uniform_bound_A_2_proof}
\|\A_{\textbf{b}_i} \|_{W^{1,\hat{p}}(\Omega'')}   \; \leq \;  C(M) M ,
\eeq
for a constant $C(M) > 0$ depending only on $\Omega'', \Omega, p, \bar{p}$ and $M$. In particular, it follows that the domain $\Omega''$ and the constant $C(M)$ are both independent of $\A_i$ and $i$.  This proves statement {(a)} of Theorem \ref{Thm_compactness} - Global (i).

\vspace{.1cm} \noindent {\bf (b)}  By \eqref{initial_bound_2_proof}, estimate \eqref{uniform_bound_U_local} of Theorem \ref{Thm_existence} provides a uniform bound on each $U_i$ in $W^{1,2p}(\Omega'')$. The Banach Alaouglu compactness theorem \cite{Evans} then implies the existence of a subsequence of the $U_i$ which converges weakly in $W^{1,2p}(\Omega'')$ to some $U \in W^{1,2p}(\Omega'')$. Since weak convergence in $W^{1,2p}$ implies strong convergence in $L^{2p}$, which is a convergence regular enough to pass limits through product, it follows that $U$ lies pointwise in $SO(r,s)$ as well. 

\vspace{.1cm} \noindent {\bf (c)} To proof the main conclusion of Theorem \ref{Thm_compactness} - Global (i), note that by the uniform bound \eqref{uniform_bound_A_2_proof} on the $\A_{\textbf{b}_i}$ the Banach Alaouglu theorem yields the existence of a subsequence of the $\A_{\textbf{b}_i}$ converging weakly in $W^{1,\hat{p}}$ to some $\A_{\textbf{b}} \in W^{1,\hat{p}}(\Omega'')$, which again implies strong convergence in $L^{\hat{p}}$. Moreover, by restriction to those $i \in \mathbb{N}$ for which the $U_i$ converge to $U \in W^{1,2\bar{p}}$, it follows that the limit $\A_{\textbf{b}}$ agrees with the connection coefficients of $\A$ in the gauge $\textbf{b} = U\cdot \textbf{a}$, since strong $L^{\hat{p}}$ convergence is regular enough to pass limits through products. This completes the proof of Theorem \ref{Thm_compactness} - Global(i).  

\vspace{.1cm} \noindent {\bf Global (ii):} Assume a sequence of $L^\infty$ connections $(\A_i)_{i\in \mathbb{N}}$ satisfies \eqref{initial_bound_2} for $p=\infty$, $\bar{p}  \in (1,\infty)$, and assume that each $\A_i$ satisfies the smallness condition \eqref{initial_bound_small} for $p=\infty$. Theorem \ref{Thm_existence} then implies for each $\A_i \in L^\infty$, each $p<\infty$, and on every neighborhood $\Omega''$ compactly contained in $\Omega$, the existence of gauge transformations $U_i\in W^{1,2p}(\Omega'')$. We now choose transformations $U_i\in W^{1,2\bar{p}}(\Omega'')$, so $\hat{p}=\bar{p}$.  By the above proof of Theorem \ref{Thm_compactness} - Global (i), we now conclude that (a) - (c) of case Global (i) hold for $p=\bar{p}$ and $\hat{p}=\bar{p}$. This completes the proof of Theorem \ref{Thm_compactness} - Global (ii).

\vspace{.1cm} \noindent {\bf Local (ii):} Assume a sequence of $L^\infty$ connections $(\A_i)_{i\in \mathbb{N}}$ satisfies \eqref{initial_bound_2} for $p=\infty$ and $\bar{p} \in (1,\infty)$, but not the smallness condition \eqref{initial_bound_small}. Then, by Theorem \ref{Thm_opt}, there exists for any point $q \in \Omega$ a neighborhood $\Omega'' \subset \Omega$ of $q$, independent of $\A_i$, for which statement (a) above holds for $p=\bar{p}$ and $\hat{p}=\bar{p}$. Statements (b) and (c) then follow by an argument analogous to that in the proof of Theorem \ref{Thm_compactness} - Global(i). This completes the proof of Theorem \ref{Thm_compactness}. \hfill$\Box$

\appendix

\section{Basic results from elliptic PDE theory} \label{Sec_Prelimiaries-elliptic}

We now summarize the results from elliptic PDE theory used in this paper.  We assume throughout that $1< p<\infty$, and that $\Omega \subset \mathbb{R}^n$, $n\geq 2$, is a bounded open domain, simply connected and with smooth boundary. Our proofs in this paper use only the following two theorems of elliptic PDE theory, which directly extend to matrix valued and vector valued differential forms because the Laplacian acts component-wise. As above, we take the standard weak Laplacian to be the linear functional $\Delta u[\phi] = -\langle d u, d \phi\rangle_{L^2}$ on test functions $\phi \in W_0^{1,p^*}(\Omega)$ for fixed scalar functions $u \in W^{1,p}(\Omega)$, where $W^{1,p^*}_0(\Omega)$ is the closure of $C^\infty_0(\Omega)$ with respect to the $W^{1,p^*}$-norm (so $\phi|_{\partial\Omega}=0$).  Our first theorem is based on Theorem 7.2 in \cite{Simader}, but adapted to the case of solutions of the Poisson equation with non-zero Dirichlet data, (see Appendix B in \cite{ReintjesTemple_Uhl1} for a proof).    

\begin{Thm}\label{Thm_Poisson}
Let $\Omega\subset R^n$ be a bounded open set with smooth boundary $\partial\Omega$, assume $f\in W^{-1,p}(\Omega)$ and $u_0\in W^{1,p}(\Omega)\cap C^0(\overline{\Omega})$ for $1<p<\infty$.   Then the Dirichlet boundary value problem
\begin{eqnarray}
 \Delta u[\phi]&=&f[\phi],\ \ \text{in}\ \Omega  \label{Poisson-1}\\
u&=&u_0\ \  \text{on}\ \partial\Omega, \label{Poisson-2}
\end{eqnarray}
for any $\phi \in W_0^{1,p^*}(\Omega)$, has a unique weak solution $u \in W^{1,p}(\Omega)$ with boundary data $u-u_0 \in W_0^{1,p}(\Omega)$. Moreover, any weak solution $u$ of \eqref{Poisson-1} - \eqref{Poisson-2} satisfies
\begin{eqnarray}\label{Poisson-3}
\| u\|_{W^{1,p}(\Omega)}\leq C\left( \|f\|_{W^{-1,p}(\Omega)} + \|u_0\|_{W^{1,p}(\Omega)} \right),
\end{eqnarray}
for some constant $C$ depending only on $\Omega, n, p$.
\end{Thm}

We also require the following interior elliptic estimates, which for completeness we derived from \eqref{Poisson-3} in Appendix B in \cite{ReintjesTemple_Uhl1}. 

\begin{Thm}\label{Thm_Poisson_interior}
Let $f\in W^{m-1,p}(\Omega)$, for $m\geq 0$ and $1<p<\infty$. Assume $u$ is a weak solution of \eqref{Poisson-1}. Then $u \in W^{m+1,p}(\Omega')$ for any open set $\Omega'$ compactly contained in $\Omega$ and 
\begin{eqnarray}\label{Poisson-4}
\| u \|_{W^{m+1,p}(\Omega')} \leq C \big( \| f \|_{W^{m-1,p}(\Omega)}  + \|u\|_{W^{m,p}(\Omega)} \big),
\end{eqnarray}
where $C>0$ is a constant depending only on $\Omega, \Omega', m, n, p$.
\end{Thm}

In addition, our proofs require only the following three inequalities of Analysis. \\

\noindent \textbf{H\"older inequality:}
\beq \label{Hoelder}
|\langle u, v \rangle_{L^2} |\leq \|u\|_{L^p} \|v\|_{L^{p^*}},   
\hspace{.3cm}   u\in L^p(\Omega), \ \  v\in L^{p^*}(\Omega), \ \  \tfrac{1}{p} + \tfrac{1}{p^*} =1;
\eeq
\textbf{Morrey inequality:} 
\beq \label{Morrey_textbook}
\| u \|_{C^{0,\alpha}(\overline{\Omega})}  \leq C_0 \|u\|_{W^{1,p}(\Omega)},   
   \hspace{.3cm}   u \in W^{1,p}(\Omega), \ \ \alpha \equiv 1 - \frac{n}{p}, \ \ p>n;
\eeq
\textbf{Sobolev embedding inequality:}
\beq \label{Sobolev}
\| u \|_{L^{p'}(\Omega)} \leq C_S \| u \|_{W^{1,p}(\Omega)},       
\hspace{.3cm}    u \in W^{1,p}(\Omega), \ \  p' \equiv \frac{n p}{n-p}, \ \ 1<p<n;
\eeq
where $C_0>0$ and $C_S>0$ are constants depending only on $n, p$ and $\Omega$,  c. f. \cite{Evans}.

\section*{Declarations and Statements}

\begin{itemize}
\item Funding: M. Reintjes is currently supported by CityU Start-up Grant for New Faculty (7200748) and by CityU Strategic Research Grant (7005839); and was supported by the German Research Foundation, DFG grant FR822/10-1, from June 2019 until July 2021.
\item Authors declare that there is no conflict of interest and that there are no competing interests.
\item Availability of data and materials: Does not apply.
\item Authors' contributions: Both authors contributed equally.
\end{itemize}

\section*{Acknowledgements}
We thank Craig Evans for directing us to reference \cite{Simader} which contain the basic elliptic estimates which are the starting point for our analysis here.

\end{document}